\documentclass[11pt]{amsart}
\PassOptionsToPackage{dvipsnames}{xcolor}
\usepackage[colorlinks,citecolor=red,pagebackref,hypertexnames=false]{hyperref}

\usepackage{verbatim} 
\usepackage{placeins}

\makeatletter
\def\@tocline#1#2#3#4#5#6#7{\relax
	\ifnum #1>\c@tocdepth 
	\else
	\par \addpenalty\@secpenalty\addvspace{#2}%
	\begingroup \hyphenpenalty\@M
	\@ifempty{#4}{%
		\@tempdima\csname r@tocindent\number#1\endcsname\relax
	}{%
		\@tempdima#4\relax
	}%
	\parindent\z@ \leftskip#3\relax \advance\leftskip\@tempdima\relax
	\rightskip\@pnumwidth plus4em \parfillskip-\@pnumwidth
	#5\leavevmode\hskip-\@tempdima
	\ifcase #1
	\or\or \hskip 1em \or \hskip 2em \else \hskip 3em \fi%
	#6\nobreak\relax
	\dotfill\hbox to\@pnumwidth{\@tocpagenum{#7}}\par
	\nobreak
	\endgroup
	\fi}
\makeatother

\usepackage{tikz,amsthm,amsmath,amstext,amssymb,amscd,epsfig,euscript, mathrsfs,pspicture,multicol, mathtools,graphpap,graphics,graphicx,times,enumerate,subfig,sidecap,wrapfig,color}
\IfFileExists{dsfont.sty}{\usepackage{dsfont}}{\let\mathds\mathbb}

\usepackage{xcolor,booktabs, multirow} 
\usepackage{array}
\usepackage{pgfplots}
\pgfplotsset{compat=1.17}
\RequirePackage{caption}
\usepackage{enumitem}

\usepackage{ hyperref}
\usepackage{enumerate}

\numberwithin{equation}{section}

\usetikzlibrary{arrows.meta,calc,decorations.pathreplacing,decorations.markings,patterns,shapes.geometric,positioning,intersections,fadings,shadows}
\definecolor{domainc}{RGB}{226,236,248}
\definecolor{domaind}{RGB}{186,206,232}
\definecolor{solidA}{RGB}{158,166,178}
\definecolor{solidAtop}{RGB}{196,202,212}
\definecolor{solidAside}{RGB}{128,136,150}
\colorlet{solidB}{NavyBlue}
\colorlet{solidBtop}{NavyBlue!55!white}
\colorlet{solidBside}{NavyBlue!70!black}
\colorlet{solidC}{LimeGreen}
\colorlet{selc}{NavyBlue}
\definecolor{hlc}{RGB}{196,32,40}
\definecolor{flowc}{RGB}{34,74,150}
\colorlet{cutc}{NavyBlue!30!white}
\tikzset{
	lab/.style={font=\footnotesize},
	dimline/.style={<->,>={Stealth[length=4pt]},line width=.4pt},
	mapsto/.style={->,>={Stealth[length=6pt]},line width=.8pt},
	stream/.style={flowc,line width=.55pt},
	equip/.style={gray!55,line width=.35pt,dash pattern=on 1.5pt off 1.2pt},
}

\def\bT{{\mathbb{T}}}

\def\mR{{\mathcal{R}}}
\def\mS{{\mathcal{S}}}

\def\PP{{\mathcal{P}}}
\def\M{{\mathcal{M}}}

\def\cA{{\mathscr{A}}}

\def\one{\mathds{1}}

\def\ve{\varepsilon}

\def\Rn1{\mathbb{R}^{n+1}}

\def\dist{\textup{dist}} 						
\def\ker{\mathop\mathrm{ker }}						
\def\supp{\operatorname{supp}}					

\DeclareMathOperator*{\esssup}{ess\,sup}			
\renewcommand{\div}{\operatorname{div}}			

\def\Xint#1{\mathchoice
	{\XXint\displaystyle\textstyle{#1}}%
	{\XXint\textstyle\scriptstyle{#1}}%
	{\XXint\scriptstyle\scriptscriptstyle{#1}}%
	{\XXint\scriptscriptstyle\scriptscriptstyle{#1}}%
	\!\int}
\def\XXint#1#2#3{{\setbox0=\hbox{$#1{#2#3}{\int}$ }
		\vcenter{\hbox{$#2#3$ }}\kern-.58\wd0}}

\def\avint{\Xint-}

\theoremstyle{plain}
\newtheorem{theorem}{Theorem}

\newtheorem{corollary}[theorem]{Corollary}

\newtheorem{lemma}[theorem]{Lemma}

\newtheorem{proposition}[theorem]{Proposition}

\theoremstyle{definition}

\newtheorem{definition}[theorem]{Definition}
\newtheorem{remark}[theorem]{Remark}

\numberwithin{equation}{section}
\numberwithin{theorem}{section}

\IfFileExists{mathb10.tfm}{
	\DeclareFontFamily{U}{mathb}{\hyphenchar\font45} 
	\DeclareFontShape{U}{mathb}{m}{n}{
		<5> <6> <7> <8> <9> <10> gen * mathb
		<10.95> mathb10 <12> <14.4> <17.28> <20.74> <24.88> mathb12
	}{}
	\DeclareSymbolFont{mathb}{U}{mathb}{m}{n}

	\DeclareMathSymbol{\toitself}      {3}{mathb}{"FD}  

}{\providecommand{\toitself}{\circlearrowleft}}

\def\R{\mathbb{R}}

\DeclareMathOperator{\Tr}{Tr}

\usepackage{float}
\usetikzlibrary{fit}
\def\sym{\mathop\mathrm{sym}}
\def\mA{{A^\sharp}}
\def\mC{{\mathcal{C}}}
\def\mE{{\mathcal{E}}}

\definecolor{mathBlue}{RGB}{34,79,119}
\definecolor{mathOrange}{RGB}{190,98,28}
\definecolor{mathGreen}{RGB}{35,125,100}
\definecolor{lightFill}{RGB}{239,244,248}

\begin{document}

	\title[Neumann problems and single layers]{On the $L^p$-Neumann problem for divergence form elliptic equations in the upper half-space}

	\author[M. Mourgoglou]{Mihalis Mourgoglou}
	\address{Mihalis Mourgoglou\\ Departamento de Matem\'aticas, Universidad del Pa\'is Vasco (UPV/EHU), Barrio Sarriena s/n 48940 Leioa, Spain and
		Ikerbasque, Basque Foundation for Science, Bilbao, Spain.}
	\email{michail.mourgoglou@ehu.eus}

	\thanks{The author was supported by IKERBASQUE and partially supported by the project PID2024-157724NB-I00, funded by MICIU/AEI/10.13039/501100011033/FEDER, EU, by IT-1615-22, and IT1875-26 (Basque Government).}
	\date{}
	\hypersetup{pdfauthor={Mihalis Mourgoglou},pdftitle={On the Lp-Neumann problem for divergence form elliptic equations in the upper half-space}}

	\begin{abstract}
	A central question has been whether the $L^p$-Neumann problem is solvable
	for some $p\in(1,\infty)$ in the upper half-space
	$\R^{n+1}_+\coloneqq\R^n\times(0,\infty)$, $n\geq2$, for second-order
	uniformly elliptic equations in divergence form whose coefficients are either
	independent of the transverse variable or satisfy the Dahlberg--Kenig--Pipher
	(DKP) condition. We answer this question in the negative for both classes,
	even for matrices whose symmetric part is the identity. The two examples have
	the same nonzero Neumann datum $g\in C_c^\infty(\R^n)$ of mean zero. Although
	$g$ is bounded and belongs to the real Hardy space $H^1(\R^n)$, neither
	equation admits a weak solution with this datum whose gradient has a
	non-tangential maximal function in $L^p(\R^n)$ for any $1\leq p<\infty$.
	For the variational solution $u$ associated with either matrix $A$, the
	non-tangential maximal functions of $\nabla_xu$, $\partial_tu$ and
	$(A\nabla u)_{n+1}$ each have infinite $L^p$ norm for every $1\leq p\leq\infty$.
	
	In the example with coefficients independent of the transverse variable,
	the equation in the interior is the Laplace equation. The matrix in the DKP
	example has finite $L^2$-type Carleson discrepancy from the identity. This
	shows that solvability of the $L^p$-Neumann problem is not stable under
	arbitrary Carleson perturbations of the identity. For a suitable choice of
	the coefficients in the lower half-space, the $L^2$ non-tangential maximal
	and uniform slice estimates fail for both the tangential gradient and the
	transverse derivative of the associated single layer potential. The $L^2$
	estimate for the tangential gradient of its boundary trace also fails. We explain the role of the
	coefficients below the boundary and prove the $L^2$ layer potential estimates
	for the reflected extension when the DKP norm is sufficiently small.
	\end{abstract}

	\maketitle

	\setcounter{tocdepth}{1}
	\tableofcontents

	\section{Introduction}\label{sec:introduction}\label{hs:rellich}
	
	Let $L=-\div(A\nabla)$ be a second-order elliptic operator in divergence form in the upper half-space $\R^{n+1}_+$, whose coefficients are real, bounded and uniformly elliptic. In the Dirichlet and regularity problems one prescribes the boundary trace of a solution $u$ of $Lu=0$, and in the Neumann problem one prescribes its outward conormal derivative. In the Dirichlet problem one asks for an $L^p$ bound for the non-tangential maximal function of $u$. In the regularity and Neumann problems one asks for the corresponding bound for $\nabla u$, in terms of the $L^p$ norm of the tangential gradient of the boundary datum and of the Neumann datum, respectively. For general bounded measurable coefficients even the Dirichlet estimate can fail. Indeed, L. Caffarelli, E. Fabes and C. Kenig \cite{CFK81} and L. Modica and S. Mortola \cite{MM81} constructed operators whose elliptic measure is singular with respect to the surface measure, so that the Dirichlet problem fails for every finite exponent. Two classical hypotheses on the dependence of the coefficients on the transverse variable are
	\begin{enumerate}[label=\textup{(\roman*)}]
		\item $A$ is independent of the transverse variable $t$;
		\item $A$ satisfies the Dahlberg--Kenig--Pipher (DKP) condition defined in \eqref{eq:dkp-definition}.
	\end{enumerate}
	In the upper half-space, the Dirichlet problem is solvable in $L^q$ for some $q<\infty$ in both cases, see \cite{HKMP15} and \cite{KP01}, respectively. The regularity problem is solvable in $L^p$ for some $p>1$ for $t$-independent coefficients \cite{HKMP15b}. For DKP coefficients, Dindo\v{s}, Hofmann and Pipher proved this for domains above Lipschitz graphs \cite[Theorem~1.1]{DHP23}. Independently, in joint work with B.~Poggi and X.~Tolsa \cite[Corollary~1.35]{MPT25}, we proved $L^p$-regularity for some $p>1$ in bounded corkscrew domains in $\R^{n+1}$, $n\geq2$, with $n$-Ahlfors regular boundary and interior big pieces of chord-arc domains, under the DKP condition.

	For general rough coefficients the gradient estimates are formulated using non-tangential maximal functions of local $L^2$ averages. In the plane the $L^p$-Neumann problem is solvable for some $p>1$ in both cases \cite{KR09,DHP23}. In higher dimensions it had remained open whether either hypothesis suffices for $L^p$-Neumann solvability for some $1<p<\infty$.
	
	In this paper we show that \textit{the answer is negative for both classes}, for every $n\geq 2$. Our main result, Theorem~\ref{thm:main}, provides two real uniformly elliptic matrices in $\R^{n+1}_+$, $n\geq2$, one in each of the classes (i) and (ii), and a single nonzero datum $g\in C_c^\infty(\R^n)$ with zero mean, for which the $L^p$-Neumann problem fails for every $1<p<\infty$. Both matrices have symmetric part equal to the identity. The first is independent of $t$. Its skew-symmetric part contributes zero to $\div(A\nabla u)$ in the interior, so the equation is still the Laplace equation. It does, however, change the conormal derivative. The second matrix satisfies the DKP condition and has finite Carleson discrepancy from the identity. Thus solvability of the Neumann problem need not be preserved under large Carleson perturbations, see Corollary~\ref{cor:large-perturbations}. This has consequences for single layer potentials, whose estimates fail for a suitable extension of the coefficients to the lower half-space (Theorem~\ref{sl:negative}), and hold for the reflected extension if the DKP norm is small (Theorem~\ref{sl:small-layer}). Our construction is inspired by the geometric counterexample in joint work with C. Puliatti and X. Tolsa \cite{Geo}, where we proved that the $L^p$-Neumann estimate for the variational solution of the Laplace equation can fail in a bi-Lipschitz image of a ball. In that paper we modified the boundary of the domain. Here we keep the half-space fixed and modify the coefficients. The main ideas are described in Subsection~\ref{sec:intro-neumann-proof}.
	
	\subsection{The main results}\label{sec:main-results}
	Throughout the paper $n\geq2$. Let $A$ be a real $(n+1)\times(n+1)$ matrix with bounded measurable entries in $\R^{n+1}_+$. We assume that for some $0<\lambda\leq\Lambda<\infty$ and for almost every $(x,t)\in\R^{n+1}_+$,
	\begin{equation}\label{eq:uniform-ellipticity}
		A(x,t)\xi\cdot\xi\geq\lambda|\xi|^2,
		\qquad |A(x,t)\xi|\leq\Lambda|\xi|
		\qquad\text{for every }\xi\in\R^{n+1}.
	\end{equation}
	We set $L=L_A\coloneqq-\div(A\nabla)$.
	We say that a function $u \in Y^{1,2}_{\mathrm{loc}}(\R^{n+1}_+)$ (see Section~\ref{sec:sobolev})  is a {\it weak solution} of $Lu=0$ in $\R^{n+1}_+$ if
	\begin{equation}
		\int_{\R^{n+1}_+} A\nabla u\cdot\nabla\varphi\,dx\,dt=0 \quad\text{for every } \varphi \in C^\infty_c(\R^{n+1}_+).
	\end{equation}
	
	For $a>1$ and $x\in\R^n$, let
	\[
	\gamma_a(x)=\{(y,t)\in\R^{n+1}_+: |(y,t)-(x,0)|<at\}.
	\]
	The non-tangential maximal function of a measurable function $f$ is
	\[
	N_a(f)(x)=\esssup_{(y,t)\in\gamma_a(x)}|f(y,t)|,
	\qquad N=N_2.
	\]
	We also use the averaged non-tangential maximal function
	\[
	\widetilde N_a(f)(x)=\sup_{(y,t)\in\gamma_a(x)}
	\left(\avint_{B((y,t),t/2)}|f|^2\right)^{1/2}.
	\]
	
	\subsubsection*{The Neumann problem}
	{
		Given $g\in L^1_{\mathrm{loc}}(\R^n)$, we say that a function
		$u\in Y^{1,2}_{\mathrm{loc}}(\R^{n+1}_+)$, whose gradient is locally
		integrable up to the boundary, is a {\it weak solution of the Neumann
			problem} with datum $g$ if
		\begin{equation}\label{eq:weak-neumann}
			\int_{\R^{n+1}_+}A\nabla u\cdot\nabla\varphi\,dx\,dt
			=\int_{\R^n}g(x)\varphi(x,0)\,dx
			\qquad\text{for all }\varphi\in C_c^\infty(\R^{n+1}).
		\end{equation}
	}
	\begin{definition}\label{def:neumann-solvability}
		{
			Let $p\in (1,\infty)$. We say that the {\it $L^p$-Neumann problem} is solvable if there exists $C_p<\infty$ such that, for every	$g\in L^p(\R^n)$, there exists a weak solution $u$ of \eqref{eq:weak-neumann} satisfying
			\begin{equation}\label{eq:Np}
				\|N(\nabla u)\|_{L^p(\R^n)}	\leq C_p\|g\|_{L^p(\R^n)}.
			\end{equation}
		}
	\end{definition}
	{
		The maximal function bound implies local integrability of the gradient
		up to the boundary, so the equality \eqref{eq:weak-neumann} is meaningful.
		We say that the problem is {\it well posed} if, in addition, the solution
		is unique up to an additive constant among weak solutions with
		$N(\nabla u)\in L^p(\R^n)$.
	}
	
	We first define the two Carleson measures used in the statement. For a positive measure $\mu$ on $\R^{n+1}_+$, we set
	\[
	\|\mu\|_{\mC}
	=\sup_{x\in\R^n,\,r>0}r^{-n}
	\mu(B(x,r)\times(0,r)).
	\]
	If $A$ is locally Lipschitz in $\R^{n+1}_+$, let $\nabla A$ denote its full gradient. We set
	\begin{equation}\label{eq:dkp-definition}
		\,d\mu_A(x,t)
		=t\esssup_{(y,s)\in B((x,t),t/2)}|\nabla A(y,s)|^2\,dx\,dt,
	\end{equation}
	and say that $A$ satisfies the \textit{Dahlberg--Kenig--Pipher condition} if
	\[
	\|\mu_A\|_{\mC}<\infty.
	\]
	
	For the discrepancy from the identity, we set
	\begin{equation}\label{eq:discrepancy-definition}
		a_A(x,t)=\esssup_{(y,s)\in B((x,t),t/2)}|A(y,s)-I_{n+1}|,
		\qquad \,d\nu_A(x,t)=a_A(x,t)^2\frac{\,dx\,dt}{t}.
	\end{equation}
	
	If $A$ satisfies the DKP condition, then
	$t\|\nabla A\|_{L^\infty(B((y,t),t/2))}\lesssim_n\|\mu_A\|_{\mC}^{1/2}$.
	Thus, for every weak solution of $L_Au=0$, the interior gradient estimate
	\cite[Theorem~8.17 and Corollary~8.36]{GT} gives
	\[
	\|\nabla u\|_{L^\infty(B((y,t),t/4))}\lesssim
	\left(\avint_{B((y,t),t/2)}|\nabla u|^2\right)^{1/2},
	\]
	with a constant depending only on $n,\lambda,\Lambda$ and $\|\mu_A\|_{\mC}$.
	Since $B((y,t),t/2)\subset\gamma_{2a+1}(x)$ whenever $(y,t)\in\gamma_a(x)$, we have
	\begin{equation}\label{eq:gradient-maximal-comparison}
		N_a(\nabla u)\lesssim\widetilde N_a(\nabla u)\leq N_{2a+1}(\nabla u).
	\end{equation}
	By change of aperture (see \cite[proof of Lemma~1.27]{HT20}), for every $a>1$ and $0<p<\infty$ we have
	$\|N_a(\nabla u)\|_p\approx\|\widetilde N_a(\nabla u)\|_p$.
	The same holds for $p=\infty$, since $\|N_a(\nabla u)\|_{L^\infty(\R^n)}=\|\nabla u\|_{L^\infty(\R^{n+1}_+)}$
	for every $a>1$.
	
	For a bounded real vector field $c:\R^{n+1}_+\to\R^n$, we set
	\begin{equation}\label{eq:matrix}
		S(c)\coloneqq\begin{pmatrix}0&c\\-c^T&0\end{pmatrix}
		\qquad \text{ and }\qquad A_c\coloneqq I_{n+1}+S(c).
	\end{equation}
	We write $\sym A\coloneqq(A+A^T)/2$ for the symmetric part of $A$.
	By \eqref{eq:matrix},
	\[
	(A_c\nabla u)_{n+1}=\partial_tu-c\cdot\nabla_xu.
	\]
	The outward normal is $-e_{n+1}$, so the outward conormal derivative is
	the weak boundary trace of $-(A\nabla u)_{n+1}$.
	
	Our main result is the following.
	\begin{theorem}\label{thm:main}
		Let $n\geq2$. There exist real uniformly elliptic matrices $A^{(1)}$, $A^{(2)}$, with $\sym A^{(i)}=I_{n+1}$, and a single $0\ne g\in C_c^\infty(\R^n)$ with $\int g=0$, such that:
		\begin{enumerate}[label=\textup{(\roman*)}]
			\item $A^{(1)}$ is independent of $t$ and $L_{A^{(1)}}=-\Delta$;
			\item $A^{(2)}$ satisfies the DKP condition and $\|\nu_{A^{(2)}}\|_{\mC}<\infty$.
		\end{enumerate}
		For $i=1,2$, let $u^{(i)}$ be the variational solution of \eqref{eq:variational} with matrix $A^{(i)}$ and datum $g$. Then for every $a>1$ and $1\leq p\leq\infty$,
		\begin{equation}\label{eq:main-components}
			\|N_a(\nabla_xu^{(i)})\|_p
			=\|N_a(\partial_tu^{(i)})\|_p
			=\|N_a((A^{(i)}\nabla u^{(i)})_{n+1})\|_p=\infty.
		\end{equation}
		For either coefficient, no weak solution $v$ with datum $g$ satisfies $N_a(\nabla v)\in L^p(\R^n)$ for any $a>1$ and $1\leq p<\infty$.
	\end{theorem}
	
	At $p=1$, the datum $g$ is a multiple of a Hardy space atom, as we verify after \eqref{eq:fixed-datum}. For this datum no weak solution $v$ satisfies
	\[
	\|N(\nabla v)\|_{L^1(\R^n)}\leq C\|g\|_{H^1(\R^n)}.
	\]
	At $p=\infty$, $g$ is bounded and the three maximal functions in \eqref{eq:main-components} have infinite $L^\infty$ norm for the variational solution. The comparison with an arbitrary weak solution uses $p<\infty$.
	
	The matrices and the datum are independent of $p$. Thus Theorem~\ref{thm:main} answers the question at the beginning of this section for both classes. For comparison, the Dirichlet and regularity problems for the same two matrices have the following properties.
	\begin{itemize}
		\item For $A^{(1)}$, the equation itself is the Laplace equation, so that the Dirichlet and regularity problems for $L_{A^{(1)}}$ are those for the Laplacian, which are solvable for every $1<p<\infty$.
		\item For $A^{(2)}$, the Dirichlet and regularity problems in the upper half-space are solvable for some exponents by \cite{KP01} and \cite[Theorem~1.1]{DHP23}, respectively.
	\end{itemize}
	The assumption $n\geq2$ cannot be removed. Indeed, in the plane the $L^p$-Neumann problem is solvable for some $p>1$ both for real $t$-independent coefficients \cite{KR09} and under the DKP condition \cite{DHP23}, see Subsection~\ref{sec:intro-neumann-history}.
	
	\begin{remark}\label{rem:main}
		(a) The Neumann datum prescribes the outward conormal derivative, while \eqref{eq:Np} asks for control of the full gradient. Theorem~\ref{thm:main} shows that even the maximal functions of $\partial_tu$ and $(A\nabla u)_{n+1}$ have infinite $L^p$ norm for the variational solution. This holds for fixed coefficients and a fixed datum.
		
		(b) In the weak sense of \eqref{eq:weak-neumann}, the boundary value of $(A\nabla u)_{n+1}$ is $-g$. In the open half-space the equation gives
		\[
			\partial_t\bigl((A\nabla u)_{n+1}\bigr)=-\div_x\bigl((A\nabla u)_x\bigr),
		\]
		in the sense of distributions. For the variational solution with $A=I_{n+1}$ we have $\partial_tu(\cdot,t)=-P_tg$, so the $L^p$ bound for $N(\partial_tu)$, $1<p<\infty$, follows from the properties of the Poisson kernel. In our DKP example, $A=I_{n+1}$ near every nonzero boundary point. Thus both $(A\nabla u)_{n+1}$ and $\partial_tu$ have boundary value $-g$ almost everywhere, although their non-tangential maximal functions have infinite $L^p$ norm for every $1\leq p\leq\infty$.
		
		(c) In both examples the one-sided Rellich estimate $\|\nabla_x\Tr u\|_p\lesssim\|g\|_p$ fails, even though $A\nabla u\cdot\nabla u=|\nabla u|^2$. Indeed, summing \eqref{eq:main-selected} over $j$ gives $\int_{E_i}|\nabla_x\Tr u|\,dx=\infty$, and $|E_i|<\infty$.
		
	\end{remark}
	
	The following corollary of Theorem~\ref{thm:main} concerns Carleson perturbations of the identity.
	\begin{corollary}\label{cor:large-perturbations}
		{
			There exists a real uniformly elliptic matrix $A$ in $\R^{n+1}_+$,
			with $\sym A=I_{n+1}$, satisfying the DKP condition and
			$\|\nu_A\|_{\mC}<\infty$, for which the $L^p$-Neumann problem is
			not solvable for any $1<p<\infty$.
			For the identity matrix, the $L^p$-Neumann problem is solvable
			for every $1<p<\infty$.
		}
		
	\end{corollary}
	\begin{proof}
		For the identity matrix, the estimate follows from the Poisson representation and the $L^p$ boundedness of the Riesz transforms, see Remark~\ref{rem:laplace-neumann}. On the other hand, it suffices to take $A=A^{(2)}$ from Theorem \ref{thm:main}.
	\end{proof}
	
	For the Dirichlet problem, solvability in $L^q$ for some $q<\infty$ is preserved under Carleson perturbations of arbitrary size, possibly with a change of exponent. See \cite{FKP91} and \cite[Theorem~1.2]{DFM23}. The same holds for the regularity problem in $L^p$ for some $p>1$, with the maximal function of local $L^2$ averages of the gradient \cite[Corollary~1.6 and Remark~1.7]{DFM23}. On the other hand, the known perturbation results for the Neumann problem require the Carleson norm of the discrepancy to be small \cite{KP95,AA11}. Up to the choice of Whitney regions, $\|\nu_A\|_{\mC}<\infty$ is equivalent to finiteness of the modified Carleson norm of $A-I_{n+1}$ in \cite{AA11}. Corollary~\ref{cor:large-perturbations} shows that the smallness assumption cannot be removed from the Neumann estimate, even for perturbations of the identity. In \cite{AA11} the maximal function is defined using Whitney $L^2$ averages of the gradient. The nonexistence statement of Theorem~\ref{thm:main} also holds with $\widetilde N_a$. For $A^{(2)}$ this follows from \eqref{eq:gradient-maximal-comparison}, and for $A^{(1)}$ it follows from the mean value property for the harmonic components of the gradient.
	
	In bounded one-sided chord-arc domains, J. Feneuil and L. Li
	\cite[Propositions~1.17 and~1.18]{FL24} proved that, for $1<p<\infty$ and $p'=p/(p-1)$,
	solvability of the $L^{p'}$ Poisson--Neumann problem for $L^*$ implies solvability of the
	$L^p$-Neumann problem for $L$. The converse holds if the $L^{p'}$-Dirichlet problem for
	$L^*$ is solvable. In \cite[Theorem~1.14]{FL26} they proved that $L^p$ Poisson--Neumann
	solvability in these domains is preserved under sufficiently small Carleson perturbations
	of the coefficients, provided that the elliptic measure of the unperturbed operator belongs
	to $A_\infty$ with respect to surface measure.

	\subsubsection*{Layer potentials}
	Let $\cA$ be a real bounded uniformly elliptic matrix on $\R^{n+1}$, with ellipticity constant $\lambda$ and operator norm bounded by $\Lambda$. For $f\in\dot H^{-1/2}(\R^n)$, we define the \textit{single layer potential} $\mS_{\cA}f$ to be the unique function in $Y^{1,2}(\R^{n+1})$ such that
	\begin{equation}\label{sl:variational}
		\int_{\R^{n+1}}\cA\nabla\mS_{\cA}f\cdot\nabla\varphi\,dx\,dt
		=\langle f,\Tr\varphi\rangle,
		\qquad\text{ for every }\varphi\in Y^{1,2}(\R^{n+1}).
	\end{equation}
	For the definitions of the Sobolev and trace spaces we refer to Section \ref{sec:sobolev}. The trace estimate and the Lax--Milgram theorem give existence and uniqueness, and testing with $\varphi=\mS_{\cA}f$ yields
	\begin{equation}\label{sl:variational-bound}
		\|\nabla\mS_{\cA}f\|_{L^2(\R^{n+1})}
		\leq C_{n,\lambda}\|f\|_{\dot H^{-1/2}}.
	\end{equation}
	Note that test functions supported in either open half-space have zero trace, so $L_{\cA}\mS_{\cA}f=0$ there and, since $\cA$ is real, the De Giorgi--Nash estimates hold for $L_{\cA}$ and its adjoint. By \cite[Theorems~3.1 and~3.2]{HK07}, there exists a fundamental solution $\Gamma_{\cA}$, continuous off the diagonal, satisfying
	\[
	|\Gamma_{\cA}(z,w)|\leq C_{n,\lambda,\Lambda}|z-w|^{1-n}.
	\]
	For $f\in C_c^\infty(\R^n)$, we have the representation formula
	\begin{equation}\label{sl:kernel}
		\mS_{\cA}f(x,t)
		=\int_{\R^n}\Gamma_{\cA}((x,t),(y,0))f(y)\,dy,
		\qquad\text{ for }t\ne0,
	\end{equation}
	which is proved in Lemma~\ref{sl:kernel-variational} for the continuous representative in each open half-space.
	
	Following \cite[p.~450]{AM14}, we extend $A$ to the lower half-space by changing the signs of the mixed entries under reflection. In our coordinates $(x,t)$, the resulting matrix is
	\begin{equation}\label{sl:reflection-matrix}
		\mathcal{J}=\operatorname{diag}(I_n,-1),\qquad
		\mA(x,t)=
		\begin{cases}
			A(x,t),&t>0,\\
			\mathcal{J} A(x,-t)\mathcal{J},&t<0.
		\end{cases}
	\end{equation}
	We write $N_a^+=N_a$ for the maximal function in the upper half-space and $N_a^-$ for the corresponding maximal function in the lower half-space, with cones obtained from $\gamma_a$ by reflection. The superscripts in $\widetilde N_a^\pm$ have the same meaning. In the lower half-space the averaging ball is $B((y,t),|t|/2)$.
	
	If $u$ is the variational Neumann solution for $A$ with datum $g$, then
	\begin{equation}\label{sl:intro-even-identity}
		\mS_{\mA}(2g)(x,t)=u(x,|t|),\qquad t\ne0,
	\end{equation}
	as we prove in Lemma~\ref{sl:reflection}. Indeed, the even extension has outward conormal derivative $g$ on each side, so the density in \eqref{sl:variational} is $2g$.
	Theorem~\ref{thm:main} and \eqref{sl:intro-even-identity} give the following result.
	\begin{theorem}\label{sl:negative}\label{hs:bad-layers}
		Let $A=A^{(i)}$, $i=1,2$, and $g$ be as in Theorem~\ref{thm:main}, and let $\mA$ be its reflected extension \eqref{sl:reflection-matrix}. For $u=\mS_{\mA}(2g)$, every $a>1$ and $1\leq p\leq\infty$,
		\begin{equation}\label{sl:negative-maximal}
			\|N_a^\pm(\nabla_xu)\|_p
			=\|N_a^\pm(\partial_tu)\|_p
			=\|N_a^\pm((\mA\nabla u)_{n+1})\|_p=\infty.
		\end{equation}
	\end{theorem}
	
	The density $2g$ belongs to $H^1(\R^n)\cap L^\infty(\R^n)$. Thus the $H^1 \to L^1$ and the $L^\infty \to L^\infty$ maximal estimates for the single layer both fail in each half-space.
	
	In each half-space, the reflected matrices still belong to the corresponding class in Theorem~\ref{thm:main}: for $i=1$ they are independent of $t$, and for $i=2$ they satisfy the DKP condition and have finite Carleson discrepancy from the identity. The three equalities in \eqref{sl:negative-maximal} also hold with $\widetilde N_a^\pm$. The slice and tangential trace estimates fail as well, see Subsection~\ref{lp:failure-section}. Moreover, for every $0<\varepsilon<1/2$, the $\varepsilon$-truncation of the tangential boundary kernel applied to $2g$ has infinite $L^2$ norm on a set where the potential is smooth, see Corollary~\ref{sl:boundary-truncations}.
	
	The choice of the lower coefficients matters. For $A^{(1)}=I+S(b)$, one may instead keep the same coefficients on both sides of the boundary, that is,
	\begin{equation}\label{eq:unchanged_extension}
		A_b(x,t)\coloneqq I+S(b(x)),\qquad\text{ for }t\in\R.
	\end{equation}
	Since $L_{A_b}=-\Delta$ in the whole space, $\mS_{A_b}$ is just the classical single layer potential for the Laplacian, whose gradient satisfies the $L^2$ maximal, slice and tangential trace estimates. However, its upper conormal operator is not onto $L^2$, see Remark~\ref{rem:two-extensions}. So the layer estimates hold for \eqref{eq:unchanged_extension} and fail for \eqref{sl:reflection-matrix}, although both extensions coincide in the upper half-space, and hence give the same Neumann problem there.
	
	For the reflected extension \eqref{sl:reflection-matrix}, the Neumann estimate in \cite{DPR17,DP24} gives the following result when the DKP norm is sufficiently small.
	
	\begin{theorem}\label{sl:small-layer}\label{hs:small-layers}
		Let $A$ be a real locally Lipschitz matrix in $\R^{n+1}_+$, with ellipticity constant $\lambda$ and operator norm bounded by $\Lambda$. There exists $\ve_0=\ve_0(n,\lambda,\Lambda)>0$ such that, if $\|\mu_A\|_{\mC}<\ve_0$, its reflected extension $\mA$ satisfies, for every $f\in L^2(\R^n)\cap\dot H^{-1/2}(\R^n)$ and $a>1$,
		\begin{equation}\label{sl:small-maximal}
			\|N_a^+(\nabla\mS_{\mA}f)\|_2
			+\|N_a^-(\nabla\mS_{\mA}f)\|_2
			\leq C_{n,\lambda,\Lambda,a}\|f\|_2
		\end{equation}
		and
		\begin{equation}\label{sl:small-slice-trace}
			\sup_{t\ne0}\|\nabla\mS_{\mA}f(\cdot,t)\|_2
			+\|\nabla_x\Tr\mS_{\mA}f\|_2
			\leq C_{n,\lambda,\Lambda}\|f\|_2.
		\end{equation}
		The maximal estimate also holds with $\widetilde N_a^\pm$.
	\end{theorem}
	
	By completion, these estimates define the layer gradient for every $f\in L^2(\R^n)$. In Subsection~\ref{sl:small-section} we prove the estimates and identify the tangential derivative of the trace and the two outward conormal derivatives. In Subsection~\ref{sec:intro-layer-history} we compare Theorem~\ref{sl:small-layer} with earlier results on layer potentials.
	
	\subsection{History of the Neumann problem}\label{sec:intro-neumann-history}
	For the Laplacian in the upper half-space, the Poisson extension and the $L^p$ boundedness of the Riesz transforms give Neumann solvability for every $1<p<\infty$ as discussed in the classical literature, see for instance \cite{Stein}. On bounded Lipschitz domains, the corresponding theory gives solvability for $1<p<2+\varepsilon$ for Neumann data of mean zero \cite{JK81,DK87}; see also \cite{Ver84,Ke94}. Flattening a Lipschitz graph transforms the Laplacian into a divergence form operator in the upper half-space with real symmetric $t$-independent coefficients, and transforms the normal derivative into a conormal derivative. The resulting matrices form a particular subclass of the symmetric matrices; see \cite{AAM10}. The general real symmetric $t$-independent case also has Neumann solvability for $1<p<2+\varepsilon$ \cite{KP93}. For merely measurable coefficients, the gradient estimates use local $L^2$ averages in the non-tangential maximal function.
	
	Beyond Lipschitz domains, much less is known about the Neumann problem for the Laplacian; see, for instance, \cite{HMT,MT24b} for positive results. In a bounded corkscrew domain $\Omega\subset\R^{n+1}$, $n\geq2$, with Ahlfors regular boundary, solvability of the $L^p$-Neumann problem for the Laplacian by $W^{1,2}(\Omega)$ solutions, for some $1<p<\infty$, implies that the boundary is uniformly rectifiable \cite[Theorem~1.8]{MMP}. On the other hand, the example in \cite{Geo} is a bi-Lipschitz image of a ball where the $L^p$-Neumann estimate for the variational solution fails for every $1<p<\infty$, although the Dirichlet and regularity problems are solvable there for some exponents. Pulling the Laplacian back to the unit ball gives a real symmetric matrix with bounded measurable coefficients for which the variational estimate fails in the same range. The present paper shows that failure can also occur within the two classes of coefficients considered here, for which the Dirichlet and regularity problems are solvable.
	
	For real $t$-independent matrices, the Dirichlet problem is solvable for some finite exponent, without symmetry or horizontal regularity, since the associated elliptic measure belongs to $A_\infty$ \cite{HKMP15}; in the plane, this goes back to \cite{KKPT00}. The regularity problem is also solvable in $L^p$ for some $p>1$ in this generality \cite{HKMP15b}. The $L^2$-Neumann problem is well posed for sufficiently small complex $t$-independent $L^\infty$ perturbations of real symmetric, block, or constant $t$-independent matrices. More generally, $L^2$ well-posedness persists under sufficiently small $t$-independent $L^\infty$ perturbations \cite{AAH08,AAM10}. In the first order approach, $L^2$ well-posedness of the Neumann problem requires that the map taking the boundary conormal gradient to its normal component be invertible. Related formulations in terms of one-sided Rellich inequalities allow Neumann solvability to be extrapolated to smaller exponents \cite{AM14}, and the representation and uniqueness of solutions in this framework are studied in \cite{AM}. Thus, control of the boundary gradient from the conormal derivative remains a separate issue.
	
	Carleson conditions enter through adapted changes of variables, which flatten a Lipschitz boundary while controlling derivatives of the coefficients, and through perturbation theory, which measures the discrepancy between two matrices \cite{FKP91,KP95,KP01}. The stability of the Dirichlet, regularity and Neumann problems under Carleson perturbations is discussed after Corollary~\ref{cor:large-perturbations}. A finite derivative Carleson condition gives Dirichlet solvability for some finite exponent, also without symmetry \cite{KP01}. For a DKP operator $L$ on a bounded corkscrew domain $\Omega\subset\R^{n+1}$, $n\geq2$, with uniformly $n$-rectifiable boundary, we proved in \cite[Theorem~1.33]{MPT25} that $L^{p'}$-Dirichlet solvability for $L^*$ implies $L^p$-regularity for $L$, where $1<p<\infty$ and $p'=p/(p-1)$. In this estimate the boundary datum is measured in the homogeneous Haj\l asz--Sobolev seminorm.

	For the Neumann problem, solvability at every fixed $1<p<\infty$ follows when the oscillation Carleson norm and the Lipschitz constant of the boundary are sufficiently small \cite{DPR17}, see also \cite[Theorem~2.2]{DP24}. The smallness depends on $p$, and the boundary Lipschitz constant is zero in the upper half-space. Dindo\v{s} and Pipher \cite[Section~6]{DP24} observed that, in dimensions $n+1\geq3$, it was open even in smooth domains whether a finite Carleson norm suffices for $L^p$-Neumann solvability for some $p>1$. For a parabolic result under a smallness assumption, see \cite{DLP26}.
	
	Dimension two admits stronger conclusions. Real $t$-independent matrices give Neumann solvability for some $p>1$ without symmetry \cite{KR09}, and a finite oscillation Carleson norm also gives solvability for $p>1$ sufficiently close to one \cite{DHP23}. Conjugate solutions relate the Neumann and regularity problems in these arguments. For $t$-independent coefficients, failure at any prescribed exponent can occur as the matrix varies \cite{KR09}, and the analogous examples for the Dirichlet problem are given in \cite{KKPT00}. These examples belong to the family $I_2+S(k\operatorname{sgn}(x))$, $k\in\R$, in the notation of \eqref{eq:matrix}, so that their skew-symmetric part has a jump across $\{x=0\}$; see also \cite{Axe10}. In the plane, the condition $\div_xc=0$ forces $c$ to depend only on $t$, so that the localized divergence-free perturbations used in this paper are not available; see Subsection~\ref{sec:intro-neumann-proof}.
	
	\subsection{History of layer potentials}\label{sec:intro-layer-history}
	The method of layer potentials reduces a boundary value problem to an equation for a density on the boundary. Its use requires both bounds for the potential and invertibility of the corresponding boundary operator. For the Laplacian on Lipschitz domains, singular integral estimates, including the boundedness of the Cauchy integral proved in \cite{CMM82}, and Rellich identities lead to the layer potential theory of \cite{Ver84}. Flattening a Lipschitz graph produces a real symmetric matrix independent of the transverse variable, so the upper half-space theory with measurable coefficients also contains this geometric setting.
	
	For rough coefficients, square function estimates and first order methods allow one to obtain bounds without pointwise estimates for derivatives of the fundamental solution. The Kato square root theorem and the functional calculus of perturbed Dirac operators \cite{Kato02,AKM06} provided tools for this approach. Stability of well-posedness and of layer potential invertibility was established in \cite{AAM10,AAAHK11}, including real symmetric matrices independent of the transverse variable and small complex perturbations of them. For real, possibly non-symmetric, $t$-independent coefficients, the solutions of the Dirichlet and regularity problems admit layer potential representations \cite{HKMP15b}. For general bounded accretive systems independent of the transverse variable, Ros\'en \cite[Theorem~7.1 and the estimates following it]{Rosen13} obtained $L^2$ slice, square function, and averaged non-tangential estimates for the single layer gradient, independently of well-posedness. The averaged estimate does not assert the ordinary pointwise maximal bound for arbitrary rough coefficients. Under suitable interior estimates, the single layer potential defined by the functional calculus agrees with the integral against the fundamental solution. For the construction of the fundamental solution for rough coefficients, see \cite{HK07}.
	
	Once $L^2$ bounds are available, interior estimates permit extrapolation to other boundary spaces and the derivation of jump relations. This gives the $L^p$ and endpoint theory of \cite{HMiMo15}, with invertibility remaining a separate requirement for solvability. For coefficients depending on the transverse variable, Hofmann, Mayboroda, and the author \cite[Theorem~1.12]{HMM15} proved layer potential bounds under a small Carleson condition on the Whitney supremum of the difference $A(x,t)-A_0(x)$, where $A_0$ is independent of $t$. The operator with matrix $A_0$ and its adjoint are assumed to satisfy interior De Giorgi--Nash estimates. The conclusions include averaged maximal, square function, and boundary trace bounds. The coefficients may be complex and need not be differentiable. For the exponents considered in \cite[Theorem~1.35]{HMM15}, solvability by layer potentials follows when the corresponding boundary operators for $A_0$ are invertible.
	
	Our positive result instead assumes smallness of the DKP norm for a real matrix and uses the reflected extension \eqref{sl:reflection-matrix}. It gives ordinary maximal, slice, and tangential trace estimates without comparison with a prescribed $t$-independent matrix. The counterexample shows that a finite DKP norm does not suffice. The single layer potential depends on the coefficients in both half-spaces, while the upper Neumann problem depends only on the coefficients in the upper half-space. If $A$ is independent of $t$ and we use the same matrix in the lower half-space, the coefficients are independent of $t$ in the whole space and \cite{Rosen13} applies. Its reflected extension can have different layer estimates while the upper Neumann problem remains unchanged. For the planar family in Subsection~\ref{sec:intro-neumann-history}, Axelsson \cite{Axe10} showed that the boundary equation method and the Lax--Milgram method can give different solutions when the skew-symmetric part of the matrix is sufficiently large.
	
	\subsection{The proof of Theorem~\ref{thm:main}}\label{sec:intro-neumann-proof}
	Let us now describe the main ideas behind the proof of Theorem \ref{thm:main}. Our strategy is inspired by the one in \cite{Geo}. There, we modify the domain near a
	flat piece of the boundary on which the Neumann datum vanishes, in such a way that
	the gradient of the solution changes little in $L^2$, while the integral of the
	absolute value of its tangential gradient increases on a new collection of flat
	pieces of smaller measure. Since these new pieces are flat as well, we can then
	repeat the construction on them.
	
	\smallskip
	\noindent\textit{The choice of the coefficient.}
	 Here we keep the half-space fixed and modify the coefficients. For the iteration to work, we need the boundary
	integral to increase on cubes above which the coefficient is still the identity, so
	that the same construction can be performed again inside them.
	
	The natural candidates for these modifications are skew-symmetric perturbations of
	the identity. This choice is suggested by the free boundary argument in
	\cite[Sections~8.1--8.2]{MMP}, where a skew-adjoint term involving tangential
	derivatives preserves the accretivity of the testing functions. Indeed, for the
	matrix $A_c$ in \eqref{eq:matrix} we have $A_c\xi\cdot\xi=|\xi|^2$, so that the
	perturbation does not affect the coercivity of the operator, while the conormal
	derivative acquires the term $c\cdot\nabla_xu$. Moreover, if $c$ is independent of
	$t$ and $\div_xc=0$, then $L_{A_c}=-\Delta$ in the open half-space, see
	\cite[Remark~6.7]{AMcM13} and \cite[Remark~10.3]{AM14}. In this case the perturbation
	does not change the equation at all, and it only acts through the boundary condition.
	
	\smallskip
	\noindent\textit{The affine function and the periodic perturbation.}
	Consider a boundary ball on which $g=0$ and
	near which $A=I_{n+1}$. The even reflection of $u$ is then harmonic across the ball
	and, since $\partial_tu=0$ there, the first order Taylor polynomial of $u$ at a
	boundary point depends only on the horizontal variables. If the tangential gradient at this point is nonzero, we may take the point as the
	origin, subtract $u(0,0)$, rotate the horizontal coordinates and divide by
	$|\nabla_xu(0,0)|$. The first order Taylor polynomial is then $x_1$. We first perturb the affine function $x_1$, keeping its derivative in the $x_1$ direction equal to one and making its derivative in the $x_2$ direction large on part of the boundary. To this end, we take
	$c=Kh(x_2)\zeta(t)e_1$ and $u_0=x_1+Kw(x_2,t)$, so that $\div_xc=0$ and
	$\partial_1u_0=1$. Writing $s=x_2$, the equation for $u_0$ and its zero conormal
	boundary condition then reduce to
	\begin{equation}\label{eq:intro-periodic}
		w_{ss}+w_{tt}=h(s)\zeta'(t)
		\qquad\text{and}\qquad w_t(s,0)=h(s)\zeta(0).
	\end{equation}
	We choose $h=f'$, where $f$ is a smooth, nonnegative and nonzero periodic function
	supported in a proper subinterval of each period, and we take $\zeta=1$ for the
	$t$-independent matrix and a fixed nonnegative, nonzero $\zeta\in C_c^\infty((1,2))$
	in the DKP example. For both choices of $\zeta$, Lemma~\ref{lem:profiles} provides an
	interval $J$ away from $\supp f$ and a number $d>0$ such that, for $s\in J$, we have
	\[
	|w_s(s,0)|\ge d\qquad \text{ and } \qquad |w_t(s,3)|\ge d.
	\]
	Thus, for $x_2\in J$, we have $|\partial_2u_0(x,0)|\geq Kd$ and
	$|\partial_tu_0(x,3)|\geq Kd$. Since $f=h=0$ on $J$, the matrix equals the
	identity above this interval. In particular, $(A_c\nabla u_0)_{n+1}=\partial_tu_0$ there. This is where we use that $n\geq2$.
	With only one horizontal variable, the condition $\div_xc=0$ would force $c$ to be constant in
	$x$, and so $c$ would vanish whenever its horizontal support is bounded.
	
	\smallskip
	\noindent\textit{One perturbation in a ball.}
	We now make the periodic perturbation in a ball. On a small
	ball $B$ the direction of $a=\nabla_x\Tr u$ is close to a fixed one which, after a
	rotation, we may assume to be $e_1$. Then, we localize the periodic perturbation by means
	of a cutoff $\phi$ in $B$, and we rescale it with a small period $\rho$. More
	precisely, we set
	\[
	\psi_\rho(x,t)\coloneqq\rho K\, \phi(x)\, f\Bigl(\frac{x_2}{\rho}\Bigr)\zeta\Bigl(\frac{t}{\rho}\Bigr)
	\qquad \text{ and }\qquad c_\rho\coloneqq\partial_2\psi_\rho\,e_1-\partial_1\psi_\rho\,e_2.
	\]
	
	Since the mixed derivatives of $\psi_\rho$ agree, the field $c_\rho$ is still
	divergence-free. Notice also that, when a derivative falls on $\phi$, the resulting term
	contains $f$. Since $f$ and $f'$ vanish near $J$ and are periodic, $c_\rho(x,t)=0$ whenever $x_2/\rho\in J+\mathbb Z$. 
	
	We denote by $A^{(0)}$ the original matrix and by $u_\rho$ the solution with matrix
	$A_\rho=A^{(0)}+S(c_\rho)$ and the same datum. As in the calculation for $x_1$, we
	expect $u_\rho$ to be close to $u+v_\rho$, where
	$v_\rho=\rho K\phi(x)a_1(x)w(x_2/\rho,t/\rho)$. Although $\|\nabla v_\rho\|_2\lesssim\sqrt\rho$, on the boundary we have
	\[
	\partial_2v_\rho(x,0)
	=K\phi(x)a_1(x)w_s(x_2/\rho,0)
	+\rho K\partial_2(\phi a_1)(x)w(x_2/\rho,0).
	\]
	On the selected cubes, $\phi=1$ and $|w_s(x_2/\rho,0)|\geq d$. Thus the first term has absolute value at least $Kd|a_1(x)|$, while the second is bounded by $C\rho$.
	
\begin{figure}[!htb]
		\centering
		\begin{tikzpicture}[font=\small,>=Latex,line cap=round,line join=round]
			\path[fill=mathOrange!9] (4.55,0) rectangle (6.37,3.28);
			\draw[mathOrange!75!black,densely dotted] (4.55,0)--(4.55,3.28)
			(6.37,0)--(6.37,3.28);
			\path[draw=mathBlue!80!black,fill=mathBlue!17] (.75,.96) rectangle (2.78,1.92);
			\node at (1.76,1.44) {DKP};
			\draw[->] (0,0)--(0,3.41) node[above] {$t$};
			\draw[->] (0,0)--(7.16,0) node[below right] {$x_2$};
			\foreach \z/\lab in {.96/\rho,1.92/2\rho,2.88/3\rho}{
				\draw (-.08,\z)--(0,\z) node[left=2pt] {$\lab$};
			}
			\node[below left] at (0,0) {$0$};
			\draw[gray,densely dotted] (.06,2.88)--(4.90,2.88);
			\node at (5.46,1.58) {$A_\rho=I_{n+1}$};
			\draw[mathOrange!80!black,line width=1.25pt] (4.90,0)--(6.02,0)
			(4.90,2.88)--(6.02,2.88);
			\draw[mathBlue,thin] (5.46,2.88) circle (.22);
			\draw[mathBlue,thin] (5.68,0) arc[start angle=0,end angle=180,radius=.22];
			\draw[mathBlue,thin,dashed] (5.24,0) arc[start angle=180,end angle=360,radius=.22];
			\fill[mathBlue] (5.46,0) circle (1.2pt) (5.46,2.88) circle (1.2pt);
			\node[below] at (5.46,-.27) {$S_k$};
			\draw[black!60] (6.05,2.88)--(7.12,2.88);
			\node[anchor=west] at (7.13,2.88) {$|\partial_t u_\rho(x,3\rho)|$};
			\draw[black!60] (6.05,.05)--(6.75,.46)--(7.12,.46);
			\node[anchor=west] at (7.13,.46) {$|\nabla_x\Tr u_\rho(x)|$};
		\end{tikzpicture}
		\caption{A vertical section through a selected cube $S_k$. The two balls lie
			over the same $x\in S_k$ and are drawn with their radii enlarged. The dashed
			semicircle shows the even reflection. On $S_k$ the normal derivative vanishes
			at $t=0$. At the upper marked point, $\partial_tu_\rho=(A_\rho\nabla u_\rho)_{n+1}$
			since $A_\rho=I_{n+1}$ above $S_k$. In this period the DKP perturbation is supported in the marked
			rectangle. For the $t$-independent example the corresponding strip extends
			from $t=0$ to infinity.}
		\label{fig:hs-two-heights}
	\end{figure}
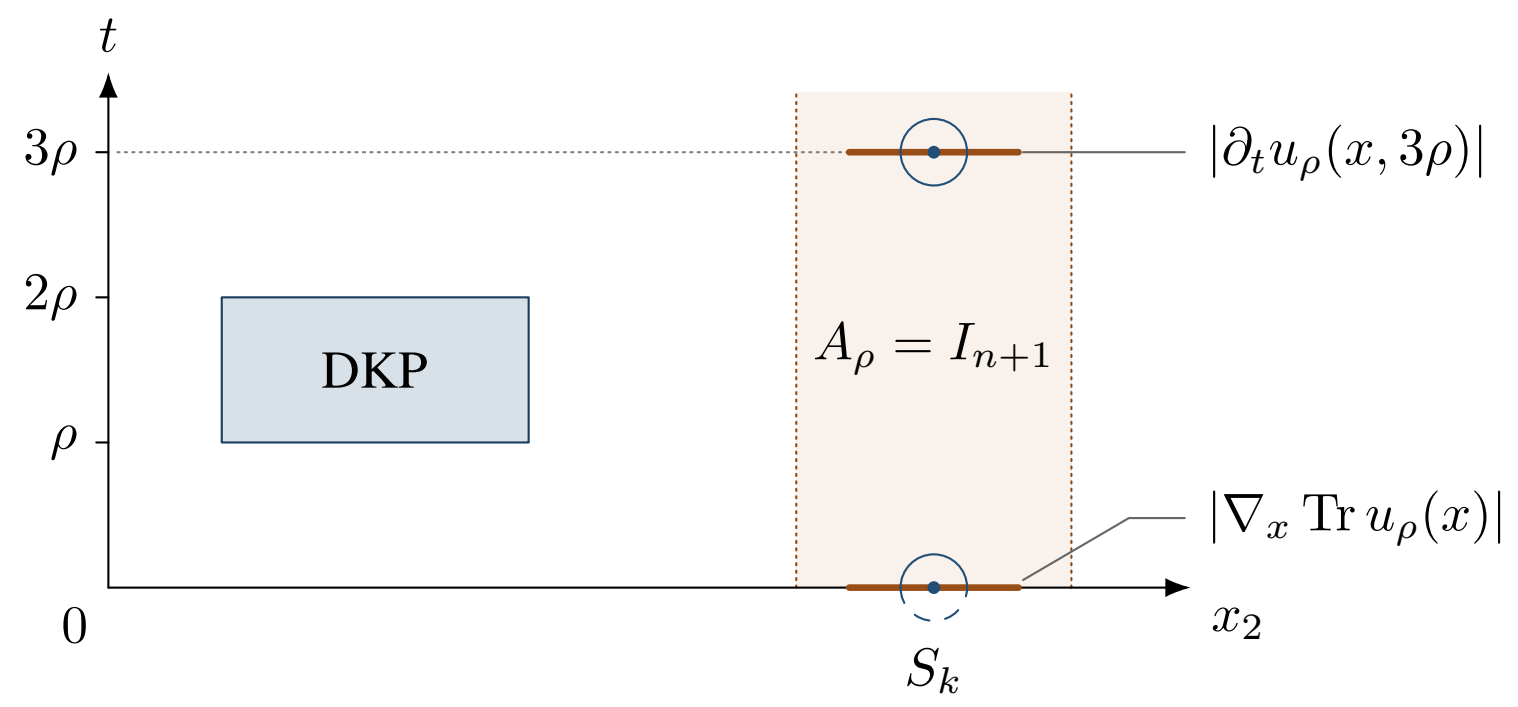

	By \eqref{eq:intro-periodic}, the leading terms in the
	equation for $u+v_\rho$ cancel, and so the error $\|\nabla(u_\rho-u-v_\rho)\|_2$ is
	bounded by $C\rho$. To obtain the lower bounds for $u_\rho$, on each period cube
	with center $x_k$ we compare $u_\rho$ with the function in \eqref{eq:affine} with $a=a(x_k)$.
	Both functions are harmonic above the selected subcube and extend by even reflection
	across the plane there, so we can apply the mean value inequality on balls of radius
	comparable to $\rho$ centered at heights $0$ and $3\rho$, see
	Figure~\ref{fig:hs-two-heights}. Denoting by $S_\rho\Subset B$ the union of the
	selected subcubes, we have $|S_\rho|\leq\theta|B|$ for a fixed $\theta<1$, and
	Lemma~\ref{lem:local} gives
	\[
	\min\Bigl\{\int_{S_\rho}|\nabla_x\Tr u_\rho|\,dx,\,
	\int_{S_\rho}|\partial_tu_\rho(x,3\rho)|\,dx\Bigr\}
	\geq3\int_B|\nabla_x\Tr u|\,dx,
	\qquad \|\nabla(u_\rho-u)\|_2\lesssim\sqrt\rho .
	\]
	Thus the boundary integral increases by a factor of at least three on a set of
	smaller measure. Above this set $A_\rho=I_{n+1}$, so the second integral gives the
	same lower bound for the integral of $|(A_\rho\nabla u_\rho)_{n+1}|$ at height $3\rho$.
	In the geometric example of \cite{Geo}, the estimate for the harmonic function in
	the exterior of a square gives a factor slightly larger than one. Here we obtain
	the factor three by choosing $K$ large enough.
	
	\noindent\textit{Iteration and passage to the limit.}
	Since the coefficient matrix is still the identity in the columns above the selected cubes,
	we can repeat the construction inside these columns. To this end, we choose disjoint
	balls inside the cubes so that the integral of $|\nabla_x\Tr u|$ over the union of the balls
	exceeds two thirds of its integral over the union of the cubes. By Lemma~\ref{lem:local},
	the integral over the cubes selected at the next step is at least three times the integral over the balls.
	Thus the boundary integral more than doubles at each step. The perturbations introduced at different steps have disjoint
	horizontal supports, and this keeps the matrix bounded, independently of the number of
	steps. Moreover, after each step we choose the next period so small that the $L^2$ norms of the differences of successive gradients have finite sum, see Proposition~\ref{prop:finite}. In the DKP case
	we also have to control the two Carleson measures. Each field is supported between
	the heights $\rho$ and $2\rho$, and its derivative is bounded by $C/\rho$, so that
	both $dt/t$ and $t\,dt/\rho^2$ have bounded integral over this interval,
	independently of $\rho$. Since at each step the new cubes occupy at most a fixed
	fraction $\theta$ of their parents, we obtain the packing condition in
	Lemma~\ref{lem:packing}, which allows us to sum the contributions of all the steps in
	Proposition~\ref{prop:dkp}.
	
	To construct the final coefficients, we keep the same nonzero smooth datum $g$
	throughout, and we choose disjoint balls $B_j$ outside $\supp g$ with
	$\sum_j|B_j|<\infty$. For the $t$-independent matrix, $\dist(0,B_j)\to\infty$ as $j\to\infty$. For the DKP matrix, the balls accumulate at the origin. Let $u_{j-1}$ be the solution after the
	perturbations in the first $j-1$ balls. In order to start the construction in $B_j$,
	we need a cube $Q\Subset B_j$ of positive side length with
	$\int_Q|\nabla_x\Tr u_{j-1}|\,dx>0$, and this is guaranteed by Lemma~\ref{lem:seed}.
	Indeed, in the $t$-independent case $u_{j-1}$ is harmonic in the whole upper
	half-space, while in the DKP case each of the finitely many perturbations vanishes for $0<t<\rho$, where $\rho$ is the scale used for that perturbation. So, in both cases, $u_{j-1}$ is harmonic in $\R^n\times(0,h)$ for some $h>0$
	and, since $A_{j-1}=I_{n+1}$ above $B_j$ and $g=0$ on $B_j$, its even extension is
	harmonic across $Q$. If the integral above were zero, the trace of $u_{j-1}$ would be
	constant on $Q$ and, by evenness, $\partial_tu_{j-1}$ would vanish there. By the
	uniqueness of harmonic functions with zero Cauchy data, $u_{j-1}$ would then be equal
	to this constant in the whole slab, and testing the Neumann identity with
	$g(x)\varphi(t)$, where $\varphi\in C_c^\infty((-h,h))$ and $\varphi(0)=1$, would
	give $0=\int_{\R^n}g^2\,dx$, contrary to $g\ne0$.
	
	We can now apply Proposition~\ref{prop:finite} on $Q$ with $H=4$. Denoting by $u_j$
	the resulting solution and by $\rho_j$ the last period used in $B_j$, we obtain a
	selected set $S_j$ on which the integrals of $|\nabla_x\Tr u_j|$ and of
	$|\partial_tu_j(\cdot,t_j)|$, with $t_j=3\rho_j$, are both at least four. The
	coefficient equals $I_{n+1}$ above $S_j$, and all the later perturbations are
	supported outside $B_j$. Thus $(A_k\nabla u_k)_{n+1}=\partial_tu_k$ above $S_j$
	for every $k\geq j$, and $(A\nabla u)_{n+1}=\partial_tu$ there for the limiting
	solution. We need a lower bound uniform in $k$. To this end,
	using harmonic reflection near $S_j$ and interior estimates at height $t_j$, we
	choose the later periods so small that each of the two integrals changes by less than
	one when we pass from $u_j$ to $u_k$, uniformly in $k>j$. In this way,
	Proposition~\ref{prop:construction} gives
	\[
	\int_{S_j}|\nabla_x\Tr u|\,dx\geq1
	\qquad\text{and}\qquad \int_{S_j}|(A\nabla u)_{n+1}(x,t_j)|\,dx\geq1.
	\]
	Let $E=\bigcup_j S_j$. For $x\in S_j$ and $a>1$, the point $(x,t_j)$ lies in
	$\gamma_a(x)$, and $(A\nabla u)_{n+1}=\partial_tu$ is continuous near this point. Since the sets $S_j$ are
	disjoint, we deduce that
	\[
	\int_E N_a((A\nabla u)_{n+1})\,dx\geq\sum_j\int_{S_j}|(A\nabla u)_{n+1}(x,t_j)|\,dx=\infty,
	\]
	and the same holds for $N_a(\partial_tu)$, because $(A\nabla u)_{n+1}=\partial_tu$ there. For
	$N_a(\nabla_xu)$ we argue in the same way, taking the boundary limit on $S_j$ and
	using the first lower bound above. Finally, since $|E|\leq\sum_j|B_j|<\infty$, all
	three maximal functions have infinite $L^p$ norm for every $1\leq p\leq\infty$.
	
	It remains to consider an arbitrary weak solution. Suppose that $v$ is a
	weak solution with the same datum $g$ and that $N_a(\nabla v)\in L^p$ for some
	$1\leq p<\infty$. We will show in Section~\ref{sec:global} that $\nabla v=\nabla u$,
	which is impossible by \eqref{eq:main-components}. Lemma~\ref{lem:boundary-local-energy} gives $v\in Y^{1,2}_{\rm loc}(\overline{\R^{n+1}_+})$. Since $u\in Y^{1,2}(\R^{n+1}_+)$, their difference also belongs to $Y^{1,2}_{\rm loc}$ up to the boundary and has zero Neumann datum. We can therefore reflect $v-u$ evenly across the boundary. The
	reflected matrix equals the identity outside $(\bigcup_jB_j)\times\R$ and, since
	$\sum_j|B_j|<\infty$, the average of $|\mA-I_{n+1}|$ over $B(0,R)$ tends to zero
	as $R\to\infty$. On the other hand, the $L^p$ bound on $N_a(\nabla v)$ and the $L^2$
	bound on $\nabla u$ show that the $L^2$ mean oscillation of the even extension of $v-u$ on
	$B(0,R)$, $R\geq1$, is bounded by a constant times $R^{\max\{0,1-n/p\}}$. For finite $p$ this
	exponent is less than one, and so Lemma~\ref{lem:liouville-at-infinity} gives
	$\nabla v=\nabla u$.
	
	Concerning the structure of the paper, the periodic functions are constructed in
	Section~\ref{sec:profiles} and the perturbation in a ball in Section~\ref{sec:local},
	while Sections~\ref{sec:iteration} and~\ref{sec:carleson} contain the finite
	iteration and the Carleson estimates, respectively.
	
	\subsection{The proofs of Theorems~\ref{sl:negative} and~\ref{sl:small-layer}}\label{sec:intro-layer-proofs}
	Both layer potential results are based on the identity
	\eqref{sl:intro-even-identity}, which relates the single layer potential for the
	reflected coefficients to the Neumann problem in the upper half-space. More
	precisely, if $u_g$ is the variational Neumann solution for $A$ with datum $g$ and
	$\mA$ is the reflected coefficient in \eqref{sl:reflection-matrix}, then
	Lemma~\ref{sl:reflection} shows that
	\(
	\mS_{\mA}(2g)(x,t)=u_g(x,|t|).
	\)
	In other words, the even extension of the Neumann solution is a single layer
	potential, and so the maximal, slice and trace estimates for $\mS_{\mA}(2g)$ follow from the corresponding estimates for $u_g$.
	
	Let us first consider Theorem~\ref{sl:negative} and write $v=\mS_{\mA}(2g)$.
	Since $v$ is even in $t$, its derivative $\partial_t v$ is odd. By
	\eqref{sl:reflection-matrix}, $(\mA\nabla v)_{n+1}$ is also odd. Hence Theorem~\ref{thm:main}
	gives the failure of both maximal estimates in each half-space. On the smaller balls in the proof of Lemma~\ref{lem:normal-components} we have $A=I_{n+1}$, so the same lower bound applies to the local averages of both $\partial_tu_g$ and $(A\nabla u_g)_{n+1}$. For the
	tangential derivative we use the disjoint sets $S_j$ of
	Proposition~\ref{prop:construction}, whose union has finite measure and which satisfy
	$\int_{S_j}|\nabla_x\Tr u_g|\,dx\geq1$ for every $j$. Summing over $j$, we obtain the
	failure of the tangential trace estimate. Moreover, since $u_g$ extends harmonically
	across each $S_j$, the same lower bounds also give the failure of the slice and
	averaged tangential maximal estimates. By \eqref{eq:fixed-datum} and \eqref{eq:initial-balls}, the sets $S_j$ are at distance greater than $1/2$ from $\supp g$. Hence every truncation with $0<\varepsilon<1/2$ agrees with $\nabla_x\Tr v$ on $\bigcup_jS_j$ and has infinite $L^2$ norm there, see Corollary~\ref{sl:boundary-truncations}. Let us emphasize that
	the choice of the reflected coefficients is essential here. Indeed, the Neumann
	solution only depends on the coefficients in the upper half-space, and its even
	extension agrees with the single layer for the reflected coefficients. However, a
	different choice below the boundary can change the single layer estimates, as we show
	in Remark~\ref{rem:two-extensions} for the extension \eqref{eq:unchanged_extension}.
	
	Concerning Theorem~\ref{sl:small-layer}, the same identity transfers the Neumann
	estimate for small Carleson norm to the single layer. We apply
	the averaged Neumann estimate of \cite{DP24} with datum $f/2$ and use
	\eqref{eq:gradient-maximal-comparison} to obtain the bound for $N_a(\nabla u)$. This gives a uniform $L^2$ bound for $\nabla u(\cdot,t)$. Testing $\nabla_xu(\cdot,t)$ against a smooth compactly supported vector field and letting $t\to0^+$ gives the tangential trace estimate. Finally, by approximation of the density, we obtain the layer gradient for all densities in $L^2(\R^n)$. To recover the potential up to an additive constant, we normalize the approximating potentials to have mean zero on the unit ball. All this is done in Subsection~\ref{sl:small-section}.
	
	\subsection*{Acknowledgements}
	I completed this work during my visit to the Centre de Recerca Matem\`atica
	in Bellaterra (Barcelona) for the program ``Analysis of Free Boundary
	Problems'' in the fall of 2026. I thank the institute for its hospitality.

	\subsection*{AI use}
	My starting point was the geometric counterexample in my joint work with
	C.~Puliatti and X.~Tolsa \cite{Geo}. There we kept the Neumann datum fixed
	and changed the domain so that the integral of the absolute value of the
	tangential gradient increased on smaller flat pieces of the boundary.
	I wanted to obtain the same increase in the half-space by changing the
	coefficients. The skew-adjoint term in the construction of the testing function in the free boundary argument of
	\cite{MMP} suggested using a skew-symmetric perturbation of the identity.
	By the observation in \cite{AMcM13,AM14}, if the horizontal field defining
	the perturbation is divergence-free and independent of $t$, the equation
	in the upper half-space is still the Laplace equation, while the conormal
	derivative changes.

	I started with the affine solution $u(x,t)=x_1$ and looked for a
	perturbation that would keep $\partial_1 u=1$ and make $|\partial_2 u|$
	large on smaller boundary cubes. I needed the matrix to equal the identity
	above those cubes so that I could repeat the construction there. To do
	this at each step, I needed a uniform bound for the coefficients and an
	error small enough to preserve the increase in the boundary integral.
	For the DKP example I also needed a cutoff in the transverse variable.

	I used large language models to help find explicit periodic functions and
	cutoffs with these properties and work through calculations and details of the
	proofs. I also used them to draft and revise the text
	and the \LaTeX{} source, and to prepare the figures.
	I checked the calculations and the proofs, chose
	the arguments presented here, and take full responsibility for the paper.

	\section{Preliminaries and notation}\label{sec:sobolev}
	
	Throughout, $C$ and $c$ denote positive constants depending on the dimension and the fixed parameters.
	Their values may change at each occurrence. We write $a\lesssim b$ for $a\leq Cb$ and
	$a\approx b$ for $a\lesssim b\lesssim a$. A subscript in $\lesssim_\gamma$ indicates dependence on
	$\gamma$. Constants with a subscript retain their values. We denote the characteristic function
	of a measurable set $E$ by $\chi_E$ and the horizontal gradient by $\nabla_x=(\partial_1,\ldots,\partial_n)$.
	The outward unit normal to the upper half-space is $-e_{n+1}$. For a function $f$ on $\R^{n+1}_+$, we write
	\[
	\supp_x f=\overline{\{x\in\R^n:f(x,t)\ne0\text{ for some }t>0\}}.
	\]
	
	\subsection{Sobolev spaces}
	We set $2^*\coloneqq2(n+1)/(n-1)$. For an open set $D\subset\R^{n+1}$, we define
	$Y^{1,2}(D)$ to be the space of weakly differentiable functions $v\in L^{2^*}(D)$ such that
	$\nabla v\in L^2(D)$, with norm
	\[
	\|v\|_{Y^{1,2}(D)}\coloneqq\|v\|_{L^{2^*}(D)}+\|\nabla v\|_{L^2(D)}.
	\]
	We follow \cite[Definition~2.5]{HK07}. By
	$v\in Y^{1,2}_{\mathrm{loc}}(D)$ we mean that $v$ and $\nabla v$ are square integrable on compact subsets
	of $D$. We write $v\in Y^{1,2}_{\mathrm{loc}}(\overline{\R^{n+1}_+})$ if the same holds on every bounded
	half-ball up to the boundary.
	We also set
	\[
	Y_0^{1,2}(D)\coloneqq\overline{C_c^\infty(D)}^{\,Y^{1,2}(D)}.
	\]
	For $D=\R^{n+1}$ or $D=\R^{n+1}_\pm$, the Sobolev inequality in the whole space, applied after even reflection in the half-space case, gives
	\[
	\|v\|_{L^{2^*}(D)}\leq C_n\|\nabla v\|_{L^2(D)},
	\qquad\text{ for every }v\in Y^{1,2}(D).
	\]
	Thus, $\|\nabla v\|_2$ is an equivalent norm on each of these spaces. With this norm, $Y^{1,2}(D)$ is a Hilbert space with inner product
	$\langle u,v\rangle\coloneqq\int_D\nabla u\cdot\nabla v$.
	Moreover, $Y_0^{1,2}(\R^{n+1})=Y^{1,2}(\R^{n+1})$, and restrictions of functions in
	$C_c^\infty(\R^{n+1})$ are dense in $Y^{1,2}(\R^{n+1}_\pm)$.
	
	We use the unitary Fourier transform on $\R^n$, given for $f\in L^1(\R^n)$ by
	\[
	\widehat f(\xi)\coloneqq(2\pi)^{-n/2}\int_{\R^n}e^{-ix\cdot\xi}f(x)\,dx.
	\]
	With this normalization, $\widehat{\partial_j f}=i\xi_j\widehat f$ for tempered distributions.
	We define $\dot H^{1/2}(\R^n)$ as the completion of the real-valued functions in $C_c^\infty(\R^n)$ under the norm
	\[
	\|f\|_{\dot H^{1/2}}^2\coloneqq\int_{\R^n}|\xi||\widehat f(\xi)|^2\,d\xi,
	\]
	and use the homogeneous Sobolev embedding to identify its elements with their $L^{2n/(n-1)}(\R^n)$ representatives.
	We denote its real continuous dual by $\dot H^{-1/2}(\R^n)$. In terms of the Fourier transform, we have
	\[
	\|g\|_{\dot H^{-1/2}}^2=\int_{\R^n}|\xi|^{-1}|\widehat g(\xi)|^2\,d\xi
	\qquad\text{ and }\qquad
	\langle g,f\rangle=\int_{\R^n}\widehat g(\xi)\overline{\widehat f(\xi)}\,d\xi.
	\]
	For $g\in\dot H^{-1/2}$ and $f\in\dot H^{1/2}$, the Cauchy--Schwarz inequality with weights $|\xi|^{-1}$
	and $|\xi|$ gives absolute convergence of the last integral. The pairing is real for real $f$ and $g$.
	
	Restriction to the boundary extends to a bounded surjective trace operator
	\[
	\Tr:Y^{1,2}(\R^{n+1}_\pm)\longrightarrow\dot H^{1/2}(\R^n),
	\qquad\text{ with }\qquad Y_0^{1,2}(\R^{n+1}_\pm)=\ker\Tr,
	\]
	and, for $v\in Y^{1,2}(\R^{n+1}_\pm)$, we have
	\begin{equation}\label{eq:trace}
		\|\Tr v\|_{\dot H^{1/2}}\leq\|\nabla v\|_2
		\qquad\text{ and }\qquad
		\|\Tr v\|_{L^{2n/(n-1)}}\lesssim_n\|\nabla v\|_2.
	\end{equation}
	By duality in the boundary Sobolev embedding $\dot H^{1/2}\hookrightarrow L^{2n/(n-1)}$, we obtain
	\begin{equation}\label{eq:negative}
		\|g\|_{\dot H^{-1/2}}\lesssim_n\|g\|_{L^{2n/(n+1)}},
		\qquad\text{ for every }g\in L^{2n/(n+1)}(\R^n).
	\end{equation}
	For the trace and density statements see \cite[Section~2]{AM14}, and for the homogeneous Sobolev embedding see \cite[Chapter~V, Section~1]{Stein}. The trace in \cite{AM14} is defined modulo constants. Here the $L^{2n/(n-1)}(\R^n)$ representative is unique.
	To obtain the constant one in the first inequality of \eqref{eq:trace}, we take $v\in
	C_c^\infty(\overline{\R^{n+1}_+})$ and integrate the derivative of $|\widehat v(\xi,t)|^2$ in $t$. This
	gives
	\[
	|\xi||\widehat v(\xi,0)|^2
	=-2\operatorname{Re}\int_0^\infty|\xi|\widehat v_t\,\overline{\widehat v}\,dt
	\leq\int_0^\infty\bigl(|\widehat v_t|^2+|\xi|^2|\widehat v|^2\bigr)\,dt.
	\]
	Integrating in $\xi$ and applying Plancherel's theorem gives the first inequality of \eqref{eq:trace} for
	smooth $v$. Density gives it on $Y^{1,2}(\R^{n+1}_+)$, and reflection gives it on $Y^{1,2}(\R^{n+1}_-)$.
	For $u\in Y^{1,2}(\R^{n+1})$, we denote its common upper and lower trace by $\Tr u$.
	
	We write $\M$ for the Hardy--Littlewood maximal operator on $\R^n$, $P_t$ for the Poisson semigroup, and $R_j$ for the Riesz transform. Our Fourier multipliers are
	\[
	\widehat{P_tf}(\xi)=e^{-t|\xi|}\widehat f(\xi),\qquad
	\widehat{R_jf}(\xi)=\frac{i\xi_j}{|\xi|}\,\widehat f(\xi),\quad 1\leq j\leq n.
	\]
	
	\subsection{Variational Neumann problem}\label{sec:variational-neumann}
	Given $g\in\dot H^{-1/2}(\R^n)$, we say that $u\in Y^{1,2}(\R^{n+1}_+)$ is a weak solution of the
	{\it variational Neumann problem} with datum $g$ if
	\begin{equation}\label{eq:variational}
		\int_{\R^{n+1}_+}A\nabla u\cdot\nabla v\,dx\,dt
		=\langle g,\Tr v\rangle,
		\qquad\text{ for every }v\in Y^{1,2}(\R^{n+1}_+).
	\end{equation}
	Under \eqref{eq:uniform-ellipticity}, the form on the left-hand side is bounded and coercive in the
	gradient norm. By \eqref{eq:trace}, the right-hand side is a bounded functional on the same space. The
	Lax--Milgram theorem \cite[Section~6.2.1, Theorem~1]{Evans} therefore gives a unique solution. Taking $v=u$ and using \eqref{eq:trace}, we obtain
	\[
	\|\nabla u\|_{L^2(\R^{n+1}_+)}\leq\lambda^{-1}\|g\|_{\dot H^{-1/2}(\R^n)}.
	\]
	For the matrices used in the construction, $\sym A=I_{n+1}$. Testing \eqref{eq:variational} with $v=u$ and using \eqref{eq:trace}, we obtain
	\[
	\|\nabla u\|_2^2=\langle g,\Tr u\rangle
	\leq\|g\|_{\dot H^{-1/2}}\|\nabla u\|_2,
	\]
	and hence $\|\nabla u\|_2\leq\|g\|_{\dot H^{-1/2}}$.
	We take the outward normal in the definition of the conormal derivative. This gives the opposite sign from
	the convention in \cite[Section~2]{AM14}.
	
	When $\sym A=I_{n+1}$, the same calculation allows us to compare $u$ with an approximate solution.
	If $w\in Y^{1,2}(\R^{n+1}_+)$ satisfies
	\[
	\left|\int_{\R^{n+1}_+}A\nabla w\cdot\nabla v\,dx\,dt-\langle g,\Tr v\rangle\right|
	\leq\delta\|\nabla v\|_2,
	\qquad\text{ for every }v\in Y^{1,2}(\R^{n+1}_+),
	\]
	where $\delta\geq0$, then $\|\nabla(w-u)\|_2\leq\delta$.
	Indeed, subtracting \eqref{eq:variational} and testing with $v=w-u$ gives
	\[
	\|\nabla(w-u)\|_2^2\leq\delta\|\nabla(w-u)\|_2.
	\]
	
	Let $\cA$ be a real bounded uniformly elliptic matrix on $\R^{n+1}$, with operator norm bounded by $\Lambda$,
	and let $u\in Y^{1,2}(\R^{n+1})$ solve $L_{\cA}u=0$ weakly in both open half-spaces.
	We use the variational conormal derivative, as in \cite[Section~7.1]{HMM15}.
	For $h\in\dot H^{1/2}(\R^n)$, we choose $\varphi_\pm\in Y^{1,2}(\R^{n+1}_\pm)$ with trace $h$ and define
	\begin{equation}\label{sl:weak-conormals}
		\langle\partial_{\nu_{\cA}}^{\pm}u,h\rangle
		\coloneqq\int_{\R^{n+1}_\pm}\cA\nabla u\cdot\nabla\varphi_\pm\,dx\,dt.
	\end{equation}
	The difference of any two extensions of $h$ lies in $Y_0^{1,2}(\R^{n+1}_\pm)$. By density we can use this
	difference as a test function in the weak equation. Thus,
	\eqref{sl:weak-conormals} is independent of the choice of $\varphi_\pm$.
	Taking $\varphi_\pm(x,t)=P_{|t|}h(x)$, we obtain by Plancherel's theorem $\|\nabla\varphi_\pm\|_2=\|h\|_{\dot H^{1/2}}$. The Cauchy--Schwarz inequality therefore gives
	\[
	\|\partial_{\nu_{\cA}}^{\pm}u\|_{\dot H^{-1/2}}\leq\Lambda\|\nabla u\|_{L^2(\R^{n+1}_\pm)}.
	\]
	
	In the local perturbation argument we will use the following weighted estimate for $v$.
	
	\begin{lemma}\label{lem:weighted}
		Let $E\subset\R^n$ be bounded and measurable, $\gamma>0$, and $0<\rho<1$. For every $v\in Y^{1,2}(\R^{n+1}_+)$, we have
		\begin{equation}\label{eq:weighted}
			\int_E\int_0^\infty e^{-2\gamma t/\rho}|v(x,t)|^2\,dt\,dx
			\leq C_{n,\gamma}(1+|E|^{1/n})\rho\|\nabla v\|_2^2.
		\end{equation}
		The constant depends only on $n$ and $\gamma$.
	\end{lemma}
	\begin{proof}
		For $v\in C_c^\infty(\overline{\R^{n+1}_+})$, the fundamental theorem of calculus and the Cauchy--Schwarz inequality give
		\begin{align*}
			v(x,t)&=v(x,0)+\int_0^t v_t(x,s)\,ds,\\
			|v(x,t)|^2&\leq2|\Tr v(x)|^2
			+2t\int_0^t|v_t(x,s)|^2\,ds.
		\end{align*}
		We set $b\coloneqq 2\gamma/\rho$. Since
		\[
		\int_0^\infty e^{-bt}\,dt=b^{-1},\qquad
		\int_s^\infty te^{-bt}\,dt
		=e^{-bs}(s/b+b^{-2})\leq b^{-2},
		\]
		Tonelli's theorem gives
		\[
		\int_E\int_0^\infty e^{-2\gamma t/\rho}|v|^2\,dt\,dx
		\leq\frac\rho\gamma\|\Tr v\|_{L^2(E)}^2
		+\frac{\rho^2}{2\gamma^2}\|v_t\|_2^2.
		\]
		By H\"older's inequality and \eqref{eq:trace}, we have
		\[
		\|\Tr v\|_{L^2(E)}^2
		\leq|E|^{1/n}\|\Tr v\|_{L^{2n/(n-1)}}^2
		\lesssim |E|^{1/n}\|\nabla v\|_2^2.
		\]
		Combining the last two bounds and using $\rho^2\leq\rho$ proves \eqref{eq:weighted} for smooth $v$.
		Now let $v\in Y^{1,2}(\R^{n+1}_+)$ and choose $v_k\in C_c^\infty(\overline{\R^{n+1}_+})$ such that
		$\|\nabla(v_k-v)\|_2\to0$. By the Sobolev inequality, $v_k\to v$ in $L^{2^*}(\R^{n+1}_+)$, so a subsequence
		converges almost everywhere. Fatou's lemma then gives
		\[
		\int_E\int_0^\infty e^{-2\gamma t/\rho}|v|^2\,dt\,dx
		\leq\liminf_k C_{n,\gamma}(1+|E|^{1/n})\rho\|\nabla v_k\|_2^2
		=C_{n,\gamma}(1+|E|^{1/n})\rho\|\nabla v\|_2^2.
		\]
	\end{proof}
	
	For $v\in Y^{1,2}(\R^{n+1}_+)$, we define the tangential derivative of its trace as a distribution by
	\[
	\langle\nabla_x\Tr v,\varphi\rangle
	\coloneqq -\langle\Tr v,\div\varphi\rangle,
	\qquad \varphi\in C_c^\infty(\R^n;\R^n).
	\]
	
	To pass to the limit in integrals of the boundary gradient, we will use even reflection on the part of the
	boundary where the coefficients equal the identity and the datum vanishes.
	
	\begin{lemma}\label{lem:reflection}
		Let $S\subset\R^n$ be compact and let $r>0$. For $s>0$, write $\Omega_s\coloneqq
		\{z\in\R^{n+1}:\dist(z,S\times\{0\})<s\}$. Let $u\in Y^{1,2}(\R^{n+1}_+)$ satisfy
		\[
		\int_{\Omega_{4r}\cap\R^{n+1}_+}\nabla u\cdot\nabla\psi=0
		\qquad(\psi\in C_c^\infty(\Omega_{4r})).
		\]
		The even extension $u^{\sharp}(x,t)\coloneqq u(x,|t|)$ is then smooth and harmonic in $\Omega_{4r}$, with
		\[
		u_t^{\sharp}(x,0)=0,\qquad
		N_a(\nabla u)(x)\geq|\nabla_x\Tr u(x)|
		\quad(x\in S,\ a>1).
		\]
		Further, if $u,v\in Y^{1,2}(\R^{n+1}_+)$ satisfy the same weak identity, we have
		\begin{equation}\label{eq:reflection}
			\sup_{\Omega_r}|\nabla(u^{\sharp}-v^{\sharp})|
			\leq C_nr^{-(n+1)/2}\|\nabla(u-v)\|_{L^2(\R^{n+1}_+)}.
		\end{equation}
		Here $C_n$ depends only on $n$. In particular, the weak solution of \eqref{eq:variational} satisfies the
		hypothesis if $A=I_{n+1}$ in $\Omega_{4r}\cap\R^{n+1}_+$ and $g=0$ on $\{x\in\R^n:\dist(x,S)<4r\}$ in the sense of distributions.
	\end{lemma}
	\begin{proof}
		Approximating $u$ in $Y^{1,2}(\R^{n+1}_+)$ by restrictions of smooth functions and reflecting the
		approximations evenly, we obtain $u^{\sharp}\in Y^{1,2}_{\rm loc}(\Omega_{4r})$, with weak derivatives
		\[
		\nabla_xu^{\sharp}(x,t)=\nabla_xu(x,|t|),\qquad
		u_t^{\sharp}(x,t)=\operatorname{sgn}(t)u_t(x,|t|)
		\quad\text{a.e.}
		\]
		Given $\psi\in C_c^\infty(\Omega_{4r})$, we change variables by $t\mapsto-t$ and use the hypothesis to obtain
		\[
		\int_{\Omega_{4r}}\nabla u^{\sharp}\cdot\nabla\psi
		=\int_{\Omega_{4r}\cap\R^{n+1}_+}\nabla u\cdot
		\nabla\bigl(\psi(x,t)+\psi(x,-t)\bigr)=0.
		\]
		Thus, $u^{\sharp}$ is weakly harmonic in $\Omega_{4r}$ and hence smooth there by interior regularity
		\cite[Section~6.3.1, Theorem~3]{Evans}. By evenness, $u_t^{\sharp}=0$ on the plane. If $v$ satisfies the
		same weak identity, each component of $\nabla(u^{\sharp}-v^{\sharp})$ is harmonic. For $z\in\Omega_r$, the
		ball $B(z,r)$ lies in $\Omega_{4r}$, so the mean-value inequality gives
		\[
		|\nabla(u^{\sharp}-v^{\sharp})(z)|^2
		\leq\avint_{B(z,r)}|\nabla(u^{\sharp}-v^{\sharp})|^2
		\leq C_nr^{-n-1}\|\nabla(u-v)\|_2^2.
		\]
		This proves \eqref{eq:reflection}. To obtain the boundary bound, we fix $a>1$ and $x\in S$, and choose
		$c_a\in(0,1)$ so that $B((x,t),c_at)\subset\gamma_a(x)$ for all $t>0$. The average of $|\nabla u|^2$ over
		each such ball is bounded by $N_a(\nabla u)(x)^2$. As $t \to 0^+$, continuity gives convergence of these
		averages to $|\nabla u^{\sharp}(x,0)|^2=|\nabla_x\Tr u(x)|^2$.
	\end{proof}
	
	\begin{lemma}\label{lem:cauchy}
		Let $D\subset\R^{n+1}_+$ be connected and open, let $B\subset\R^n$ be a ball, and let $r>0$ satisfy
		$B\times(0,r)\subset D$. Suppose that $u$ is harmonic in $D$, extends harmonically to $B\times(-r,r)$, and
		satisfies $u(x,0)=u_t(x,0)=0$ for $x\in B$. Then $u=0$ in $D$.
	\end{lemma}
	\begin{proof}
		If $\alpha$ is a tangential multi-index and $k\geq0$ is an integer, then harmonicity and the fact that $u,u_t$ vanish on $B$ imply
		\[
		\partial_x^\alpha\partial_t^{2k}u(x,0)
		=\partial_x^\alpha(-\Delta_x)^ku(x,0)=0,\qquad
		\partial_x^\alpha\partial_t^{2k+1}u(x,0)
		=\partial_x^\alpha(-\Delta_x)^ku_t(x,0)=0.
		\]
		All the Taylor coefficients of the extension vanish at each point of $B\times\{0\}$. The real analyticity
		of harmonic functions gives $u=0$ in a neighborhood of each of these points. Since $D$ is connected, the
		identity theorem gives $u=0$ throughout $D$ \cite[Theorems~1.27--1.28]{ABR}.
	\end{proof}

	\subsection{Two auxiliary facts}
	
	To compare an arbitrary weak Neumann solution with the variational solution in Section~\ref{sec:global}, we
	first need square integrability of its gradient up to the boundary. We will then use a Liouville theorem
	for the difference after reflection.
	
	{
		\begin{lemma}\label{lem:boundary-local-energy}
			Let $c\in C^\infty(\R^{n+1}_+;\R^n)\cap L^\infty$ satisfy $\div_xc=0$, and let $A=I_{n+1}+S(c)$. Suppose
			that, in a neighborhood of each boundary point $(x_0,0)$ with $x_0\ne0$, $c$ is smooth up to the
			boundary and independent of $t$. Let $g\in C_c^\infty(\R^n)$ vanish near the origin. If $v\in
			Y^{1,2}_{\mathrm{loc}}(\R^{n+1}_+)$ satisfies
			\begin{equation}\label{eq:ambient-local-neumann}
				\int_{\R^{n+1}_+}A\nabla v\cdot\nabla\varphi
				=\int_{\R^n}g(x)\varphi(x,0)\,dx
				\qquad(\varphi\in C_c^\infty(\R^{n+1}))
			\end{equation}
			and $N_a(\nabla v)\in L^1_{\mathrm{loc}}(\R^n)$ for some $a>1$, then
			$v\in Y^{1,2}_{\mathrm{loc}}(\overline{\R^{n+1}_+})$.
		\end{lemma}
		\begin{proof}
			Almost every $(x,t)\in\R^{n+1}_+$ is a Lebesgue point of $\nabla v$
			in the interior of $\gamma_a(x)$. Thus $|\nabla v(x,t)|\leq N_a(\nabla v)(x)$.
			By absolute continuity on almost every vertical line, the limit
			$f(x)=\lim_{t\to0^+}v(x,t)$ exists for almost every $x$ and
			\[
			|v(x,t)-f(x)|\leq tN_a(\nabla v)(x).
			\]
			For each bounded cube $Q\subset\R^n$, choose $t_0>0$ such that $v(\cdot,t_0)\in L^1(Q)$.
			The last inequality then gives $f\in L^1_{\mathrm{loc}}(\R^n)$ and
			$v(\cdot,t)\to f$ in $L^1_{\mathrm{loc}}(\R^n)$ as $t\to0^+$.
			In particular, $v$ and $\nabla v$ are integrable on every bounded half-ball.
			For $\psi\in C_c^\infty(\R^n;\R^n)$, integration by parts on almost every horizontal slice and passage to the limit give
			\[
			\left|\int_{\R^n} f\,\div_x\psi\,dx\right|
			\leq\int_{\R^n}N_a(\nabla v)|\psi|\,dx.
			\]
			Hence $\nabla_x f$ is a measure whose total variation is bounded by $N_a(\nabla v)\,dx$.
			It follows that $f\in W^{1,1}_{\mathrm{loc}}(\R^n)$ and, by the Sobolev inequality,
			$f\in L^{n/(n-1)}_{\mathrm{loc}}(\R^n)$.

			Fix $x_0\ne0$ and take $r>0$ small enough that $A$ is smooth up to the plane
			in a neighborhood of the closed upper half of $B((x_0,0),6r)$.
			For $\varphi\in C_c^\infty(B((x_0,0),6r))$, integration in $t$ gives
			\[
			\int_{\R^{n+1}_+} A\nabla v\cdot\nabla\varphi
			=\int_{\R^n}g(x)\varphi(x,0)\,dx
			=-\int_{\R^{n+1}_+}g(x)\,\partial_t\varphi(x,t)\,dx\,dt.
			\]
			This is the weak conormal identity in \cite[(1.2)]{Choi20} with vector field $-g(x)e_{n+1}$,
			which is square integrable on the half-ball. Since $\nabla_x A$ is bounded there,
			the mean oscillation of $A$ in the tangential variables at scale $s$ is at most $Cs$,
			and $\int_0^r Cs\,ds/s<\infty$. We already have $v\in W^{1,1}$ on the half-ball.
			After translation and dilation, \cite[Theorem~1.2(a)]{Choi20} with exponent $2$
			gives $v\in W^{1,2}(B((x_0,0),r)\cap\R^{n+1}_+)$.

			It remains to consider the origin. Write $v^\sharp(x,t)=v(x,|t|)$,
			and let $\mA$ be as in \eqref{sl:reflection-matrix}. Choose $r_0>0$ such that
			$g(x)=0$ for $|x|<2r_0$. By even reflection, $v^\sharp$ belongs to
			$W^{1,2}_{\mathrm{loc}}(B(0,2r_0)\setminus\{0\})$. For
			$\varphi\in C_c^\infty(B(0,2r_0)\setminus\{0\})$, testing
			\eqref{eq:ambient-local-neumann} with $\varphi(x,t)+\varphi(x,-t)$ gives
			\[
			\int_{B(0,2r_0)}\mA\nabla v^\sharp\cdot\nabla\varphi=0.
			\]
			We take the locally H\"older continuous representative of $v^\sharp$ on the punctured ball
			\cite[Theorem~8.22]{GT}. For $|z|=r<r_0$, the ball $B(z,r/2)$ lies in this punctured ball.
			The local boundedness estimate \cite[Lemma~2.4, \textup{(2.10)}]{HK07} with exponent one gives
			\[
			|v^\sharp(z)|\leq C r^{-n-1}\int_{B(z,r/2)}|v^\sharp(y,t)|\,dy\,dt.
			\]
			On $B(z,r/2)$ we have $|y|<2r$, $|t|<2r$, and each vertical section has length at most $r$.
			Using $|v^\sharp(y,t)|\leq |f(y)|+|t|N_a(\nabla v)(y)$ and then H\"older's inequality, we obtain
			\[
			\begin{aligned}
			r^{n-1}|v^\sharp(z)|
			&\leq C r^{-1}\int_{|y|<2r}|f(y)|\,dy
			+C\int_{|y|<2r}N_a(\nabla v)(y)\,dy\\
			&\leq C\left(\int_{|y|<2r}|f(y)|^{n/(n-1)}\,dy\right)^{(n-1)/n}
			+C\int_{|y|<2r}N_a(\nabla v)(y)\,dy .
			\end{aligned}
			\]
			Both terms tend to zero as $r\to0^+$. Thus
			$\sup_{|z|=r}|v^\sharp(z)|=o(r^{1-n})$. Since $n\geq2$ and $\mA$ is
			real, bounded and uniformly elliptic, Serrin's theorem
			\cite[Theorem~6, p.~236]{Ser65} applies. Indeed, the last bound gives
			$v^\sharp(z)\geq-|z|^{1-n}$ for small $|z|$ and excludes a singularity
			comparable in absolute value to $|z|^{1-n}$. Hence $v^\sharp$ extends
			across the origin as a $W^{1,2}_{\mathrm{loc}}$ weak solution.
			Restricting to the upper half-space proves the lemma.
		\end{proof}
		
		\begin{lemma}\label{lem:liouville-at-infinity}
			Let $\cA$ be a bounded real uniformly elliptic matrix on $\R^d$, $d\geq3$, such that
			\begin{equation}\label{eq:identity-at-infinity}
				\lim_{R\to\infty}\avint_{B(0,R)}|\cA-I_d|=0.
			\end{equation}
			Suppose that $w\in Y^{1,2}_{\mathrm{loc}}(\R^d)$ satisfies $\div(\cA\nabla w)=0$ and, for some $0\leq\beta<1$,
			\[
			\Phi_w(R)\coloneqq \inf_{a\in\R}
			\left(\avint_{B(0,R)}|w-a|^2\right)^{1/2}
			\leq CR^\beta,\qquad R\geq1.
			\]
			Then $w$ is constant.
		\end{lemma}
		\begin{proof}
			We fix $\alpha\in(\beta,1)$. We claim that there exist $\theta\in(0,1/4)$ and $\delta>0$, depending only on
			$\alpha$, $d$, and the ellipticity bounds, such that
			\begin{equation}\label{eq:harmonic-improvement}
				\avint_{B(0,1)}|\widetilde{\cA}-I_d|<\delta
				\quad\Longrightarrow\quad
				\Phi_v(\theta)\leq\theta^\alpha\Phi_v(1)
			\end{equation}
			for every weak solution of $\div(\widetilde{\cA}\nabla v)=0$ in $B(0,1)$, with the same ellipticity bounds.
			
			For every harmonic function $h$ in $B(0,3/4)$ and every $0<\theta<1/4$, the interior gradient estimate gives
			\[
			\Phi_h(\theta)\leq 2\theta\sup_{B(0,1/4)}|\nabla h|
			\leq C_d\theta\left(\avint_{B(0,1/2)}|h|^2\right)^{1/2}.
			\]
			Choose $\theta$ so small that $2^{d/2}C_d\theta<\theta^\alpha/2$. We may subtract the mean of $v$ over $B(0,1)$ and normalize $\Phi_v(1)=1$, since the claim is immediate for constant solutions. Suppose that \eqref{eq:harmonic-improvement} fails for every $\delta>0$. There exist matrices $\cA_k\to I_d$ in
			$L^1(B(0,1))$ and normalized solutions $v_k$ such that $\Phi_{v_k}(\theta)>\theta^\alpha$. By Caccioppoli's
			inequality and Rellich compactness, after passing to a subsequence, $v_k$ converges weakly in
			$Y^{1,2}(B(0,3/4))$ and strongly in $L^2(B(0,1/2))$ to a function $v$. The uniform bound on the matrices
			also gives $\cA_k-I_d\to0$ in $L^2(B(0,1))$. For every $\varphi\in C_c^\infty(B(0,3/4))$,
			\[
			\left|\int_{B(0,3/4)}(\cA_k-I_d)\nabla v_k\cdot\nabla\varphi\right|
			\leq \|\nabla\varphi\|_\infty\|\cA_k-I_d\|_{L^2(B(0,3/4))}\|\nabla v_k\|_{L^2(B(0,3/4))}\longrightarrow0.
			\]
			Passing to the limit in the weak equation gives $\Delta v=0$ in $B(0,3/4)$. Since $\avint_{B(0,1)}|v_k|^2=1$, the strong convergence gives $(\avint_{B(0,1/2)}|v|^2)^{1/2}\leq 2^{d/2}$. Our choice of $\theta$ now yields $\Phi_v(\theta)<\theta^\alpha/2$. Strong convergence on $B(0,\theta)$ implies
			$\Phi_{v_k}(\theta)\to\Phi_v(\theta)$, contradicting $\Phi_{v_k}(\theta)>\theta^\alpha$. This proves
			\eqref{eq:harmonic-improvement}.
			
			By \eqref{eq:identity-at-infinity}, the average of $|\cA-I_d|$ on $B(0,R)$ is less than $\delta$ for every
			sufficiently large $R$. Applying \eqref{eq:harmonic-improvement} to $w(R\,\cdot)$ and $\cA(R\,\cdot)$, we
			obtain
			\[
			\Phi_w(\theta R)\leq\theta^\alpha\Phi_w(R)
			\qquad(R\geq R_0)
			\]
			for some $R_0\geq1$. For $r\geq R_0$, we apply this inequality $k$ times and use the growth bound at radius $\theta^{-k}r$ to obtain
			\[
			\Phi_w(r)\leq\theta^{k\alpha}\Phi_w(\theta^{-k}r)
			\leq Cr^\beta\theta^{k(\alpha-\beta)}\longrightarrow0.
			\]
			Thus, $\Phi_w(r)=0$ for every $r\geq R_0$. The function $w$ is therefore constant on each of the nested balls $B(0,r)$, and hence on $\R^d$.
		\end{proof}
	}

	\section{Periodic boundary value problems}\label{sec:profiles}
	
	We look for a solution with zero conormal derivative whose tangential gradient is large on a set where $A_c=I$. We begin with the affine function $u(x,t)=x_1$. By adding to it a function of $x_2$ and $t$, we can change its derivative in the $e_2$ direction while keeping its derivative in the $e_1$ direction equal to one. This is the solution that we localize in Section~\ref{sec:local}.
	
	The matrix $S(c)$ is skew-symmetric, so for every $\xi\in\R^{n+1}$,
	\begin{equation}\label{eq:ellipticity}
		A_c\xi\cdot\xi=|\xi|^2,
		\qquad |A_c\xi|\leq(1+\|c\|_\infty)|\xi|,
		\qquad \sym A_c=I_{n+1}.
	\end{equation}
	
	For smooth $c:\R^{n+1}_+\to\R^n$ and $u:\R^{n+1}_+\to\R$, multiplying out the matrix in \eqref{eq:matrix} gives
	\begin{align*}
		A_c\nabla u&=(\nabla_xu+cu_t,\ u_t-c\cdot\nabla_xu),\\
		\div(A_c\nabla u)
		&=\Delta u+(\div_xc)u_t-\partial_tc\cdot\nabla_xu.
	\end{align*}
	The terms containing $c\cdot\nabla_x u_t$ cancel when we take the divergence. Hence, if $\div_xc=0$, we obtain
	\begin{equation}\label{eq:operator}
		\div(A_c\nabla u)=\Delta u-\partial_tc\cdot\nabla_xu,
		\qquad
		\partial_{\nu_{A_c}}u=-u_t+c(x,0)\cdot\nabla_xu.
	\end{equation}
	In the conormal identity we use the outward normal $-e_{n+1}$ and assume that $c$ and $u$ are smooth up to the boundary. When $c$ is independent of $t$, this is the cancellation observed in \cite[Remark~6.7]{AMcM13}.
	To see how to use the boundary term, in dimension three we take $c=b(x_2,t)e_1$ and $u_0=x_1+v(x_2,t)$. Then
	\begin{equation}\label{eq:periodic-matrix-product}
		A_c\nabla u_0=(1+bv_t,\ v_{x_2},\ v_t-b),
	\end{equation}
	and hence $\div(A_c\nabla u_0)=\Delta_{x_2,t}v-b_t$ and $\partial_{\nu_{A_c}}u_0=b-v_t$ on $t=0$.
	
	Let $K>0$, let $h$ be a smooth one-periodic function of mean zero, and let $\zeta\in C^\infty([0,\infty))$.
	We set $b(x_2,t)\coloneqq Kh(x_2)\zeta(t)$ and $v(x_2,t)\coloneqq Kw_\zeta(x_2,t)$.
	If we write $s\coloneqq x_2$, the equation for $u_0$ and the requirement that its conormal derivative vanish reduce to
	\begin{equation}\label{eq:profile-equation}
		(w_\zeta)_{ss}+(w_\zeta)_{tt}=h(s)\zeta'(t),\qquad
		(w_\zeta)_t(s,0)=\zeta(0)h(s).
	\end{equation}
	
	\subsection{The periodic datum and the cutoff functions}
	
	On $\bT=\R/\mathbb{Z}$, with normalized Lebesgue measure, take $0\leq f\in C^\infty(\bT)$,
	$f\not\equiv0$, supported in a proper subinterval of $(0,1)$. We take $h=f'$ so that $h$ has
	mean zero. We shall use $f$ as a primitive when we localize $c$ in Section~\ref{sec:local},
	where both $h$ and $f$ must vanish on the interval on which we test the tangential derivative.
	Fix $0\leq\eta\in C_c^\infty((1,2))$, $\eta\not\equiv0$, and extend it by zero to
	$[0,\infty)$. We shall use $\zeta=1$ and $\zeta=\eta$.
	
	\begin{lemma}\label{lem:profiles}
		For each $\zeta\in\{1,\eta\}$, \eqref{eq:profile-equation} has a unique solution $w_\zeta\in C^\infty(\bT\times[0,\infty))$ with $w_\zeta,\nabla w_\zeta\in L^2(\bT\times(0,\infty))$ and zero mean in $s$ at each height. For all integers $j,k\geq0$ and both solutions,
		\begin{equation}\label{eq:profile-decay}
			|\partial_s^j\partial_t^k w_\zeta(s,t)|\le C_{j,k} e^{-2\pi(t-2)},\qquad t\ge2.
		\end{equation}
		There are an open interval $J$ with $\overline J\subset(0,1)\setminus\supp f$ and a number $d>0$ such that, for both pairs
		\begin{equation}\label{eq:profiles}
			(\zeta,w)\in\{(1,w_1),(\eta,w_\eta)\},
		\end{equation}
		we have
		\begin{equation}\label{eq:profile-sign}
			|w_s(s,0)|\ge d,\qquad |w_t(s,3)|\ge d,\qquad s\in\overline J .
		\end{equation}
	\end{lemma}
	
	\begin{proof}
		Let $h_m$ be the Fourier coefficient of $h$ corresponding to $e^{2\pi i m s}$. We have $h_0=0$. We take the constant Fourier coefficient of $w$ equal to zero. For $m\ne0$, put $\kappa=2\pi|m|$ and set
		\[
		w_m(t)=-\frac{h_m\zeta(0)}{\kappa}e^{-\kappa t}
		-\frac{h_m}{2\kappa}\int_0^\infty
		\bigl(e^{-\kappa|t-r|}+e^{-\kappa(t+r)}\bigr)\zeta'(r)\,dr .
		\]
		Indeed, the derivative of the first exponential in the integral has a jump of $-2\kappa$ at $t=r$, and the derivatives of the two exponentials cancel at $t=0$ for $r>0$. Hence $w_m''-\kappa^2w_m=h_m\zeta'$ and $w_m'(0)=h_m\zeta(0)$. The formula gives $\|w_m\|_\infty\lesssim |h_m|/\kappa$ and $\|w_m'\|_\infty\lesssim |h_m|$, since $\|\zeta'\|_1$ is bounded for our two choices. The differential equation then gives $\|w_m^{(k)}\|_\infty\le C_k|h_m|\kappa^{k-1}$ for every $k\ge0$, by differentiating it and using the bounds for the derivatives of $\zeta$. Since $h$ is smooth, the coefficients $h_m$ decay faster than any power of $|m|$. The series for $w$ and all its derivatives thus converge uniformly. The resulting smooth function is real since $w_{-m}=\overline{w_m}$, and it satisfies \eqref{eq:profile-equation}.
		
		For $t\ge2$ the equation for $w_m$ is homogeneous, and its boundedness gives $w_m(t)=w_m(2)e^{-\kappa(t-2)}$. Summing the differentiated series and using $\kappa\ge2\pi$ proves \eqref{eq:profile-decay}. This also gives $w,\nabla w\in L^2(\bT\times(0,\infty))$, since they are bounded on $\bT\times[0,2]$. For uniqueness, the difference of two solutions has zero constant coefficient. For $m\ne0$, its $m$th Fourier coefficient $v_m$ is a square-integrable solution of $v_m''-\kappa^2v_m=0$ with $v_m'(0)=0$. Square integrability leaves only a multiple of $e^{-\kappa t}$, and the boundary condition makes that multiple zero.
		
		We next find a common interval for the two inequalities in \eqref{eq:profile-sign}. At $t=0$ the exponentials under the integral agree. Integration by parts gives
		\begin{equation}\label{eq:profile-trace}
			w_m(0)=-\frac{h_m}{\kappa}\left(\zeta(0)+\int_0^\infty e^{-\kappa r}\zeta'(r)\,dr\right)
			=-h_m\int_0^\infty e^{-\kappa r}\zeta(r)\,dr .
		\end{equation}
		For both choices of $\zeta$ the last integral is positive. Since $h=f'$ is nonzero and has mean zero, we have $h_m\ne0$ for some $m\ne0$, and hence $w_m(0)\ne0$ for both solutions.
		
		Fix an open interval $I$ with $\overline I\subset(0,1)\setminus\supp f$. Both solutions are harmonic in $0<t<1$ and have $w_t(s,0)=0$ on $I$. Even reflection across $I$ therefore makes $w_s(\cdot,0)$ real analytic there. If it vanished on a subinterval, $w$ would have constant trace and zero normal derivative there. Extend $w$ periodically in $s$ to $\R\times(0,1)$ and subtract this constant. The proof of Lemma~\ref{lem:cauchy} also applies in two dimensions, so the difference vanishes throughout the strip. This contradicts \eqref{eq:profile-trace}.
		
		At height $3$ the function $w_t(\cdot,3)$ is real analytic because $w$ is harmonic for $t>2$. For $\zeta=1$ its $m$th coefficient is $h_m e^{-3\kappa}$. For $\zeta=\eta$ and $t\ge2$, the formula for $w_m$ and integration by parts give
		\[
		w_m(t)=-\frac{h_m}{\kappa}e^{-\kappa t}
		\int_1^2\cosh(\kappa r)\eta'(r)\,dr
		=h_m e^{-\kappa t}\int_1^2\sinh(\kappa r)\eta(r)\,dr.
		\]
		The last integral is positive, so $w_m'(3)\ne0$ whenever $h_m\ne0$. Thus neither of the two functions $w_t(\cdot,3)$ vanishes identically. On any closed subinterval in $I$, the zero sets of these two functions and of the two boundary derivatives $w_s(\cdot,0)$ are finite by real analyticity. We may choose a closed interval $\overline J$ of positive length in $I$ disjoint from all four zero sets. The minimum of the four absolute values on $\overline J$ is positive. This proves \eqref{eq:profile-sign}.
	\end{proof}
	
	\subsection{The rescaled solutions}
	
	Fix one closed cube with sides parallel to the coordinate axes and positive side length such that
	\begin{equation}\label{eq:moments}
		\widehat Q_0\Subset Q_0\cap\{x:x_2\in J\},\qquad Q_0=(0,1)^n,
		\qquad \theta=|\widehat Q_0|<1.
	\end{equation}
	Choose $K\ge1$ once for both pairs, large enough that
	\begin{equation}\label{eq:choice-K}
		0.9\,\theta\,(0.9Kd-0.1)>4 .
	\end{equation}
	For $n\ge2$, $a\in\R^n$, and $\rho>0$, we write $(s,\tau)=(x_2/\rho,t/\rho)$ and set
	\begin{equation}\label{eq:affine}
		c_*(x,t)=Kh(s)\zeta(\tau)e_1,\qquad
		v(x,t)=a\cdot x+\rho Ka_1 w(s,\tau).
	\end{equation}
	We have $\partial_1v=a_1$, $\nabla_xv=a+Ka_1w_s e_2$ and $v_t=Ka_1w_t$. Moreover $\div_x c_*=0$. By \eqref{eq:operator}, the divergence of $A_{c_*}\nabla v$ is $(Ka_1/\rho)(w_{ss}+w_{tt}-h\zeta')=0$, and the outward conormal derivative on the boundary is $Ka_1(h\zeta(0)-w_t(\cdot,0))=0$. Whenever $s\in J+\mathbb Z$, the matrix equals the identity at every height and \eqref{eq:profile-sign} gives
	\[
	|\nabla_x v(x,0)|\ge Kd|a_1|-|a_2|,\qquad
	|(A_{c_*}\nabla v)_{n+1}(x,3\rho)|=|v_t(x,3\rho)|\ge Kd|a_1|.
	\]
	For $a_1\ge0.9|a|$ and $|a_2|\le0.1|a|$, both lower bounds are at least $(0.9Kd-0.1)|a|$.
	
	\begin{figure}[H]
		\centering
		\begin{tikzpicture}[font=\small,>=Latex]
			\begin{scope}[shift={(0,0)}]
				\draw[->] (0,0)--(4.4,0) node[right] {$s$};
				\draw[->] (0,0)--(0,2.1);
				\draw[mathBlue,line width=1pt]
				(.45,0)..controls(.65,.02)and(.65,1.5)..(1.15,1.5)
				..controls(1.65,1.5)and(1.65,.02)..(1.85,0);
				\node[mathBlue] at (1.1,1.9) {$f\geq0$};
				\path[fill=mathOrange!16] (2.6,-.14) rectangle (3.8,.14);
				\draw[mathOrange,line width=1.5pt] (2.6,0)--(3.8,0);
				\node[above,mathOrange] at (3.2,.14) {$J$};
				\draw[mathOrange,line width=1pt] (2.6,-.65)--(3.8,-.65);
				\node[below,align=center] at (3.2,-.7)
				{$f=h=0$\\$|w_s(s,0)|,\ |w_t(s,3)|\geq d$};
				\node[below] at (0,-.1) {$0$};
				\node[below] at (4.1,-.1) {$1$};
				\node[align=center] at (1.0,-.65) {$h=f'$};
			\end{scope}
			\begin{scope}[shift={(6.5,-.2)}]
				\path[fill=lightFill] (0,0) rectangle (4,3.15);
				\path[fill=mathBlue!10] (0,.45) rectangle (4,1.35);
				\path[fill=mathOrange!10] (0,2.0) rectangle (4,2.85);
				\draw[mathBlue] (0,0) rectangle (4,3.15);
				\draw[densely dashed,mathOrange] (0,2.0)--(4,2.0) (0,2.85)--(4,2.85);
				\path[fill=mathOrange!32,draw=mathOrange,line width=.9pt]
				(1.4,2.15) rectangle (2.25,2.7);
				\node at (1.82,2.43) {$\widehat Q_0$};
				\node[align=center,anchor=west] at (4.15,2.43)
				{$A=I$\\near $\widehat Q_0\times[0,\infty)$};
				\foreach \xx in {.4,1.6,2.8}{
					\draw[->,mathBlue,line width=.8pt] (\xx,.65)--({\xx+.65},.65);
					\draw[->,mathBlue,line width=.8pt] ({\xx+.65},1.13)--(\xx,1.13);
				}
				\node[mathBlue,fill=white,inner sep=1.5pt] at (2,1.63) {$Kh(y_2)\zeta(\tau)e_1$};
				\node[anchor=south east] at (4,3.17) {$Q_0$};
				\draw[->] (0,-.35)--(1,-.35) node[right] {$e_1$};
				\draw[->] (-.35,0)--(-.35,1) node[above] {$e_2$};
			\end{scope}
		\end{tikzpicture}
		\caption{The periodic function $f$ and one horizontal period cell in dimension
			three, in the coordinates $y=x/\rho$ and $\tau=t/\rho$. The arrows show the direction
			of $c_*(\rho y,\rho\tau)=Kh(y_2)\zeta(\tau)e_1$ at a height where
			$\zeta(\tau)>0$. On the strip $y_2\in J$, we have $c_*=0$ at all heights,
			while $|w_s(y_2,0)|\geq d$ and $|w_t(y_2,3)|\geq d$. The cube $\widehat Q_0$
			is chosen inside this strip.}
		\label{fig:hs-cell}
	\end{figure}
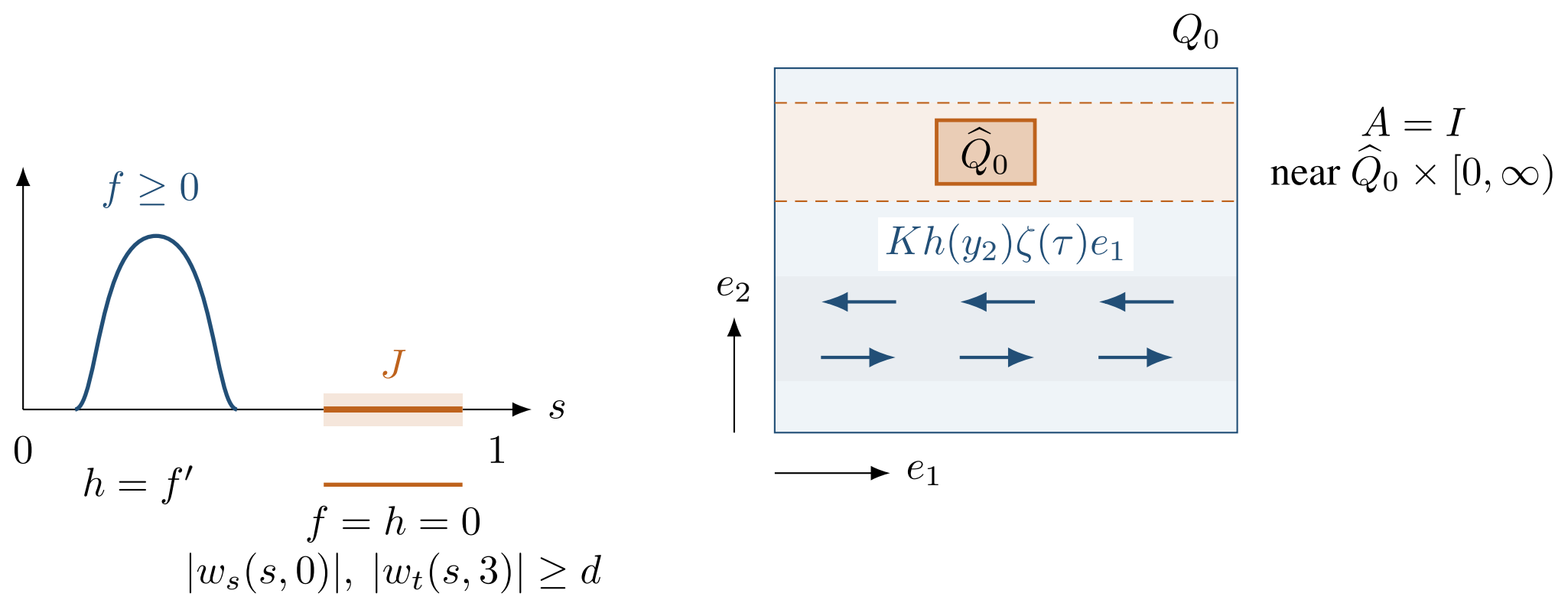
	
	\FloatBarrier
	
	We will also use the following consequence of \eqref{eq:profile-decay}. For every bounded measurable $E\subset\R^n$, $0<\gamma<\pi$ and integers $j,k\ge0$,
	\begin{equation}\label{eq:profile-weighted}
		\int_E\int_0^\infty e^{2\gamma t/\rho}
		|(\partial_s^j\partial_t^kw)(x_2/\rho,t/\rho)|^2\,dt\,dx
		\le C_{j,k,\gamma}|E|\rho .
	\end{equation}
	Indeed, after $t=\rho\tau$ we bound the derivative uniformly on $[0,2]$ and use \eqref{eq:profile-decay} for $\tau\ge2$. The resulting integral in $\tau$ is finite since $2\gamma<4\pi$. In particular, multiplying the correction in \eqref{eq:affine} by a smooth compactly supported horizontal cutoff gives a function in $Y^{1,2}(\R^{n+1}_+)$. The function and its gradient are square integrable by \eqref{eq:profile-weighted}, and the function is bounded, hence also belongs to $L^{2(n+1)/(n-1)}$.

	\section{A local perturbation}\label{sec:local}
	
	Throughout this section, we fix $n\geq2$, a pair $(\zeta,w)$ from
	\eqref{eq:profiles}, and the functions and constants chosen in
	Section~\ref{sec:profiles}.
	
	Let $A^{(0)}\in C^\infty(\overline{\R^{n+1}_+};\R^{(n+1)\times(n+1)})
	\cap L^\infty$ satisfy $\sym A^{(0)}=I_{n+1}$, and let
	$g\in C_c^\infty(\R^n)$. We denote by $u\in Y^{1,2}(\R^{n+1}_+)$ the solution of
	\[
	\int_{\R^{n+1}_+}A^{(0)}\nabla u\cdot\nabla v\,dx\,dt
	=\langle g,\Tr v\rangle
	\qquad\text{for every }v\in Y^{1,2}(\R^{n+1}_+).
	\]
	Let $B_\nu=B(p_\nu,r_\nu)$, $1\leq\nu\leq N$, be balls in $\R^n$
	with pairwise disjoint closures. We assume that $A^{(0)}=I_{n+1}$
	in a neighborhood of each $\overline B_\nu\times[0,\infty)$ and that
	$g$ vanishes near $\bigcup_\nu\overline B_\nu$. For $R_\nu\in O(n)$,
	we write $T_\nu(y)\coloneqq p_\nu+R_\nu y$ and set
	\[
	\widetilde B_\nu=T_\nu^{-1}(B_\nu),\qquad
	u_\nu(y,t)=u(T_\nu(y),t),\qquad
	a_\nu=\nabla_y\Tr u_\nu.
	\]
	We assume that the coordinates on each ball have been chosen so that
	\begin{equation}\label{eq:local-alignment}
		I_\nu\coloneqq \int_{\widetilde B_\nu}|a_\nu(y)|\,dy>0,\qquad
		(a_\nu)_1\geq0.9|a_\nu|,\qquad
		|(a_\nu)_2|\leq0.1|a_\nu|\quad\hbox{on }\widetilde B_\nu.
	\end{equation}
	
	We describe the construction on one ball. We write
	$B=\widetilde B_\nu$, $a=a_\nu$ and $I_B=I_\nu$. In these coordinates
	we write $u$ for $u_\nu$ and $x$ for the horizontal variable. By Lemma~\ref{lem:reflection},
	$a$ is smooth near $\overline B$. We choose a concentric ball
	$B'\Subset B$ such that
	\begin{equation}\label{lp:retained-integral}
		\int_{B'}|a|\,dx\geq0.9I_B,
	\end{equation}
	and $\phi\in C_c^\infty(B)$, $0\leq\phi\leq1$, with $\phi=1$
	near $\overline{B'}$. Since $a$ is smooth near $\overline B$ and
	$\phi$ has compact support in $B$, we extend $\phi a_1$ by zero to a smooth
	function on $\R^n$.
	For $\rho>0$, we write $s\coloneqq x_2/\rho$ and $\tau\coloneqq t/\rho$.
	We perturb the coefficient by $S(c_\rho)$
	and use $v_\rho$ to approximate the corresponding change in $u$, where
	\begin{align}
		c_\rho(x,t)
		&=K\zeta(\tau)\bigl[\phi(x)h(s)e_1
		+\rho f(s)(\phi_2(x)e_1-\phi_1(x)e_2)\bigr],
		\label{eq:local-field}\\
		v_\rho(x,t)&=\rho K\phi(x)a_1(x)w(s,\tau).
		\label{eq:local-correction}
	\end{align}
	Here $\phi_j=\partial_j\phi$ and below $\phi_{ij}=\partial_i\partial_j\phi$. The terms involving derivatives
	of $\phi$ in \eqref{eq:local-field} ensure that $\div_xc_\rho=0$.
	Indeed, with $\psi_\rho\coloneqq\rho K\phi(x)f(s)\zeta(\tau)$, we have
	\[
	c_\rho=\partial_2\psi_\rho\,e_1-\partial_1\psi_\rho\,e_2,
	\qquad \div_xc_\rho=0.
	\]
	
	For each $\nu$, we denote the cutoff on $\widetilde B_\nu$ by
	$\phi_\nu$ and the functions in
	\eqref{eq:local-field}--\eqref{eq:local-correction} by
	$\widetilde c_{\nu,\rho},\widetilde v_{\nu,\rho}$.
	Returning to the original coordinates, we set
	\[
	\begin{gathered}
		c_{\nu,\rho}(T_\nu(y),t)=R_\nu\widetilde c_{\nu,\rho}(y,t),
		\qquad v_{\nu,\rho}(T_\nu(y),t)=\widetilde v_{\nu,\rho}(y,t),\\
		A_\rho=A^{(0)}+\sum_{\nu=1}^N S(c_{\nu,\rho}),\qquad
		v_\rho=\sum_{\nu=1}^Nv_{\nu,\rho}.
	\end{gathered}
	\]
	For fixed $\rho>0$, the matrix $A_\rho$ is bounded and has symmetric
	part $I_{n+1}$. We denote by $u_\rho\in Y^{1,2}(\R^{n+1}_+)$ the solution
	of \eqref{eq:variational} with matrix $A_\rho$ and the same datum $g$.
	We write $e_\rho\coloneqq u_\rho-u-v_\rho$ for the error in approximating
	$u_\rho$ by $u+v_\rho$.
	
	\begin{lemma}\label{lem:local}
		Under the hypotheses above, including \eqref{eq:local-alignment}, there exist
		$C_*,C<\infty$, $0<b<1/4$ and $\rho_0>0$ such that, for every $0<\rho<\rho_0$,
		the following hold.
		\begin{enumerate}[label=\textup{(\roman*)}]
			\item For each $1\leq\nu\leq N$, we have
			\[
			\div_xc_{\nu,\rho}=0,\qquad
			\supp_xc_{\nu,\rho}\Subset B_\nu.
			\]
			Also, $v_\rho\in Y^{1,2}(\R^{n+1}_+)$, $\sym A_\rho=I_{n+1}$, and
			\begin{equation}\label{lp:coefficient-bounds}
				\begin{gathered}
					\|c_{\nu,\rho}\|_\infty\leq C_*,\qquad
					\|\nabla c_{\nu,\rho}\|_\infty\leq C_*/\rho,\\
					\|A_\rho\|_\infty\leq\max\{\|A^{(0)}\|_\infty,1+C_*\}.
				\end{gathered}
			\end{equation}
			For $\zeta=1$, $c_{\nu,\rho}$ is independent of $t$, and for $\zeta=\eta$
			it vanishes unless $\rho<t<2\rho$.
			\item For the solutions with datum $g$, we have
			\begin{equation}\label{lp:local-conclusions}
				\|\nabla(u_\rho-u)\|_2\leq C\sqrt\rho.
			\end{equation}
			\item In the coordinates of each ball, take the period cubes
			$Q_k=\rho(k+[0,1]^n)\subset B'$ and the closed cubes
			$S_k=\rho(k+\widehat Q_0)$, $k\in\mathbb Z^n$. Let $S_{\nu,\rho}$
			be the union of these $S_k$ in the original coordinates and set
			$S_\rho=\bigcup_\nu S_{\nu,\rho}$. The selected cubes are pairwise disjoint,
			$|S_{\nu,\rho}|\leq\theta |B_\nu|$, and
			$A_\rho(x,t)=I_{n+1}$ whenever $\dist(x,S_\rho)<b\rho$ and $t\geq0$.
			Moreover,
			\begin{equation}\label{eq:local-conormal-conclusion}
				\begin{aligned}
					\int_{S_\rho}|\nabla_x\Tr u_\rho|\,dx&\geq3\sum_{\nu=1}^N I_\nu,\\
					\int_{S_\rho}|\partial_tu_\rho(x,3\rho)|\,dx
					=\int_{S_\rho}|(A_\rho\nabla u_\rho)_{n+1}(x,3\rho)|\,dx
					&\geq3\sum_{\nu=1}^N I_\nu .
				\end{aligned}
			\end{equation}
		\end{enumerate}
		The constants $C_*,b$ depend only on $n$ and the fixed data of Section~\ref{sec:profiles}.
		The constants $C$ and $\rho_0$ may also depend on $u$, the balls
		$B_\nu$, and the cutoffs. For each $\nu$, every period cube meeting $\supp\phi_\nu$ lies in $\widetilde B_\nu$.
	\end{lemma}
	
\begin{figure}[htbp]
		\centering
		\begin{tikzpicture}[font=\small,>=Latex,line join=round]
			\node at (2.05,4.75) {Before: $A^{(0)},u$};
			\node at (8.75,4.75) {After: $A_\rho,u_\rho$};
			\path[fill=lightFill,draw=mathBlue,line width=.9pt]
			(2.05,2.65) circle (1.62);
			\node[mathBlue] at (.65,3.95) {$B$};
			\node at (2.05,3.65) {$A^{(0)}=I_{n+1}$};
			\foreach \yy in {1.85,2.45,3.05}{
				\draw[->,mathBlue,line width=.85pt] (1.35,\yy)--(2.75,\yy);
			}
			\node at (2.05,1.45) {$a_1\geq0.9|a|$};
			\draw[->,line width=.9pt] (4.05,2.65)--(6.7,2.65)
			node[midway,above=5pt] {add $S(c_\rho)$}
			node[midway,below=5pt] {same datum $g$};
			
			\begin{scope}
				\clip (8.75,2.65) circle (1.62);
				\fill[lightFill] (7.1,1) rectangle (10.4,4.3);
				\foreach \yy in {1.12,1.93,2.74,3.55}{
					\fill[mathOrange!22] (7.1,\yy) rectangle (10.4,{\yy+.28});
				}
				\foreach \xx in {8.10,8.75,9.40}{
					\foreach \yy in {1.66,2.47,3.28}{
						\path[fill=mathGreen!17,draw=mathGreen!80!black,line width=.65pt]
						({\xx-.145},{\yy-.145}) rectangle ({\xx+.145},{\yy+.145});
					}
				}
			\end{scope}
			\draw[mathBlue,line width=.9pt] (8.75,2.65) circle (1.62);
			\node[mathBlue] at (7.35,3.95) {$B$};
			\draw[<->,mathOrange!85!black] (6.95,1.93)--(6.95,2.74)
			node[midway,left=2pt] {$\rho$};
			\draw[mathOrange!85!black] (9.65,3.7)--(10.85,4.1);
			\node[anchor=west,align=left,text width=3.1cm] at (10.9,4.1)
			{possible support\\of $c_\rho$};
			\draw[mathGreen!80!black] (9.55,2.47)--(10.85,1.9);
			\node[anchor=west,align=left,text width=3.1cm] at (10.9,1.9)
			{selected cubes $S_\rho$\\$A_\rho=I_{n+1}$ near\\$S_\rho\times[0,\infty)$};
			
			\node[align=center] at (2.05,.35)
			{$\displaystyle I_B=\int_B|\nabla_x\Tr u|$};
			\node[align=center] at (8.75,.35)
			{$\displaystyle\int_{S_\rho}|\nabla_x\Tr u_\rho|\geq3I_B$\\[4pt]
				$|S_\rho|\leq\theta|B|$};
			\node[align=center] at (6.4,-.75)
			{$u_\rho=u+v_\rho+e_\rho,\qquad
				\|\nabla e_\rho\|_2\leq C\rho,\qquad
				\|\nabla(u_\rho-u)\|_2\leq C\sqrt\rho$};
		\end{tikzpicture}
		\caption{The perturbation on one ball in two horizontal variables.
			The cubes in $S_\rho$ lie between the strips where $c_\rho$ is supported.
			Above these cubes the coefficient is again the identity. We estimate the tangential gradient on the cubes and the normal derivative at height $3\rho$ above them.}
		\label{fig:hs-local-perturbation}
	\end{figure}
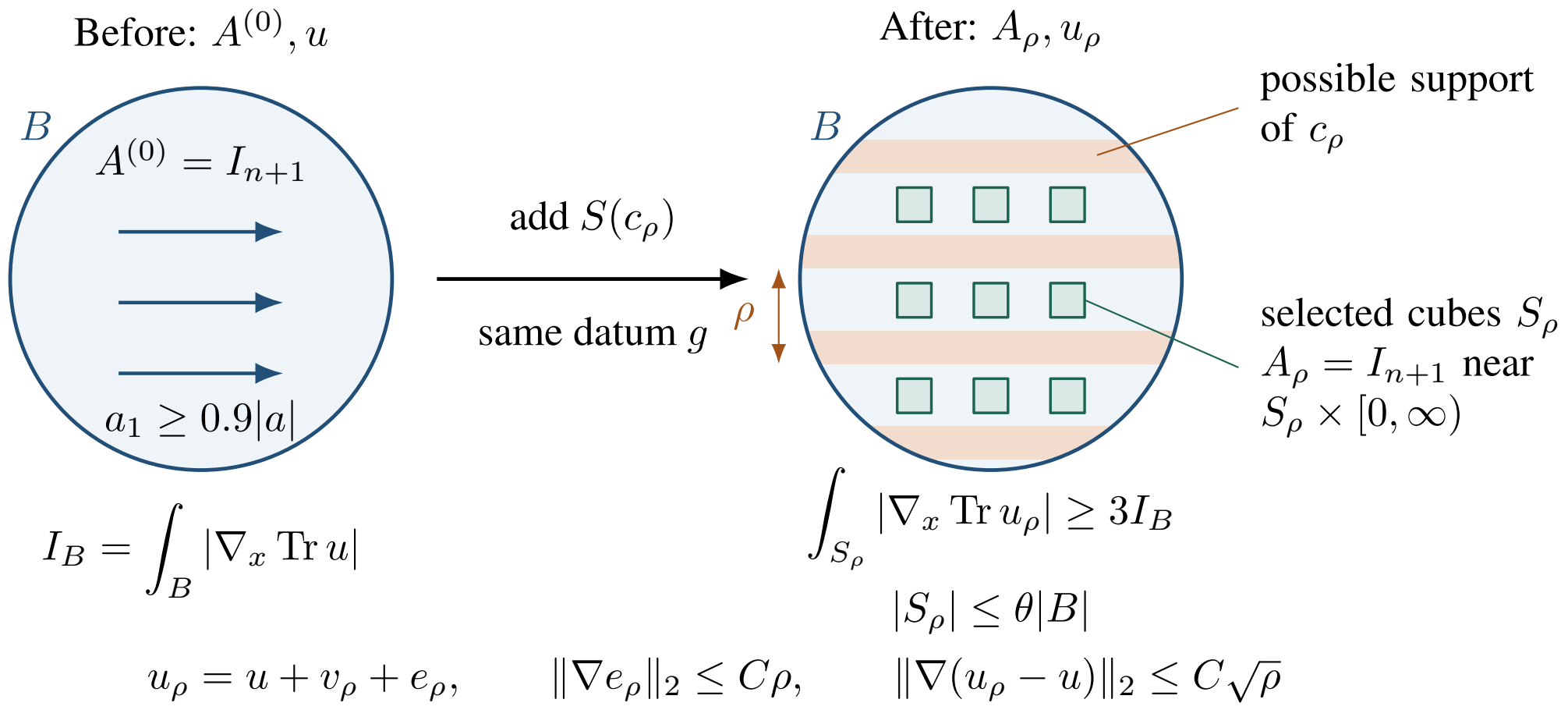

	Figure~\ref{fig:hs-local-perturbation} shows where the cubes are chosen.
	
	\begin{proof}
		The constants in this proof may depend on $u$, the balls, the cutoffs and the
		data fixed in Section~\ref{sec:profiles}, but are uniform in $\rho$.
		Put $\widehat R_\nu=\operatorname{diag}(R_\nu,1)$. Since $R_\nu$ is orthogonal,
		\[
		\begin{gathered}
			\nabla u_\nu=\widehat R_\nu^T(\nabla u)\circ(T_\nu,\mathrm{id}),\quad
			\widehat R_\nu^TS(c)\widehat R_\nu=S(R_\nu^Tc),\quad
			|\det\widehat R_\nu|=1,\\
			\div_y(R_\nu^Tc\circ T_\nu)=(\div_xc)\circ T_\nu,\quad
			\|\nabla u_\nu\|_2=\|\nabla u\|_2.
		\end{gathered}
		\]
		Here the divergence identity follows from
		$\sum_j(R_\nu)_{ij}(R_\nu)_{kj}=\delta_{ik}$. Changing variables in
		the weak equation thus conjugates the coefficient by $\widehat R_\nu$ and
		replaces $g$ by $g\circ T_\nu$. It preserves all the norms in the lemma.
		
		To estimate $c_\rho$ and $v_\rho$, we fix one ball and omit the index $\nu$.
		In the rotated coordinates we still denote the solution by $u$, the coefficient by
		$A^{(0)}$, and the datum by $g$.
		The even reflection of $u$ is harmonic near $\overline B\times\{0\}$,
		since $A^{(0)}=I_{n+1}$ on the column and $g=0$ near $\overline B$.
		In particular, for some $h_0>0$, $m_u<\infty$ and every $(x,t)\in
		\overline B\times[0,h_0]$, we have
		\begin{equation}\label{lp:old-gradient}
			|\nabla_xu(x,t)|\leq m_u,\qquad
			|\nabla_xu(x,t)-a(x)|
			\leq t\sup_{\overline B\times[0,h_0]}|\partial_t\nabla_xu|
			\leq m_ut.
		\end{equation}
		
		In the following formulas, $f,h$ and their derivatives are evaluated at
		$s=x_2/\rho$, $\zeta$ and $\zeta'$ at $\tau=t/\rho$, and $w$ and its
		derivatives at $(s,\tau)$. Differentiating \eqref{eq:local-field} and
		using $f'=h$ gives, for $1\leq j\leq n$,
		\begin{align*}
			\partial_jc_\rho
			&=K\zeta\bigl[\phi_j h e_1+\rho f(s)\bigl(\phi_{2j}e_1-\phi_{1j}e_2\bigr)
			+\delta_{j2}
			(\rho^{-1}\phi h'e_1+h(\phi_2e_1-\phi_1e_2))\bigr],\\
			\partial_tc_\rho
			&=K\zeta'\bigl[\rho^{-1}\phi h e_1+f(s)\bigl(\phi_2e_1-\phi_1e_2\bigr)\bigr].
		\end{align*}
		Here $\delta_{j2}=1$ for $j=2$ and $\delta_{j2}=0$ otherwise.
		With $\phi$ fixed, we choose $\rho$ small enough that
		\begin{equation}\label{lp:scale-conditions}
			\begin{gathered}
				0<\rho<1,\quad 4\rho<h_0,\quad
				\rho\|\nabla\phi\|_\infty\leq1,\quad
				\rho^2\|D^2\phi\|_\infty\leq1,\\
				\sqrt n\rho<\min\{\dist(\overline{B'},\{\phi\ne1\}),
				\dist(\supp\phi,B^c)\}.
			\end{gathered}
		\end{equation}
		By \eqref{eq:local-field} and the derivative formulas above,
		\eqref{lp:scale-conditions} gives
		$\|c_\rho\|_\infty\leq C_*$ and $\|\nabla c_\rho\|_\infty\leq C_*/\rho$.
		In particular, the derivatives of the cutoffs enter the choice of $\rho$,
		while $C_*$ depends only on $n$ and the fixed periodic functions and constants.

		By the choice of $\zeta$, we also have
		\[
		\begin{array}{ll}
			\partial_tc_{\nu,\rho}=0,&\zeta=1,\\
			c_{\nu,\rho}(x,t)=0\ (t\leq\rho\text{ or }t\geq2\rho),
			&\zeta=\eta.
		\end{array}
		\]
		By \eqref{eq:local-correction}, the derivatives of $v_\rho$ are
		\[
		\nabla_xv_\rho=K\phi a_1 w_se_2+\rho Kw\nabla(\phi a_1),
		\qquad (v_\rho)_t=K\phi a_1 w_\tau.
		\]
		For the rest of the proof, we fix $0<\gamma<\pi$ and allow the implicit
		constants to depend on $\gamma$. Applying \eqref{eq:profile-weighted} with
		$E=B$, we obtain, for all fixed nonnegative integers $j,k$,
		\[
		\bigl\|e^{\gamma t/\rho}
		(\partial_s^j\partial_\tau^kw)(x_2/\rho,t/\rho)
		\bigr\|_{L^2(B\times(0,\infty))}
		\leq C_{j,k,B}\sqrt\rho.
		\]
		Applying this estimate to $w,w_s,w_\tau$, and using boundedness of
		$w$, $\phi a_1$ and $\nabla(\phi a_1)$, we deduce that
		\begin{equation}\label{lp:correction-norm}
			\|e^{\gamma t/\rho}\nabla v_\rho\|_2\leq C\sqrt\rho,\qquad
			\|v_\rho\|_2\leq C\rho^{3/2},\qquad
			\|v_\rho\|_\infty\leq C\rho.
		\end{equation}
		Since $2^*=2(n+1)/(n-1)>2$, we also have
		$\|v_\rho\|_{2^*}^{2^*}\leq\|v_\rho\|_\infty^{2^*-2}\|v_\rho\|_2^2<\infty$.
		This and \eqref{lp:correction-norm} give $v_\rho\in Y^{1,2}(\R^{n+1}_+)$.

		To estimate the error in the weak equation, we write $\widetilde u\coloneqq u+v_\rho$ and
		$F_\rho\coloneqq\nabla v_\rho+S(c_\rho)\nabla\widetilde u$. Since $A^{(0)}=I_{n+1}$ on $\supp\nabla v_\rho$, we have
		\[
		(A^{(0)}+S(c_\rho))\nabla\widetilde u-A^{(0)}\nabla u=F_\rho,
		\qquad (F_\rho)_x=\nabla_xv_\rho+c_\rho\,\partial_t\widetilde u,
		\qquad (F_\rho)_{n+1}=(v_\rho)_t-c_\rho\cdot\nabla_x\widetilde u.
		\]
		When we take the divergence, the two terms containing
		$c_\rho\cdot\nabla_x\partial_t\widetilde u$ cancel. Since $\div_xc_\rho=0$, we obtain
		\[
		\div F_\rho=\Delta v_\rho-\partial_tc_\rho\cdot\nabla_x\widetilde u\eqqcolon f_\rho,
		\qquad -e_{n+1}\cdot F_\rho(x,0)
		=c_\rho(x,0)\cdot\nabla_x\widetilde u(x,0)-(v_\rho)_t(x,0)\eqqcolon g_\rho.
		\]
		For $\varphi\in C_c^\infty(\overline{\R^{n+1}_+})$, we subtract the weak
		equation for $u$ and integrate $F_\rho\cdot\nabla \varphi$ by parts to obtain
		\begin{equation}\label{lp:weak-error-identity}
			\begin{aligned}
				\mR_\rho(\varphi)
				&\coloneqq \int_{\R^{n+1}_+}(A^{(0)}+S(c_\rho))\nabla(u+v_\rho)
				\cdot\nabla \varphi\,dx\,dt-\langle g,\Tr \varphi\rangle\\
				&=-\int_{\R^{n+1}_+}f_\rho \varphi\,dx\,dt
				+\int_{\R^n}g_\rho\Tr \varphi\,dx.
			\end{aligned}
		\end{equation}
		The cutoffs are supported in $B$, so there is no term from $\partial B$.
		Since $\varphi(x,t)=0$ for all sufficiently large $t$, the same integration by parts applies
		when $\zeta=1$.
		
		We first estimate $f_\rho$. Substituting
		$w_{ss}+w_{\tau\tau}=h\zeta'$ into the expression for $\Delta v_\rho$
		gives
		\begin{equation}\label{lp:bulk-error-expansion}
			\begin{aligned}
				\Delta v_\rho
				&=(K/\rho)\phi a_1 h\zeta'+2K\partial_2(\phi a_1)w_s
				+\rho K\Delta_x(\phi a_1)w,\\
				f_\rho
				&=2K\partial_2(\phi a_1)w_s+\rho K\Delta_x(\phi a_1)w
				-(K/\rho)\phi h\zeta'(\partial_1u-a_1)\\
				&\quad-Kf\zeta'(\phi_2\partial_1u-\phi_1\partial_2u)
				-\partial_tc_\rho\cdot\nabla_xv_\rho.
			\end{aligned}
		\end{equation}
		The terms with $\rho^{-1}$ from $\Delta v_\rho$ and $-\partial_tc_\rho\cdot\nabla_xu$
		sum to $-(K/\rho)\phi h\zeta'(\partial_1u-a_1)$. By \eqref{lp:old-gradient},
		$|\partial_1u-a_1|\leq2m_u\rho$ on the support of $\zeta'$.
		For the last term in \eqref{lp:bulk-error-expansion}, the derivative
		formulas for $c_\rho$ and $v_\rho$ give
		\[
		\begin{array}{ll}
			|(\partial_tc_\rho)_1|\lesssim\rho^{-1}|\zeta'|,
			&|\partial_1v_\rho|\lesssim\rho|w|,\\
			|(\partial_tc_\rho)_2|\lesssim|\zeta'|,
			&|\partial_2v_\rho|\lesssim|w_s|+\rho|w|.
		\end{array}
		\]
		Multiplying the bounds in each row and adding gives the first inequality below.
		The second follows from \eqref{lp:old-gradient}.
		\[
		|\partial_tc_\rho\cdot\nabla_xv_\rho|
		\lesssim|\zeta'|(|w|+|w_s|),
		\qquad
		\rho^{-1}|\partial_1u-a_1|\leq2m_u
		\quad\text{on }B\times(0,2\rho).
		\]
		By \eqref{eq:profile-weighted}, the weighted $L^2$ norms of
		$2K\partial_2(\phi a_1) w_s$ and $\rho K\Delta_x(\phi a_1) w$ are bounded by
		$C\sqrt\rho$ and $C\rho^{3/2}$, respectively. Every other term in
		$f_\rho$ contains $\zeta'$. These terms vanish for $\zeta=1$. For
		$\zeta=\eta$, the estimates above bound them uniformly on
		their support in $B\times[\rho,2\rho]$. This set has volume at most
		$2|B|\rho$, and $e^{\gamma t/\rho}\leq e^{2\gamma}$ there. In both
		cases it follows that
		\begin{equation}\label{lp:bulk-error-bound}
			\|e^{\gamma t/\rho}f_\rho\|_2\leq C\sqrt\rho.
		\end{equation}
		
		At the boundary, we have $\nabla_xu(x,0)=a(x)$, while
		$w_\tau(s,0)=\zeta(0)h(s)$ gives $(v_\rho)_t(x,0)=K\zeta(0)\phi h a_1$.
		This cancels the leading term of $c_\rho(x,0)\cdot a(x)$, leaving
		\begin{equation}\label{lp:boundary-error-expansion}
			g_\rho=\rho Kf\zeta(0)(\phi_2a_1-\phi_1a_2)
			+c_\rho(x,0)\cdot\nabla_xv_\rho(x,0).
		\end{equation}
		At $t=0$ we have $|\partial_1v_\rho|+|(c_\rho)_2|\lesssim\rho$,
		whereas $|(c_\rho)_1|+|\partial_2v_\rho|\lesssim1$. Thus both products
		in $c_\rho\cdot\nabla_xv_\rho$ are $O(\rho)$. It follows that $|g_\rho|\leq C\rho$ on $\supp\phi$, and $g_\rho$ vanishes
		outside this set. By the dual Sobolev inequality,
		\[
		\|g_\rho\|_{\dot H^{-1/2}}
		\leq C_n\|g_\rho\|_{L^{2n/(n+1)}}
		\lesssim\rho|\supp\phi|^{(n+1)/(2n)}\lesssim\rho.
		\]
		In fact $g_\rho=0$ for $\zeta=\eta$, since $\zeta(0)=0$.
		For the integral over the half-space in \eqref{lp:weak-error-identity}
		we use \eqref{lp:bulk-error-bound} and \eqref{eq:weighted} with $E=B$.
		For the boundary integral we use the trace inequality of Section~\ref{sec:sobolev}. This gives
		\begin{equation}\label{eq:local-residual}
			\begin{aligned}
				|\mR_\rho(\varphi)|
				&\leq\|e^{\gamma t/\rho}f_\rho\|_2
				\|e^{-\gamma t/\rho}\varphi\|_{L^2(B\times(0,\infty))}
				+\|g_\rho\|_{\dot H^{-1/2}}\|\Tr \varphi\|_{\dot H^{1/2}}\\
				&\leq C\rho\|\nabla \varphi\|_2.
			\end{aligned}
		\end{equation}
		The estimate extends \eqref{lp:weak-error-identity} and \eqref{eq:local-residual}
		to $Y^{1,2}(\R^{n+1}_+)$ by density.
		
		We return to the original coordinates. Since the closures of the balls are
		disjoint, the correction on $B_\mu$ contributes nothing on the support of
		the perturbation in $B_\nu$ for $\mu\ne\nu$. Hence
		\begin{equation}\label{eq:finite-cross}
			S(c_{\nu,\rho})\nabla v_{\mu,\rho}=0\quad(\nu\ne\mu).
		\end{equation}
		On the support of $c_{\nu,\rho}$ we have
		$A_\rho=I_{n+1}+S(c_{\nu,\rho})$, while outside all these supports
		$A_\rho=A^{(0)}$. The bound for $\|c_{\nu,\rho}\|_\infty$ therefore
		gives the remaining estimate in \eqref{lp:coefficient-bounds}. For
		each $\nu$, we write
		\[
		\mR_{\nu,\rho}(\varphi)
		=\int_{\R^{n+1}_+}(A^{(0)}+S(c_{\nu,\rho}))
		\nabla(u+v_{\nu,\rho})\cdot\nabla \varphi\,dx\,dt
		-\langle g,\Tr \varphi\rangle.
		\]
		Estimate \eqref{eq:local-residual}, after the change of coordinates, gives
		$|\mR_{\nu,\rho}(\varphi)|\leq C_\nu\rho\|\nabla \varphi\|_2$ for each ball.
		Using the weak equation for $u$ and \eqref{eq:finite-cross}, we obtain
		\begin{align*}
			\mE_\rho(\varphi)
			&\coloneqq \int_{\R^{n+1}_+}A_\rho\nabla(u+v_\rho)\cdot\nabla \varphi\,dx\,dt
			-\langle g,\Tr \varphi\rangle\\
			&=\sum_{\nu=1}^N\mR_{\nu,\rho}(\varphi)
			+\sum_{\nu\ne\mu}\int_{\R^{n+1}_+}S(c_{\nu,\rho})
			\nabla v_{\mu,\rho}\cdot\nabla \varphi\,dx\,dt
			=\sum_{\nu=1}^N\mR_{\nu,\rho}(\varphi).
		\end{align*}
		Summing the bounds over the finitely many balls gives
		$|\mE_\rho(\varphi)|\leq C\rho\|\nabla \varphi\|_2$ and, by
		\eqref{lp:correction-norm}, $\|\nabla v_\rho\|_2\leq C\sqrt\rho$.
		Subtracting from the weak equation for $u_\rho$ the identity for $u+v_\rho$ gives
		\[
		\int_{\R^{n+1}_+}A_\rho\nabla e_\rho\cdot\nabla \varphi\,dx\,dt
		=-\mE_\rho(\varphi)
		\qquad(\varphi\in Y^{1,2}(\R^{n+1}_+)).
		\]
		We take $\varphi=e_\rho\in Y^{1,2}(\R^{n+1}_+)$. Since $\sym A_\rho=I_{n+1}$,
		\begin{equation}\label{eq:local-error}
			\|\nabla e_\rho\|_2^2=-\mE_\rho(e_\rho)
			\leq C\rho\|\nabla e_\rho\|_2,\qquad
			\|\nabla e_\rho\|_2\leq C\rho,\qquad
			\|\nabla(u_\rho-u)\|_2\leq C\sqrt\rho.
		\end{equation}
		For the last inequality we used $u_\rho-u=v_\rho+e_\rho$,
		\eqref{lp:correction-norm}, and $\rho<1$. This proves
		\eqref{lp:local-conclusions}.

		By \eqref{lp:scale-conditions}, every period cube meeting $\supp\phi$ lies in $B$.
		We choose $0<b<1/4$ so that the horizontal $2b$-neighborhood of
		$\widehat Q_0$ lies in $(0,1)^n\cap\{x:x_2\in J\}$. For each full period
		cube $Q_k\subset B'$, we have $\phi=1$ near $Q_k$ and $f(x_2/\rho)=h(x_2/\rho)=0$
		throughout the horizontal $2b\rho$-neighborhood of $S_k$. Hence $A_\rho(x,t)=I_{n+1}$ whenever $\dist(x,S_k)<2b\rho$ and $t>0$. The cubes $S_k$ are disjoint and have measure
		$\theta|Q_k|$, which proves the measure bound in~\textup{(iii)}.
		
		Let $x_k$ be the center of $Q_k$, set $a_k=a(x_k)$ and put
		\[
		u_k(x,t)=a_k\cdot x+\rho K(a_k)_1w(x_2/\rho,t/\rho).
		\]
		For $\dist(x,S_k)<2b\rho$, the equation for $w$ gives $\Delta u_k=0$ for $t>0$ and
		$\partial_t u_k(x,0)=0$. The same holds for $u_\rho$ since $A_\rho=I_{n+1}$
		and $g=0$ there. Thus both functions extend harmonically across $t=0$ by even reflection. Write
		\[
		D_k=\{(x,t):\dist(x,S_k)<b\rho,\quad 0<t<(3+b)\rho\}.
		\]
		If $(x,t)\in D_k$, then $x\in Q_k$, $|x-x_k|\leq\sqrt n\rho$ and
		$t\leq(3+b)\rho<h_0$. Since $\nabla u(x_k,0)=(a_k,0)$ and the even
		extension of $u$ is smooth, we have
		$|\nabla u(x,t)-(a_k,0)|+|a_1(x)-(a_k)_1|\leq C\rho$.
		Moreover, $\phi=1$ near $Q_k$. Hence
		\[
		\begin{aligned}
		\nabla(u+v_\rho-u_k)
		={}&\nabla u-(a_k,0)
		+K\bigl(a_1(x)-(a_k)_1\bigr)(w_s e_2,w_\tau)\\
		&{}+\rho K w(\nabla_x a_1,0).
		\end{aligned}
		\]
		The right-hand side has absolute value at most $C\rho$.
		In the original coordinates, the horizontal projections of the $D_k$
		are disjoint up to sets of measure zero and lie in the balls $B_\nu$.
		Consequently $\sum_k |D_k|\leq(3+b)\rho\sum_\nu |B_\nu|$ and
		\eqref{eq:local-error} gives
		\[
		\sum_k\int_{D_k}|\nabla(u_\rho-u_k)|^2
		\leq 2\|\nabla e_\rho\|_2^2+C\rho^3\leq C\rho^2 .
		\]
		Here and below the sum includes all the balls, with gradients taken in
		their corresponding coordinates.
		
		Fix $t_0\in\{0,3\rho\}$ and let $r=b\rho/2$. Write $v_k=u_\rho-u_k$ and use its even extension if $t_0=0$. For $x\in S_k$ the mean value inequality gives
		\[
		|\nabla v_k(x,t_0)|^2
		\leq C_n r^{-n-1}\int_{B((x,t_0),r)} |\nabla v_k(y,t)|^2\,dy\,dt .
		\]
		For fixed $(y,t)$ the set of $x\in S_k$ for which the ball contains $(y,t)$ lies in the horizontal ball $B(y,r)$ of measure at most $C_n r^n$. Integrating in $x$ and using Fubini's theorem therefore gives
		\[
		\int_{S_k}|\nabla v_k(x,t_0)|^2\,dx
		\leq C_n r^{-1}\int_{D_k}|\nabla v_k(y,t)|^2\,dy\,dt .
		\]
		At $t_0=0$ the integrals on the two sides of the plane agree by reflection, which contributes a factor of two to $C_n$. At $t_0=3\rho$ the balls are contained in $D_k$. Summing and applying the preceding bound, we obtain
		\begin{equation}\label{eq:two-slice-error}
			\sum_k\int_{S_k}|\nabla u_\rho(x,t_0)-\nabla u_k(x,t_0)|^2\,dx
			\leq C\rho,\qquad t_0=0,\ 3\rho .
		\end{equation}
		
		On $S_k$, our choice of $J$ and \eqref{eq:local-alignment} imply
		\[
		|\partial_2u_k(x,0)|\geq(0.9Kd-0.1)|a_k|,\qquad
		|\partial_t u_k(x,3\rho)|\geq0.9Kd|a_k|.
		\]
		Let $G_\rho$ be the union of the period cubes contained in $B'$. The complement $B'\setminus G_\rho$ lies within distance $\sqrt n\rho$ of $\partial B'$, hence its measure tends to zero. If $\omega_a$ is the modulus of continuity of $a$ on $\overline{B'}$, then
		\[
		\left|\sum_{Q_k\subset B'}\rho^n|a_k|-\int_{B'}|a|\right|
		\leq |B'|\omega_a(\sqrt n\rho)+\|a\|_{L^\infty(B')}|B'\setminus G_\rho|
		\longrightarrow 0 .
		\] We integrate the lower bounds for $|\partial_2u_k(x,0)|$ and
		$|\partial_t u_k(x,3\rho)|$ on $S_k$ and sum over $k$.
		Since $|S_k|=\theta\rho^n$, the preceding limit shows that each sum has lower limit at least
		$0.9\theta(0.9Kd-0.1)\sum_\nu I_\nu>4\sum_\nu I_\nu$.
		For $t_0=0$ and $t_0=3\rho$, Cauchy--Schwarz and
		\eqref{eq:two-slice-error} give
		\[
		\sum_k\int_{S_k}|\nabla u_\rho(x,t_0)-\nabla u_k(x,t_0)|\,dx
		\leq |S_\rho|^{1/2}(C\rho)^{1/2}
		\leq \biggl(\theta\sum_\nu|B_\nu|\biggr)^{1/2}(C\rho)^{1/2}.
		\]
		We decrease $\rho_0$ so that both sums with $u_\rho$ in place of $u_k$
		are at least $3\sum_\nu I_\nu$ for $0<\rho<\rho_0$.
		This proves both lower bounds in \eqref{eq:local-conormal-conclusion}.
		The equality of the two interior derivatives follows from $A_\rho=I_{n+1}$
		above $S_\rho$.
	\end{proof}

	\section{Finite iteration}\label{sec:iteration}
	
	We fix $(\zeta,w)$ and $K,C_*,\theta$ as in
	Sections~\ref{sec:profiles}--\ref{sec:local}. After one application of
	Lemma~\ref{lem:local}, the coefficient is the identity on the columns above
	the selected cubes. We shall use these cubes for the next application.
	At each step we first fix the balls and then take the scale small enough to
	at least double the boundary integral and make the $L^2$ norm of the difference of the gradients as small as required.
	On the last selected cubes we also use the estimate at height $3\rho$
	from the same lemma.
	
	\begin{proposition}\label{prop:finite}
		Let $n\geq2$ and let $A^{(0)}\in C^\infty(\overline{\R^{n+1}_+};\R^{(n+1)\times(n+1)})
		\cap L^\infty$ satisfy $\sym A^{(0)}=I_{n+1}$. Let $S$ be a finite union of
		closed cubes with disjoint interiors, and suppose that
		$A^{(0)}=I_{n+1}$ near $S\times[0,\infty)$. Let $g\in C_c^\infty(\R^n)$
		vanish near $S$, and let $u\in Y^{1,2}(\R^{n+1}_+)$ be the weak solution of
		\eqref{eq:variational} with coefficient $A^{(0)}$ and datum $g$. Suppose that
		\[
		I_S\coloneqq \int_S|\nabla_x\Tr u|\,dx>0.
		\]
		For every $H,\delta>0$, there exist finitely many fields $c_\nu$, constructed at scales $\rho_\nu>0$, satisfying
		\[
		\div_x c_\nu=0,\qquad \|c_\nu\|_\infty\leq C_*,
		\qquad \|\nabla c_\nu\|_\infty\leq C_*/\rho_\nu,
		\]
		and
		\[
		\supp_x c_\nu\Subset\operatorname{int}S,\qquad
		\supp_x c_\nu\cap\supp_x c_\mu=\varnothing\quad(\nu\ne\mu),
		\]
		with the following properties. For $\zeta=1$ the fields are independent of $t$.
		For $\zeta=\eta$, $c_\nu$ vanishes unless $\rho_\nu<t<2\rho_\nu$.
		The matrix
		$A'\coloneqq A^{(0)}+\sum_\nu S(c_\nu)$ satisfies
		\[
		\sym A'=I_{n+1}\qquad\text{and}\qquad
		\|A'\|_\infty\leq\max\{\|A^{(0)}\|_\infty,1+C_*\}.
		\]
		If $u'\in Y^{1,2}(\R^{n+1}_+)$ is the weak solution of
		\eqref{eq:variational} with coefficient $A'$ and the same datum $g$, then
		there exist a finite union $S'\Subset\operatorname{int}S$ of closed cubes
		and a scale $\rho>0$ for the last step such that
		\begin{equation}\label{eq:finite-conclusion}
			\begin{gathered}
				\int_{S'}|\nabla_x\Tr u'|\,dx\geq H,\qquad
				\int_{S'}|\partial_tu'(x,3\rho)|\,dx
				=\int_{S'}|(A'\nabla u')_{n+1}(x,3\rho)|\,dx\geq H,\\
				\|\nabla(u'-u)\|_2<\delta,\qquad |S'|\leq |S|.
			\end{gathered}
		\end{equation}
		Moreover, $A'(x,t)=I_{n+1}$ whenever $\dist(x,S')<b\rho/2$ and $t>0$,
		with $b$ as in Lemma~\ref{lem:local}. In particular,
		\[
		\dist(S',\supp_x(A'-I_{n+1}))\geq b\rho/2.
		\]
		Thus $A'$ and $u'$ are smooth up to the boundary near $S'$.
		
		The constant $C_*$ is the one in Lemma~\ref{lem:local}. It is independent
		of $H$, $\delta$, and the number of steps. At each step we may impose any
		finite number of positive upper bounds on the scale, provided we fix
		them before choosing the scale. We describe the associated periodic
		cubes in Section~\ref{sec:carleson}.
	\end{proposition}
	
	\begin{proof}
		We first choose disjoint balls inside $S$ on which the hypotheses of Lemma~\ref{lem:local} hold
		and on whose union the integral of $|\nabla_x\Tr u|$ exceeds $0.9I_S$. Since
		$A^{(0)}=I_{n+1}$ near $S\times[0,\infty)$ and $g$ vanishes near $S$,
		Lemma~\ref{lem:reflection} implies that $a\coloneqq\nabla_x\Tr u$ is
		smooth near $S$. We write $S=\bigcup_{\alpha=1}^M Q_\alpha$ and set
		\[
		E\coloneqq\biggl(\bigcup_{\alpha=1}^M\operatorname{int}Q_\alpha\biggr)
		\cap\{|a|>0\},\qquad
		d\mu\coloneqq\chi_S |a|\,dx.
		\]
		The boundaries of the cubes have Lebesgue measure zero, and $\mu$
		vanishes on $\{a=0\}$. Hence, $\mu(E)=I_S$. Let $x_0\in E$. We
		choose $R_{x_0}\in O(n)$ so that $R_{x_0}e_1=a(x_0)/|a(x_0)|$. By
		continuity, we can take $r_{x_0}>0$ small enough that
		$\overline{B(x_0,r_{x_0})}\Subset\operatorname{int}Q_\alpha$ for some
		$\alpha$ and
		\[
		|a(x)-a(x_0)|\leq\tfrac1{20}|a(x_0)|
		\qquad\text{for }x\in\overline{B(x_0,r_{x_0})}.
		\]
		Since $R_{x_0}^Ta(x_0)=|a(x_0)|e_1$, it follows that, on this ball,
		\[
		\begin{aligned}
			(R_{x_0}^Ta(x))_1
			&\geq\tfrac{19}{20}|a(x_0)|\geq\tfrac{19}{21}|a(x)|>0.9|a(x)|,\\
			|(R_{x_0}^Ta(x))_2|
			&\leq\tfrac1{20}|a(x_0)|\leq\tfrac1{19}|a(x)|<0.1|a(x)|.
		\end{aligned}
		\]
		All smaller closed balls with the same center satisfy these inequalities.
		By Vitali's covering theorem \cite[Theorem~2.2]{Mattila}, there are countably many
		pairwise disjoint such closed balls which cover $E$ up to a Lebesgue null set.
		Since $\mu\ll dx$, their boundaries have $\mu$-measure zero. Taking the interiors
		of finitely many of them, we obtain balls $B_\nu$ and rotations $R_\nu$,
		$1\leq\nu\leq N$, such that
		\begin{equation}\label{eq:finite-vitali}
			\begin{gathered}
				\overline B_\nu\cap\overline B_\mu=\varnothing\quad(\nu\ne\mu),
				\qquad \overline B_\nu\Subset\operatorname{int}Q_{\alpha_\nu},\\
				I_\nu\coloneqq \int_{B_\nu}|a|\,dx>0,\qquad
				\sum_{\nu=1}^NI_\nu>0.9I_S,\\
				(R_\nu^Ta)_1\geq0.9|a|,\qquad |(R_\nu^Ta)_2|\leq0.1|a|
				\quad\text{on }B_\nu.
			\end{gathered}
		\end{equation}
		
		For these balls, Lemma~\ref{lem:local} gives the matrix
		$A_\rho=A^{(0)}+\sum_\nu S(c_{\nu,\rho})$, the corresponding solution
		$u_\rho$ with datum $g$, and the union $S_1=\bigcup_\nu S_{\nu,\rho}$ of
		closed cubes. By \eqref{eq:finite-vitali} and the disjointness of the
		balls, for all sufficiently small $\rho$ we have
		\[
		\begin{gathered}
			\int_{S_1}|\nabla_x\Tr u_\rho|\,dx\geq3\sum_{\nu=1}^NI_\nu>2I_S,
			\qquad |S_1|\leq\theta\sum_\nu|B_\nu|\leq\theta|S|,\\
			S_1\Subset\operatorname{int}S,\qquad
			\|\nabla(u_\rho-u)\|_2\leq C\sqrt\rho,\qquad
			A_\rho=I_{n+1}\text{ near }S_1\times[0,\infty).
		\end{gathered}
		\]
		Each new field has horizontal support compactly contained in its ball
		and vanishes in a horizontal neighborhood of $S_1$. Since $g$ also
		vanishes near $S_1$, we may repeat the same construction with
		$S_1,A_\rho,u_\rho$ in place of $S,A^{(0)},u$.
		
		We now fix an integer $m\geq1$ such that $2^m I_S\geq H$. We write $S^{(k)}$,
		$A^{(k)}$, and $u^{(k)}$ for the set, coefficient, and solution
		obtained at step $k$, starting with $S^{(0)}=S$, $u^{(0)}=u$, and
		$I_0=I_S$. At step $k$ we choose the balls inside the individual cubes of
		$S^{(k-1)}$ as in \eqref{eq:finite-vitali}. This is possible since the
		boundaries of the cubes have zero $|\nabla_x\Tr u^{(k-1)}|\,dx$-measure.
		Once the balls, rotations, and cutoffs are fixed, we denote by $C_k$ and
		$\rho_{0,k}$ the constants in Lemma~\ref{lem:local}. If $k\geq2$, each selected
		cube from the previous step lies inside a period cube $Q$, which we call its
		parent cube. We choose $\rho_k$ so that
		\[
		\begin{gathered}
			\rho_k<\rho_{0,k},\qquad C_k\sqrt{\rho_k}<\delta/m,\\
			\sqrt n\rho_k<\min_\nu\dist(\supp\phi_{k,\nu},B_{k,\nu}^c),\\
			\rho_k\leq\ell(Q)/4
			\quad\text{for every parent cube }Q\quad(k\geq2).
		\end{gathered}
		\]
		Any period cube meeting $\supp\phi_{k,\nu}$ lies inside $B_{k,\nu}$ by
		the third inequality, since its diameter is $\sqrt n\rho_k$. The number
		of upper bounds is finite, so we can also meet the extra restrictions
		in the statement. Applying Lemma~\ref{lem:local} at this scale gives
		\[
		\begin{gathered}
			S^{(k)}\Subset\operatorname{int}S^{(k-1)},\qquad
			|S^{(k)}|\leq\theta|S^{(k-1)}|,\\
			I_k\coloneqq \int_{S^{(k)}}|\nabla_x\Tr u^{(k)}|\,dx>2I_{k-1},
			\qquad\|\nabla(u^{(k)}-u^{(k-1)})\|_2<\delta/m,\\
			\int_{S^{(k)}}|\partial_t u^{(k)}(x,3\rho_k)|\,dx
			= \int_{S^{(k)}}|(A^{(k)}\nabla u^{(k)})_{n+1}(x,3\rho_k)|\,dx
			>2I_{k-1}.
		\end{gathered}
		\]
		Iterating the inequality for $I_k$ and summing the changes in the
		gradient give
		\[
		I_m>2^mI_S\geq H,\qquad
		\|\nabla(u^{(m)}-u)\|_2
		\leq\sum_{k=1}^m\|\nabla(u^{(k)}-u^{(k-1)})\|_2<\delta.
		\]
		
		At the last step the integral of $|\partial_t u^{(m)}(x,3\rho_m)|$ on
		$S^{(m)}$ is greater than $2I_{m-1}\geq 2^m I_S\geq H$.
		The fields introduced at a given step have disjoint horizontal supports.
		Every ball at a later step lies inside a cube near which all these fields vanish. Thus at each point at most one of the fields
		$c_\nu$ from the $m$ steps is nonzero. There $A^{(0)}=I_{n+1}$ and
		$\|c_\nu\|_\infty\leq C_*$. In particular,
		\[
		\|A^{(m)}\|_\infty\leq\max\{\|A^{(0)}\|_\infty,1+C_*\}.
		\]
		If $c_\nu$ was introduced at step $k$, the nesting of the selected
		cubes and Lemma~\ref{lem:local} give
		\[
		\dist(S^{(m)},\supp_x c_\nu)
		\geq\dist(S^{(k)},\supp_x c_\nu)>0.
		\]
		Since there are finitely many fields and
		$S^{(m)}\Subset\operatorname{int}S$, we obtain
		\[
		\dist(S^{(m)},\supp_x(A^{(m)}-I_{n+1}))>0.
		\]
		Here we used that
		$A^{(0)}=I_{n+1}$ on $S\times(0,\infty)$. By Lemma~\ref{lem:local} at the last step,
		$A^{(m)}(x,t)=I_{n+1}$ whenever $\dist(x,S^{(m)})<b\rho_m/2$ and $t>0$. Finally $g=0$ near $S^{(m)}$,
		so even reflection gives smoothness of $u^{(m)}$ up to the boundary there
		by Lemma~\ref{lem:reflection}. Taking $S'=S^{(m)}$, $A'=A^{(m)}$,
		$u'=u^{(m)}$ and $\rho=\rho_m$ proves \eqref{eq:finite-conclusion}.
	\end{proof}
	
	Figure~\ref{fig:hs-packing} shows how the choice of the cubes allows
	the construction to be repeated with disjoint horizontal supports.
	
	\begin{figure}[htbp]
		\centering
		\begin{tikzpicture}[font=\small,>=Latex,line join=round]
			\begin{scope}
				\node at (1.8,4.15) {Step $k$};
				\path[fill=mathBlue!16] (0,0) rectangle (3.6,3.6);
				\path[fill=white] (.61,.61) rectangle (2.99,2.99);
				\draw[mathBlue,line width=.9pt] (0,0) rectangle (3.6,3.6);
				\draw[densely dashed,line width=.9pt] (.7,.7) rectangle (2.9,2.9);
				\node[anchor=south east,mathBlue] at (3.6,3.64) {$Q$};
				\node[fill=white,inner sep=2pt] at (1.8,1.8) {$\widehat Q$};
				\node[align=center] at (1.8,-.38)
				{the new field vanishes\\near $\widehat Q$};
			\end{scope}
			\draw[->,line width=.8pt] (3.77,1.8)--(4.43,1.8);
			\begin{scope}[shift={(4.6,0)}]
				\node at (1.8,4.15) {Step $k+1$};
				\path[fill=mathBlue!16] (0,0) rectangle (3.6,3.6);
				\path[fill=white] (.61,.61) rectangle (2.99,2.99);
				\draw[mathBlue,line width=.9pt] (0,0) rectangle (3.6,3.6);
				\draw[densely dashed,line width=.9pt] (.7,.7) rectangle (2.9,2.9);
				\node[anchor=south east,mathBlue] at (3.6,3.64) {$Q$};
				\node[anchor=north west,fill=white,inner sep=1pt] at (.76,2.85)
				{$\widehat Q$};
				\foreach \xx/\yy in {1.0/1.02,2.0/1.02,1.5/2.0}{
					\path[fill=mathGreen!23,draw=mathGreen,line width=.8pt]
					(\xx,\yy) rectangle ({\xx+.72},{\yy+.72});
					\path[fill=white]
					({\xx+.10},{\yy+.10}) rectangle ({\xx+.62},{\yy+.62});
					\draw[densely dashed,line width=.6pt]
					({\xx+.14},{\yy+.14}) rectangle ({\xx+.58},{\yy+.58});
				}
				\node[anchor=south west,fill=white,inner sep=1pt,text=mathGreen]
				at (2.24,2.06) {$Q'$};
				\node[align=center] at (1.8,-.38)
				{the fields at step $k+1$ lie in $\widehat Q$\\and vanish near each $\widehat Q'$};
			\end{scope}
			\draw[->,line width=.8pt] (8.37,1.8)--(9.73,1.8);
			\node[align=center] at (9.05,2.38) {enlarge\\one $Q'$};
			\begin{scope}[shift={(9.9,0)}]
				\node at (1.8,4.15) {Step $k+2$};
				\path[fill=mathGreen!23] (0,0) rectangle (3.6,3.6);
				\path[fill=white] (.61,.61) rectangle (2.99,2.99);
				\draw[mathGreen,line width=.9pt] (0,0) rectangle (3.6,3.6);
				\draw[densely dashed,line width=.9pt] (.7,.7) rectangle (2.9,2.9);
				\node[anchor=south east,mathGreen] at (3.6,3.64) {$Q'$};
				\node[anchor=north west,fill=white,inner sep=1pt] at (.76,2.85)
				{$\widehat Q'$};
				\foreach \xx/\yy in {1.0/1.02,2.0/1.02,1.5/2.0}{
					\path[fill=mathOrange!24,draw=mathOrange,line width=.8pt]
					(\xx,\yy) rectangle ({\xx+.72},{\yy+.72});
					\path[fill=white]
					({\xx+.10},{\yy+.10}) rectangle ({\xx+.62},{\yy+.62});
					\draw[densely dashed,line width=.6pt]
					({\xx+.14},{\yy+.14}) rectangle ({\xx+.58},{\yy+.58});
				}
				\node[anchor=south west,fill=white,inner sep=1pt,text=mathOrange]
				at (2.24,2.06) {$Q''$};
				\node[align=center] at (1.8,-.38)
				{repeat inside $\widehat Q'$\\with disjoint new supports};
			\end{scope}
			\node at (6.75,-1.24)
			{$\displaystyle |S^{(k)}|\leq\theta^k|S^{(0)}|,
				\qquad I_k\coloneqq\int_{S^{(k)}}|\nabla_x\Tr u^{(k)}|\,dx
				>2^k I_0$};
			\node at (6.75,-2.02)
			{$\displaystyle\supp_xc_\nu\cap\supp_xc_\mu=\varnothing
				\quad(\nu\ne\mu)
				\quad\Longrightarrow\quad
				\|A^{(m)}\|_\infty
				\leq\max\{\|A^{(0)}\|_\infty,1+C_*\}$};
		\end{tikzpicture}
		\caption{Three steps of the construction in the horizontal variables. The colored regions contain the supports of the new fields. Inside each dashed cube the fields of that step vanish, so we may choose the balls for the next step there. The last panel shows one of the cubes $Q'$ of the middle panel at a larger scale. For $k\geq1$, the bounds on $|S^{(k)}|$ and $I_k$ hold for the union of all the selected cubes at step $k$.}
		\label{fig:hs-packing}
	\end{figure}
	
	\FloatBarrier
	
	\section{Carleson estimates}\label{sec:carleson}
	
	Take $\zeta=\eta$. Each field constructed at scale $\rho$ is supported
	between heights $\rho$ and $2\rho$. We estimate the two Carleson
	measures using the cubes $Q$ of side $\rho$ that meet the support of the horizontal cutoff. The contribution of one such cube will be at most $C|Q|$.
	The cubes of the next step inside $Q$ are contained in $\widehat Q$,
	where $|\widehat Q|=\theta |Q|$ and $\theta=|\widehat Q_0|<1$
	(Figure~\ref{fig:hs-packing}). Summing over successive generations will therefore give a geometric series.
	In this section the constants depend only on $n,K,f,\eta$.
	
	\subsection{The periodic cubes and their packing}
	
	For each ball $B_{k,\nu}$ chosen at step $k$ of Proposition~\ref{prop:finite},
	we write $T_{k,\nu}(y)\coloneqq x_{k,\nu}+R_{k,\nu}y$ for the corresponding
	change of coordinates. We denote by $\phi_{k,\nu}$ the
	cutoff of Section~\ref{sec:local} composed with $T_{k,\nu}^{-1}$. For
	$b\in\mathbb{Z}^n$, we set
	\[
	\begin{gathered}
		Q_{k,\nu,b}\coloneqq T_{k,\nu}(\rho_k(b+[0,1]^n)),\qquad
		\widehat Q_{k,\nu,b}\coloneqq T_{k,\nu}(\rho_k(b+\widehat Q_0)),\\
		\PP_{k,\nu}\coloneqq\{Q_{k,\nu,b}:b\in\mathbb{Z}^n,
		\ Q_{k,\nu,b}\cap\supp\phi_{k,\nu}\ne\varnothing\},
		\qquad\PP_k\coloneqq\bigcup_\nu\PP_{k,\nu}.
	\end{gathered}
	\]
	If $Q\in\PP_{k,\nu}$, then $Q\Subset B_{k,\nu}$, since its diameter
	$\sqrt n\rho_k$ is smaller than the distance from $\supp\phi_{k,\nu}$ to
	$B_{k,\nu}^c$. We write $\widehat Q$ for the image of the smaller cube under the
	same change of coordinates. By Lemma~\ref{lem:local}, the selected
	cubes come from period cubes contained in the smaller ball on which
	$\phi_{k,\nu}=1$. Thus $S^{(k)}$ in Proposition~\ref{prop:finite} is a union
	of cubes $\widehat Q$ with $Q\in\PP_k$.
	
	We choose the Vitali balls for step $k+1$ in the interiors of the individual
	cubes $\widehat Q$ making up $S^{(k)}$. The boundaries of these cubes have zero
	$|\nabla_x\Tr u^{(k)}|\,dx$-measure, so removing them does not change
	the integral of $|\nabla_x\Tr u^{(k)}|$ over $S^{(k)}$. Once we have chosen
	the balls and their cutoffs, we take $\rho_{k+1}$ small enough that, for each
	ball and its corresponding cube $Q$,
	\begin{equation}\label{eq:cell-restrictions}
		\begin{gathered}
			\overline B_{k+1,\nu}\Subset\operatorname{int}\widehat Q,\qquad
			\sqrt n\rho_{k+1}<\dist(\supp\phi_{k+1,\nu},B_{k+1,\nu}^c),
			\qquad\rho_{k+1}\leq\ell(Q)/4,\\
			Q'\in\PP_{k+1,\nu}\quad\Longrightarrow\quad
			Q'\Subset B_{k+1,\nu}\Subset\operatorname{int}\widehat Q\Subset Q.
		\end{gathered}
	\end{equation}
	This is possible since at each step there are only finitely many positive
	upper bounds on $\rho_{k+1}$.
	
	We now carry out this construction in pairwise disjoint initial balls
	$B_j\subset\R^n$, with $m_j<\infty$ steps in $B_j$. We add the index $j$ to
	the preceding notation and set
	\[
	\PP_{j,k}\coloneqq\bigcup_\nu\PP_{j,k,\nu},\qquad
	\PP\coloneqq\bigcup_j\bigcup_{k=1}^{m_j}\PP_{j,k}.
	\]
	For $Q\in\PP$, we denote its center by $x_Q$ and its side length by $\ell(Q)$.
	The dilate $\lambda Q$, for $\lambda>0$, has the same center and orientation as
	$Q$ and satisfies $|\lambda Q|=\lambda^n|Q|$. The orientation is the one
	determined by the ball in which $Q$ was chosen. By construction,
	\[
	\widehat Q\Subset\operatorname{int}Q
	\qquad\text{and}\qquad |\widehat Q|=\theta|Q|.
	\]
	It follows from \eqref{eq:cell-restrictions} that every $Q'\in\PP_{j,k+1}$
	has a unique parent $Q\in\PP_{j,k}$. In fact, $Q'\Subset\operatorname{int}\widehat Q$
	and $\ell(Q')\leq\ell(Q)/4$, and the cubes in $\PP_{j,k}$ have disjoint interiors.
	For $Q\in\PP_{j,k}$, let $\mC(Q)$ consist of the cubes of $\PP_{j,k+1}$
	contained in $Q$, with $\mC(Q)=\varnothing$ if $k=m_j$.
	Since the children of $Q$ lie in $\widehat Q$ and have disjoint interiors,
	\begin{equation}\label{eq:children}
		\sum_{Q'\in\mC(Q)}|Q'|\leq|\widehat Q|=\theta|Q|.
	\end{equation}
	In particular, any two cubes in $\PP$ are nested or have disjoint interiors.
	Indeed, cubes from distinct initial balls are disjoint. Suppose that $Q\in\PP_{j,k}$ and $Q'\in\PP_{j,l}$, with $k<l$, have intersecting interiors.
	Let $P$ be the ancestor of $Q'$ in $\PP_{j,k}$. Then $\operatorname{int}P\cap\operatorname{int}Q\ne\varnothing$, so $P=Q$ and $Q'\subset Q$.
	
	\begin{lemma}\label{lem:packing}
		For every cube $Q_*\subset\R^n$, we have
		\begin{equation}\label{eq:packing}
			\sum_{\substack{Q\in\PP\\Q\subset Q_*}}|Q|
			\leq\frac{|Q_*|}{1-\theta}.
		\end{equation}
	\end{lemma}
	\begin{proof}
		For $Q\in\PP$, we set $\mC^0(Q)\coloneqq\{Q\}$ and denote by
		$\mC^{l+1}(Q)$ the family of children of all cubes in $\mC^l(Q)$.
		Iterating \eqref{eq:children} and then summing over $l$, we obtain
		\[
		\sum_{Q'\in\mC^l(Q)}|Q'|\leq\theta^l|Q|,
		\qquad
		\sum_{\substack{Q'\in\PP\\Q'\subset Q}}|Q'|
		\leq\sum_{l\geq0}\theta^l|Q|=\frac{|Q|}{1-\theta}.
		\]
		Every cube of $\PP$ contained in $Q_*$ lies in a maximal such cube, since it
		has only finitely many ancestors. The maximal cubes have disjoint interiors.
		Summing the last inequality over them gives \eqref{eq:packing}.
	\end{proof}
	
	\FloatBarrier
	\subsection{The derivative and discrepancy Carleson measures}
	
	\begin{proposition}\label{prop:dkp}\label{prop:discrepancy}
		Let $\{c_\iota\}_{\iota\in\mathcal I}$ be the fields of Lemma~\ref{lem:local}
		obtained with $\zeta=\eta$. Suppose that
		\[
		\supp_x c_\iota\cap\supp_x c_{\iota'}=\varnothing
		\qquad\text{for }\iota\ne\iota',
		\]
		that the sum $A\coloneqq I_{n+1}+\sum_{\iota\in\mathcal I}S(c_\iota)$ is
		locally finite in $\R^{n+1}_+$, and that \eqref{eq:packing} holds for the family $\PP$ of period cubes
		that meet the supports of the horizontal cutoffs.
		Then
		\begin{equation}\label{eq:dkp-bound}
			\|\mu_A\|_{\mC}\leq\frac{C}{1-\theta}
		\end{equation}
		and
		\begin{equation}\label{eq:discrepancy-bound}
			\|\nu_A\|_{\mC}\leq\frac{C}{1-\theta}<\infty.
		\end{equation}
		The constant $C$ depends only on $n,K,f,\eta$. In particular, it is uniform
		in the number of initial balls and the number of steps in each ball.
	\end{proposition}
	\begin{proof}
		For each field $c_\iota$, let $\rho_\iota$ be its scale,
		$\phi_\iota$ its horizontal cutoff, and $\PP_\iota$ the family
		of period cubes meeting $\supp\phi_\iota$. Lemma~\ref{lem:local} gives
		\[
		|c_\iota|\leq C_*,\qquad
		|\nabla c_\iota|\leq C_*/\rho_\iota,\qquad
		\supp\phi_\iota\subset\bigcup_{Q\in\PP_\iota}Q.
		\]
		Both $c_\iota$ and $\nabla c_\iota$ vanish outside
		$\supp\phi_\iota\times[\rho_\iota,2\rho_\iota]$. Local finiteness
		gives $A\in C^\infty(\R^{n+1}_+)$. The horizontal supports of the fields are pairwise disjoint, so
		\[
		\begin{aligned}
			|A(y,s)-I_{n+1}|^2
			&\leq C\sum_{\iota\in\mathcal I}
			\chi_{\supp\phi_\iota\times[\rho_\iota,2\rho_\iota]}(y,s),\\
			|\nabla A(y,s)|^2
			&\leq C\sum_{\iota\in\mathcal I}\rho_\iota^{-2}
			\chi_{\supp\phi_\iota\times[\rho_\iota,2\rho_\iota]}(y,s).
		\end{aligned}
		\]
		To take the essential suprema in \eqref{eq:dkp-definition} and
		\eqref{eq:discrepancy-definition}, we set
		\[
		B_{x,t}\coloneqq B((x,t),t/2),\qquad
		C_0\coloneqq6,\qquad a_0\coloneqq 2/3,\qquad b_0\coloneqq4.
		\]
		Suppose that $B_{x,t}$ contains a point $(y,s)$ of
		$Q\times[\ell(Q),2\ell(Q)]$. Since $t/2<s<3t/2$, we have
		$a_0\ell(Q)<t<b_0\ell(Q)$ and $|x-y|<t/2<2\ell(Q)$. In coordinates
		parallel to the sides of $Q$, each coordinate of $x-x_Q$ is therefore bounded
		in absolute value by $\ell(Q)/2+|x-y|<(1/2+2)\ell(Q)$. Thus $x\in C_0Q$.
		If $\nabla A(y,s)\ne0$, only one field contributes at $(y,s)$, and
		there is a cube $Q\in\PP_\iota$ containing $y$ with
		$|\nabla A(y,s)|\leq C/\ell(Q)$. The same assertion holds for
		$A(y,s)-I_{n+1}$ with bound $C$. Taking the essential supremum
		over $B_{x,t}$ therefore gives
		\begin{equation}\label{eq:dkp-majorant}
			d\mu_A(x,t)\leq C\sum_{Q\in\PP}
			\frac{t}{\ell(Q)^2}\chi_{C_0Q}(x)
			\chi_{(a_0\ell(Q),b_0\ell(Q))}(t)\,dx\,dt
		\end{equation}
		and
		\begin{equation}\label{eq:discrepancy-majorant}
			d\nu_A(x,t)\leq C\sum_{Q\in\PP}
			\frac1t\chi_{C_0Q}(x)
			\chi_{(a_0\ell(Q),b_0\ell(Q))}(t)\,dx\,dt.
		\end{equation}
		For each $Q$ the corresponding integrals are
		\[
		\begin{aligned}
			|C_0Q|\int_{a_0\ell(Q)}^{b_0\ell(Q)}\frac{t}{\ell(Q)^2}\,dt
			&=\frac{C_0^n}{2}(b_0^2-a_0^2)|Q|,\\
			|C_0Q|\int_{a_0\ell(Q)}^{b_0\ell(Q)}\frac{dt}{t}
			&=C_0^n\log(b_0/a_0)|Q|.
		\end{aligned}
		\]
		
		Let $Q_*$ be an axis-parallel cube with center $x_{Q_*}$ and side length $r$.
		The term for $Q$ in either \eqref{eq:dkp-majorant} or \eqref{eq:discrepancy-majorant} vanishes on $Q_*\times(0,r)$ unless
		$C_0Q\cap Q_*\ne\varnothing$ and $\ell(Q)<r/a_0$. For
		$y\in C_0Q\cap Q_*$ and $z\in Q$, the diameters of $Q$ and $C_0Q$ give
		\begin{equation}\label{eq:carleson-center}
			|z-x_{Q_*}|_\infty\leq |z-x_Q|+|x_Q-y|+|y-x_{Q_*}|_\infty
			\leq\tfrac{\sqrt n(1+C_0)}2\ell(Q)+\tfrac r2.
		\end{equation}
		Thus, every such cube $Q$ is contained in the cube
		\[
		\widetilde Q_*\coloneqq\left(1+\frac{\sqrt n(1+C_0)}{a_0}\right)Q_*,
		\]
		independently of the orientation of $Q$. Integrating \eqref{eq:dkp-majorant} and \eqref{eq:discrepancy-majorant} over
		$Q_*\times(0,r)$, we obtain by Tonelli's theorem and \eqref{eq:packing}
		\[
		\begin{aligned}
			\mu_A(Q_*\times(0,r))+\nu_A(Q_*\times(0,r))
			&\leq C\sum_{\substack{Q\in\PP\\Q\subset\widetilde Q_*}}|Q|\\
			&\leq\frac{C}{1-\theta}|\widetilde Q_*|
			\leq\frac{C}{1-\theta}|Q_*|.
		\end{aligned}
		\]
		The concentric cube of side $2r$ contains the ball of radius $r$, so the
		same bound gives \eqref{eq:dkp-bound} and \eqref{eq:discrepancy-bound}.
	\end{proof}

	\section{Proof of Theorem~\ref{thm:main}}\label{sec:global}
	
	In each ball we shall apply Proposition~\ref{prop:finite} to the solution obtained from the preceding balls, always with the same Neumann datum. The first point is to show that this solution has a nonconstant trace in the new ball.
	
	\subsection{A nonvanishing lemma}
	\begin{lemma}\label{lem:seed}
		Let $A$ be a bounded real uniformly elliptic matrix in $\R^{n+1}_+$, let
		$0\ne g\in C_c^\infty(\R^n;\R)$, and let $u\in Y^{1,2}(\R^{n+1}_+)$
		solve \eqref{eq:variational} with datum $g$. Suppose that $u$ is harmonic in
		$\R^n\times(0,h)$ for some $h>0$. Let $Q\subset\R^n$ be a closed cube of
		positive side length, disjoint from $\supp g$. Suppose that for some $r_0>0$
		we have $A(x,t)=I_{n+1}$ whenever $\dist(x,Q)<r_0$ and $0<t<r_0$. Then
		\[
		\int_Q|\nabla_x\Tr u|\,dx>0.
		\]
	\end{lemma}
	\begin{proof}
		Since $Q$ is disjoint from $\supp g$ and $A=I_{n+1}$ near $Q\times\{0\}$ in
		the upper half-space, Lemma~\ref{lem:reflection} gives an even harmonic extension
		of $u$ across $\operatorname{int}Q$. If the integral in the statement were zero,
		then $\nabla_x\Tr u=0$ on $\operatorname{int}Q$ by continuity. Hence, $\Tr u=C$
		there for some constant $C$, and $\partial_tu=0$ there by evenness.
		We fix a ball $B\Subset\operatorname{int}Q$ and choose $0<r<h$ so that the
		extension is harmonic in $B\times(-r,r)$. By Lemma~\ref{lem:cauchy} applied
		to $u-C$, we have $u=C$ throughout $\R^n\times(0,h)$.
		
		We take $\varphi\in C_c^\infty((-h,h))$ with $\varphi(0)=1$ and set
		$v(x,t)\coloneqq g(x)\varphi(t)$. Then $v$ is a smooth test function with
		compact support and $\Tr v=g$. Since $\nabla u=0$ in the slab and $v=0$
		for $t\geq h$, testing \eqref{eq:variational} with $v$ gives
		\[
		0=\int_{\R^{n+1}_+}A\nabla u\cdot\nabla v\,dx\,dt
		=\int_{\R^n}g^2\,dx,
		\]
		contrary to $g\ne0$.
	\end{proof}

	\subsection{Construction of the coefficients}
	
	\begin{proposition}\label{prop:construction}
		Let $n\geq2$. There exist bounded real vector fields $c^{(i)}:\R^{n+1}_+\to\R^n$, $i=1,2$, with $\div_xc^{(i)}(\cdot,t)=0$ for every $t>0$, such that the matrices $A^{(i)}\coloneqq I_{n+1}+S(c^{(i)})$ satisfy the following conditions.
		\begin{enumerate}[label=\textup{(\roman*)}]
			\item The field $c^{(1)}$ is independent of $t$, belongs to $C^\infty(\R^n;\R^n)$, and $|\supp c^{(1)}|<\infty$.
			\item The field $c^{(2)}$ belongs to $C^\infty(\R^{n+1}_+;\R^n)$, vanishes outside $B(0,1)\times(0,1)$, and satisfies
			\[
			\|\mu_{A^{(2)}}\|_{\mC}+\|\nu_{A^{(2)}}\|_{\mC}<\infty.
			\]
			For every $x_0\ne0$, there is $r>0$ such that $A^{(2)}=I_{n+1}$ on $B(x_0,r)\times(0,r)$.
		\end{enumerate}
		
		There exists a single nonzero datum $g\in C_c^\infty(\R^n)$ with $\int_{\R^n}g=0$ for which both variational solutions $u^{(i)}\in Y^{1,2}(\R^{n+1}_+)$ of \eqref{eq:variational} satisfy the estimates below. For each $i$ there are pairwise disjoint compact sets $S_j^{(i)}\subset\R^n$, $j\geq1$, such that
		\begin{equation}\label{eq:main-selected}
			E_i\coloneqq\bigcup_{j\geq1}S_j^{(i)},\qquad |E_i|<\infty,\qquad
			\int_{S_j^{(i)}}|\nabla_x\Tr u^{(i)}|\,dx\geq1.
		\end{equation}
		We also have $E_2\subset B(0,1)$. For every $i$ and $j$, there is a scale
		$\rho_j^{(i)}>0$ such that
		\[
		\dist(x,S_j^{(i)})< b\rho_j^{(i)}/2
		\quad\Longrightarrow\quad g(x)=0
		\quad\text{and}\quad A^{(i)}(x,t)=I_{n+1}\quad\text{for all }t>0.
		\]
		With $t_j^{(i)}=3\rho_j^{(i)}$ we have
		\begin{equation}\label{eq:conormal-selected}
			\int_{S_j^{(i)}}|(A^{(i)}\nabla u^{(i)})_{n+1}(x,t_j^{(i)})|\,dx
			=\int_{S_j^{(i)}}|\partial_tu^{(i)}(x,t_j^{(i)})|\,dx\geq1.
		\end{equation}
		The function $(A^{(i)}\nabla u^{(i)})_{n+1}$ is continuous in a neighborhood of each $S_j^{(i)}\times\{t_j^{(i)}\}$.
	\end{proposition}
	\begin{proof}
		We fix a nonnegative $0\ne\varphi\in C_c^\infty(B(0,1/8))$ and set
		\begin{equation}\label{eq:fixed-datum}
			g(x)\coloneqq\varphi(x-3e_1)-\varphi(x+3e_1).
		\end{equation}
		The translates have the same integral and disjoint supports, so $\int g=0$ and $\|g\|_\infty=\|\varphi\|_\infty>0$. Let $B=B(0,4)$. The function $g/(|B|\|g\|_\infty)$ is supported in $B$, has integral zero and has $L^\infty$ norm $|B|^{-1}$. It is therefore an atom of the real Hardy space $H^1(\R^n)$. Since $g$ is bounded and compactly supported, it also belongs to $L^p(\R^n)$ for every $1\leq p<\infty$. By \eqref{eq:negative},
		\[
		\|g\|_{\dot H^{-1/2}}\leq C_n\|g\|_{L^{2n/(n+1)}}<\infty.
		\]
		For the $t$-independent and DKP matrices, respectively, we use the balls
		\begin{equation}\label{eq:initial-balls}
			B_j^\infty\coloneqq B(4je_1,2^{-j-10}),\qquad
			B_j^0\coloneqq B(2^{-j}e_1,2^{-3j-10}),\qquad j\geq1.
		\end{equation}
		The closed balls avoid $\supp g$ and are pairwise disjoint within each family. Indeed, for $j<k$,
		\[
		\begin{gathered}
			|4je_1-4ke_1|\geq4>2^{-j-10}+2^{-k-10},\\
			|2^{-j}e_1-2^{-k}e_1|\geq2^{-j-1}
			>2^{-3j-10}+2^{-3k-10}.
		\end{gathered}
		\]
		Also, if $\omega_n=|B(0,1)|$, then
		\[
		\sum_j|B_j^\infty|=\omega_n2^{-10n}\sum_j2^{-nj}<\infty,
		\qquad
		\sum_j|B_j^0|=\omega_n2^{-10n}\sum_j2^{-3nj}<\infty.
		\]
		The family $\{B_j^\infty\}_j$ is locally finite since its centers tend to infinity. The balls $B_j^0$ lie in $B(0,1)$ and accumulate only at the origin, since $\sup_{x\in\overline B_j^0}|x|\leq2^{-j}+2^{-3j-10}\to0$.
		
		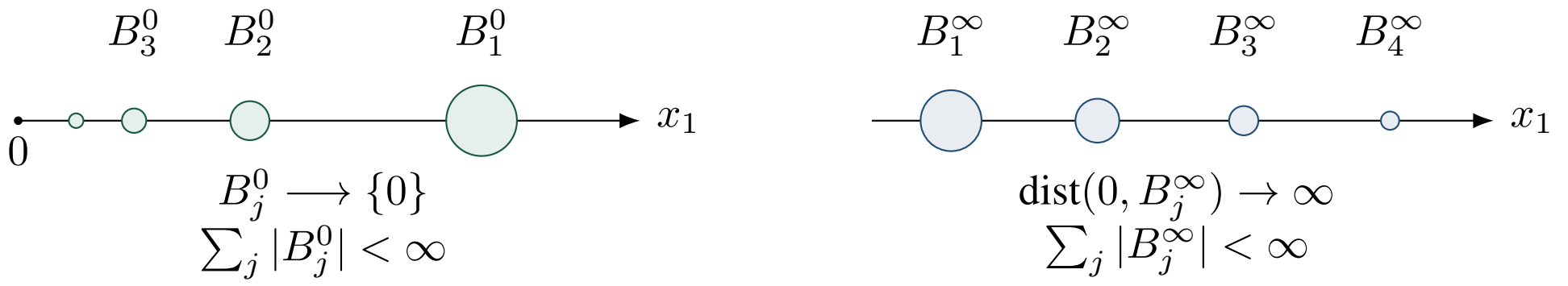
\begin{figure}[htbp]
			\centering
			\begin{tikzpicture}[font=\small,>=Latex]
				\begin{scope}
					\draw[->] (0,0)--(5.1,0) node[right] {$x_1$};
					\fill (0,0) circle (1pt) node[below] {$0$};
					\foreach \xx/\rr/\lab in {3.8/.29/1,1.9/.16/2,.95/.10/3,.475/.06/4}{
						\draw[mathGreen!70!black,fill=mathGreen!12] (\xx,0) circle (\rr);
						\ifnum\lab<4 \node[above] at (\xx,.4) {$B_{\lab}^0$};\fi
					}
					\node[align=center] at (2.5,-.85)
					{$B_j^0\longrightarrow\{0\}$\\$\sum_j|B_j^0|<\infty$};
				\end{scope}
				\begin{scope}[shift={(7.0,0)}]
					\draw[->] (0,0)--(5.1,0) node[right] {$x_1$};
					\foreach \xx/\rr/\lab in {.65/.25/1,1.85/.18/2,3.05/.12/3,4.25/.075/4}{
						\draw[mathBlue,fill=mathBlue!10] (\xx,0) circle (\rr);
						\node[above] at (\xx,.4) {$B_{\lab}^\infty$};
					}
					\node[align=center] at (2.5,-.85)
					{$\dist(0,B_j^\infty)\to\infty$\\$\sum_j|B_j^\infty|<\infty$};
				\end{scope}
			\end{tikzpicture}
			\caption{The balls in \eqref{eq:initial-balls}, with their radii enlarged. In the DKP example they accumulate at the origin (left). For the $t$-independent coefficient they tend to infinity (right).}
			\label{fig:hs-balls}
		\end{figure}
		
		\FloatBarrier
		We pair each family of balls with the corresponding function of $t$,
		\[
		(\zeta,B_j)=(1,B_j^\infty),\qquad
		(\zeta,B_j)=(\eta,B_j^0),
		\]
		and use the corresponding $w$ from \eqref{eq:profiles}. We fix one of these choices and write $B_j=B(b_j,s_j)$. The periodic functions and the constant $K$ are fixed as in Section~\ref{sec:profiles}. At each application of the local perturbation in $B_j$, we choose a scale $\rho$ such that
		\begin{equation}\label{eq:global-height-choice}
			2\rho\leq s_j,\qquad
			Q'\Subset\operatorname{int}\widehat Q,\qquad
			\ell(Q')\leq\ell(Q)/4
			\quad\text{for every child }Q'\text{ of }Q.
		\end{equation}
		To see that we can impose these conditions, let $B$ be one of the balls in Lemma~\ref{lem:local} and $\phi$ its cutoff, with $\supp\phi\Subset B$. After the first step, $\overline B$ is contained in the interior of the cube $\widehat Q$ selected in the parent period cube $Q$. Once we have fixed the balls, rotations, and cutoffs, Lemma~\ref{lem:local} gives an upper bound $\rho_{\rm loc}>0$ for this scale. We take
		\[
		0<\rho<\min\left\{\rho_{\rm loc},\frac{s_j}{2},
		\frac{\ell(Q)}4,\frac{\dist(\supp\phi,B^c)}{2\sqrt n}\right\},
		\]
		omitting $\ell(Q)/4$ at the first step. If a period cube $Q'$ of side length $\rho$ meets $\supp\phi$, then
		\[
		\dist(Q',B^c)\geq\dist(\supp\phi,B^c)-\sqrt n\rho>0.
		\]
		Thus $Q'\Subset B$. The upper bounds for $\rho$ are positive and there are only finitely many balls at this step, so one choice of $\rho$ suffices.
		
		We start with $A_0\coloneqq I_{n+1}$ and let $u_0\in Y^{1,2}(\R^{n+1}_+)$ be the variational solution with datum $g$. At step $j$ we choose a coefficient $A_j$ and its variational solution $u_j$ with the same datum $g$. We obtain a finite union of closed cubes $S_j\Subset B_j$ on which we estimate $\nabla_x\Tr u_j$ and $(A_j\nabla u_j)_{n+1}(\cdot,t_j)$ for some $t_j>0$. To compare these estimates at later steps, we will choose $r_j,C_j>0$.
		Suppose that we have made these choices for every $i<j$. We fix a closed cube $Q_j\Subset B_j$ of positive side length. Since the
		balls $B_i$ are disjoint and each of the fields already chosen is supported
		over some $B_i$ with $i<j$, we have $A_{j-1}=I_{n+1}$ on $B_j\times(0,\infty)$.
		Moreover, $g=0$ on $B_j$. If $\zeta=1$, the sum of these fields is independent
		of $t$ and has zero horizontal divergence, so $u_{j-1}$ is harmonic in
		$\R^{n+1}_+$ by \eqref{eq:operator}. If $\zeta=\eta$, only finitely many
		fields occur in $A_{j-1}$. We take $h>0$ smaller than the least of their scales,
		with $h=1$ when $j=1$. Each field constructed at scale $\rho$
		vanishes for $t\leq\rho$, so $A_{j-1}=I_{n+1}$ in $\R^n\times(0,h)$ and
		$u_{j-1}$ is harmonic there. In either case, Lemma~\ref{lem:seed} gives
		\[
		\int_{Q_j}|\nabla_x\Tr u_{j-1}|\,dx>0.
		\]
		For the error in the $L^2$ norm of the gradient at step $j$, we take
		\begin{equation}\label{eq:budget}
			\ve_j\coloneqq\frac{2^{-j-1}}{1+\max_{i<j}C_i|S_i|},
			\qquad\max\varnothing=0.
		\end{equation}
		Here the denominator involves only quantities chosen at earlier steps. In particular, $C_i|S_i|\ve_j\leq2^{-j-1}$ for each $i<j$. We choose $r_j$ and $C_j$ after fixing the coefficient in $B_j$.
		
		We apply Proposition~\ref{prop:finite} with
		\[
		(A^{(0)},u,S,H,\delta)
		=(A_{j-1},u_{j-1},Q_j,4,\ve_j).
		\]
		At each step we choose the scale to satisfy \eqref{eq:global-height-choice} and $b\rho<\dist(Q_j,B_j^c)$. Denote by $A_j$, $u_j$, and $S_j$ the coefficient, solution, and union of selected cubes given by the proposition. Write $\rho_j$ for the scale of the last step and set $t_j=3\rho_j$. We obtain
		\begin{equation}\label{eq:stage-conclusions}
			\begin{gathered}
				\|\nabla(u_j-u_{j-1})\|_2<\ve_j,\qquad
				\int_{S_j}|\nabla_x\Tr u_j|\,dx\geq4,\\
				\int_{S_j}|(A_j\nabla u_j)_{n+1}(x,t_j)|\,dx
				=\int_{S_j}|\partial_tu_j(x,t_j)|\,dx\geq4 .
			\end{gathered}
		\end{equation}
		Write $c_{j,1},\ldots,c_{j,N_j}$ for all the fields added in $B_j$. We have
		\[
		A_j=A_{j-1}+\sum_{\nu=1}^{N_j}S(c_{j,\nu}).
		\]
		By Proposition~\ref{prop:finite}, the fields in this sum have disjoint horizontal supports. Since the initial balls are also disjoint, $c_j\coloneqq\sum_{i\leq j}\sum_{\nu=1}^{N_i}c_{i,\nu}$ satisfies
		\[
		A_j=I_{n+1}+S(c_j),\qquad
		\|c_j\|_\infty\leq C_*,\qquad\div_xc_j=0.
		\]
		In particular, $A_j$ is smooth up to the boundary, and $\sym A_j=I_{n+1}$. Testing with $u_j$ and applying the trace inequality, we obtain
		\[
		\|\nabla u_j\|_2^2=\langle g,\Tr u_j\rangle
		\leq\|g\|_{\dot H^{-1/2}}\|\nabla u_j\|_2.
		\]
		
		Set $r_j=b\rho_j/8$. Since $S_j\subset Q_j$ and $b\rho_j<\dist(Q_j,B_j^c)$,
		every $x$ with $\dist(x,S_j)<4r_j$ belongs to $B_j$. All subsequent
		perturbations are supported over balls disjoint from $B_j$. Proposition~\ref{prop:finite} gives
		\begin{equation}\label{eq:fixed-neighborhood}
			A_k(x,t)=A_j(x,t)=I_{n+1},\qquad g(x)=0
			\quad\text{if }\dist(x,S_j)<4r_j,\quad t>0,\quad k\geq j.
		\end{equation}
		For $x\in S_j$, the even extension of $u_k-u_j$ is harmonic in $B((x,0),2r_j)$, by Lemma~\ref{lem:reflection}, and is harmonic in $B((x,t_j),2r_j)$. The mean value inequality for its squared gradient at the centers of these balls gives
		\[
		\sup_{x\in S_j}\max\{|\nabla_x\Tr(u_k-u_j)(x)|,|\partial_t(u_k-u_j)(x,t_j)|\}
		\leq C_j\|\nabla(u_k-u_j)\|_2,
		\]
		where $C_j=\max\{1,C_nr_j^{-(n+1)/2}\}$. We can now use $C_j$ in \eqref{eq:budget} at every subsequent step.
		
		We define the limiting field and its coefficient by
		\[
		c(x,t)\coloneqq\sum_{i\geq1}\sum_{\nu=1}^{N_i}c_{i,\nu}(x,t),
		\qquad A=I_{n+1}+S(c).
		\]
		Each $x\in\R^n$ belongs to at most one initial ball $B_i$, and only finitely many fields were added in $B_i$. Thus $c_j(x,t)=c(x,t)$ for all sufficiently large $j$, and likewise $A_j(x,t)=A(x,t)$. In particular,
		\[
		|c|\leq C_*,\qquad A\xi\cdot\xi=|\xi|^2,\qquad
		|A\xi|\leq(1+C_*)|\xi|.
		\]
		For fixed $t\geq0$ and $\psi\in C_c^\infty(\R^n)$ we can pass to the limit against $\nabla\psi$, since $|c_j|\leq C_*$. Hence
		\[
		\int_{\R^n}c(x,t)\cdot\nabla\psi(x)\,dx
		=\lim_j\int_{\R^n}c_j(x,t)\cdot\nabla\psi(x)\,dx=0.
		\]
		Thus $\div_xc(\cdot,t)=0$. To pass to the limit in the equation for $u_j$, note first that \eqref{eq:budget} and \eqref{eq:stage-conclusions} give
		\[
		\|\nabla(u_k-u_j)\|_2
		\leq\sum_{\ell=j+1}^k\ve_\ell\leq2^{-j-1}\quad(k>j).
		\]
		The gradient norm is an equivalent norm on $Y^{1,2}(\R^{n+1}_+)$, so $\{u_j\}_j$ converges in this space to a function $u$ with
		\begin{equation}\label{eq:global-gradient-tail}
			\|\nabla(u-u_j)\|_2\leq\sum_{\ell>j}\ve_\ell.
		\end{equation}
		Now fix $v\in Y^{1,2}(\R^{n+1}_+)$. We have
		\begin{align*}
			\left|\int_{\R^{n+1}_+}(A_j\nabla u_j-A\nabla u)\cdot\nabla v
			\,dx\,dt\right|
			&\leq(1+C_*)\|\nabla(u_j-u)\|_2\|\nabla v\|_2\\
			&\quad+\int_{\R^{n+1}_+}|A_j-A||\nabla u||\nabla v|
			\,dx\,dt\longrightarrow0.
		\end{align*}
		The first term tends to zero by \eqref{eq:global-gradient-tail}. In the second term the integrand tends to zero and is bounded by $2(1+C_*)|\nabla u||\nabla v|\in L^1$, so the integral tends to zero. Thus $u$ satisfies
		\[
		\int_{\R^{n+1}_+}A\nabla u\cdot\nabla v\,dx\,dt
		=\langle g,\Tr v\rangle.
		\]
		If $w$ is the difference of two such solutions, testing with $w$ gives
		$\|\nabla w\|_2^2=\int_{\R^{n+1}_+} A\nabla w\cdot\nabla w=0$. The Sobolev inequality gives $w=0$.
		
		To estimate the two integrals for $u$, we fix $i\geq1$. Passing to the limit in \eqref{eq:fixed-neighborhood} gives $A(x,t)=A_i(x,t)=I_{n+1}$ for $\dist(x,S_i)<4r_i$ and $t>0$, with $g(x)=0$ there. Thus the estimates defining $C_i$ apply to $u-u_i$ as well. By \eqref{eq:global-gradient-tail}, each of the last two integrals in \eqref{eq:stage-conclusions} changes by at most
		\[
		C_i|S_i|\sum_{j>i}\ve_j
		\leq\sum_{j>i}\frac{C_i|S_i|2^{-j-1}}{1+\max_{k<j}C_k|S_k|}
		\leq\sum_{j>i}2^{-j-1}=2^{-i-1}.
		\]
		Consequently,
		\begin{align}
			\int_{S_i}|\nabla_x\Tr u|\,dx&\geq4-2^{-i-1}\geq1,\label{eq:limit-retained}\\
			\int_{S_i}|(A\nabla u)_{n+1}(x,t_i)|\,dx=\int_{S_i}|\partial_tu(x,t_i)|\,dx
			&\geq4-2^{-i-1}\geq1.\label{eq:conormal-retained}
		\end{align}
		The sets $S_i\subset B_i$ are disjoint and $|\bigcup_i S_i|\leq\sum_i|B_i|<\infty$.
		
		When $\zeta=1$, every field in the sum defining $c$ is independent of $t$. Since the family $\{B_j^\infty\}_j$ is locally finite in $\R^n$, the sum is locally finite up to the boundary. It follows that
		\[
		\begin{gathered}
			c\in C^\infty(\R^n;\R^n)\cap L^\infty,\qquad
			\div_xc=0,\qquad\partial_tA=0,\\
			\supp c\subset\bigcup_j\overline B_j^\infty,\qquad
			|\supp c|\leq\sum_j|B_j^\infty|<\infty.
		\end{gathered}
		\]
		We used here that $\bigcup_j\overline B_j^\infty$ is closed by local finiteness. In this case the Carleson measure involving only $\partial_t A$ is identically zero.
		
		Suppose now that $(\zeta,B_j)=(\eta,B_j^0)$. By
		\eqref{eq:global-height-choice},
		\[
		\supp c_{j,\nu}\subset B_j^0\times[0,s_j],\qquad
		s_j=2^{-3j-10}\to0.
		\]
		If $K\Subset\R^{n+1}_+$ and $d_K=\min_{(x,t)\in K}t$, then every $c_{j,\nu}$ with $s_j<d_K/2$ vanishes within distance $d_K/2$ of $K$. Only finitely many indices $j$ satisfy $s_j\geq d_K/2$, and each has finitely many fields. Thus
		\[
		c\in C^\infty(\R^{n+1}_+;\R^n)\cap L^\infty,\qquad
		A=I_{n+1}\text{ on }\R^{n+1}_+\setminus(B(0,1)\times(0,1)).
		\]
		Since the initial balls are disjoint and the period cubes satisfy
		\eqref{eq:cell-restrictions} by \eqref{eq:global-height-choice},
		Lemma~\ref{lem:packing} gives \eqref{eq:packing}. The fields have disjoint
		horizontal supports, so we can apply Proposition~\ref{prop:dkp} to obtain
		\[
		\|\mu_A\|_{\mC}+\|\nu_A\|_{\mC}<\infty.
		\]
		For $x_0\neq0$, let $B_{x_0}\coloneqq B(x_0,|x_0|/4)$. If $\overline B_j^0$ meets $B_{x_0}$, then $3|x_0|/4\leq2^{-j}+2^{-3j-10}$. Hence, only finitely many initial balls meet $B_{x_0}$. Each contains the horizontal supports of finitely many fields, so
		\[
		t_{x_0}\coloneqq\min\Bigl(\{1\}\cup
		\{\rho_{j,\nu}:\supp_xc_{j,\nu}\cap B_{x_0}\ne\varnothing\}
		\Bigr)>0.
		\]
		Since $c_{j,\nu}(\cdot,t)=0$ for $t<\rho_{j,\nu}$, we have $A=I_{n+1}$ on $B_{x_0}\times(0,t_{x_0})$. At $x=0$ we also have $c(0,t)=0$, since none of the initial balls contains the origin. Finally $E_2\subset\bigcup_j B_j^0\subset B(0,1)$. This proves~\textup{(ii)}. For both coefficients the estimates \eqref{eq:limit-retained} and \eqref{eq:conormal-retained} hold for the datum \eqref{eq:fixed-datum}.
	\end{proof}

	\begin{lemma}\label{lem:normal-components}
		Let $A=A^{(i)}$ and $u=u^{(i)}$, with $i\in\{1,2\}$, be as in Proposition~\ref{prop:construction}. Then, for every $a>1$ and $1\leq p\leq\infty$,
		\[
		\|N_a(\partial_tu)\|_{L^p(\R^n)}
		=\|\widetilde N_a(\partial_tu)\|_{L^p(\R^n)}
		=\|\widetilde N_a((A\nabla u)_{n+1})\|_{L^p(\R^n)}=\infty.
		\]
		Moreover,
		\[
		\sup_{t>0}\|\partial_tu(\cdot,t)\|_2=\infty.
		\]
	\end{lemma}
	\begin{proof}
		Fix $i$ and suppress it from the notation. Let $x\in S_j$. Since $A(y,t)=I_{n+1}$ whenever $\dist(y,S_j)<b\rho_j/2$ and $t>0$, the function $\partial_tu$ is
		harmonic in $B((x,t_j),b\rho_j/2)$. The point $(x,t_j)$ belongs to the interior
		of $\gamma_a(x)$, so continuity gives
		$N_a(\partial_tu)(x)\geq|\partial_tu(x,t_j)|$. Let $B'=B((x,t_j),b\rho_j/2)$ and
		$B=B((x,t_j),t_j/2)$. Since $t_j=3\rho_j$, we have $B'\subset B$ and
		$|B'|/|B|=(b/3)^{n+1}$. For $h=\partial_tu$ or $h=(A\nabla u)_{n+1}$
		we have $h=\partial_tu$ on $B'$. The mean value inequality gives
		\[
		\widetilde N_a(h)(x)\geq\left(\avint_B|h|^2\right)^{1/2}
		\geq(b/3)^{(n+1)/2}\left(\avint_{B'}|\partial_tu|^2\right)^{1/2}
		\geq(b/3)^{(n+1)/2}|\partial_tu(x,t_j)|.
		\]
		By \eqref{eq:conormal-selected} and disjointness, each of the three maximal functions
		has infinite integral over $E$. Since $|E|<\infty$, each has infinite $L^p$ norm
		for every $1\leq p\leq\infty$.
		On each slice, Cauchy--Schwarz gives
		$\|\partial_tu(\cdot,t_j)\|_2\geq |S_j|^{-1/2}$, which tends to infinity
		because $\sum_j|S_j|<\infty$.
	\end{proof}
	
	\subsection{Completion of the proof}
	\begin{proof}[{\bf Proof of Theorem~\ref{thm:main}}]
		Take $A$ from Proposition~\ref{prop:construction}, $g$ from \eqref{eq:fixed-datum}, and let $u$ be the corresponding variational solution. Put $E=\bigcup_jS_j$. Since $(A\nabla u)_{n+1}$ is continuous near $S_j\times\{t_j\}$ and $(x,t_j)$ is interior to $\gamma_a(x)$, we have
		\[
		N_a((A\nabla u)_{n+1})(x)\geq|(A\nabla u)_{n+1}(x,t_j)|\qquad(x\in S_j,\ a>1).
		\]
		Integration over the disjoint sets $S_j$ and \eqref{eq:conormal-selected} yield
		\begin{equation}\label{eq:main-failure}
			\int_{E}N_a((A\nabla u)_{n+1})\,dx
			\geq\sum_j\int_{S_j}|(A\nabla u)_{n+1}(x,t_j)|\,dx=\infty.
		\end{equation}
		Thus $N_a((A\nabla u)_{n+1})$ cannot lie in any $L^p$, $1\leq p\leq\infty$, since $E$ has finite measure. As
		\[
		N_a((A\nabla u)_{n+1})\leq\|A\|_\infty N_a(\nabla u),
		\]
		we obtain the same conclusion for the full gradient.
		
		In a neighborhood of $S_j$ we have $A=I_{n+1}$ and $g=0$. By Lemma~\ref{lem:reflection} the tangential derivatives extend continuously across $S_j$. Taking their limit in the cone we obtain
		\[
		N_a(\nabla_xu)(x)\geq|\nabla_x\Tr u(x)|,
		\qquad x\in E.
		\]
		The lower bound \eqref{eq:main-selected} on the disjoint $S_j$ gives
		\[
		\int_E N_a(\nabla_xu)\,dx
		\geq\sum_j\int_{S_j}|\nabla_x\Tr u|\,dx=\infty.
		\]
		Since $|E|<\infty$, $N_a(\nabla_xu)\notin L^p$ for every $1\leq p\leq\infty$. Lemma~\ref{lem:normal-components} gives the same conclusion for $N_a(\partial_tu)$. This proves \eqref{eq:main-components}.
		
		Let $v\in Y^{1,2}_{\mathrm{loc}}(\R^{n+1}_+)$ be any solution of \eqref{eq:weak-neumann} with datum $g$. Suppose that $N_a(\nabla v)\in L^p(\R^n)$ for some $a>1$ and $1\leq p<\infty$. We claim $v-u$ is constant. Since $p\geq1$, $N_a(\nabla v)\in L^1_{\mathrm{loc}}$, and Lemma~\ref{lem:boundary-local-energy} gives $v\in Y^{1,2}_{\mathrm{loc}}(\overline{\R^{n+1}_+})$. Indeed, $c^{(1)}$ is smooth up to the boundary and independent of $t$, and $c^{(2)}$ vanishes near every boundary point other than the origin. The datum vanishes near the origin.
		
		Set $w=v-u$. The two weak identities give zero Neumann datum for $w$. For $\psi\in C_c^\infty(\R^{n+1})$, test with $\psi(x,t)+\psi(x,-t)$ on the upper half-space. Changing $t$ to $-t$ in the reflected term gives $\div(\mA\nabla w^\sharp)=0$ in $\R^{n+1}$, with $w^\sharp$ the even reflection and $\mA$ as in \eqref{sl:reflection-matrix}.
		
		The reflected matrix $\mA$ differs from $I_{n+1}$ only over $\bigcup_jB_j$.
		Each vertical line meets $B(0,R)$ in an interval of length at most $2R$.
		Since $\sum_j|B_j|<\infty$, we obtain
		\[
		\avint_{B(0,R)}|\mA-I_{n+1}|
		\leq \frac{2R\|A-I_{n+1}\|_\infty}{|B(0,R)|}\sum_j|B_j|
		\longrightarrow0\qquad(R\to\infty).
		\]
		Let $m_R$ be the average of $w^\sharp$ on $B(0,2R)$. The function
		$w^\sharp-m_R$ solves the same homogeneous equation. Applying the local
		boundedness estimate \cite[Lemma~2.4, \textup{(2.10)}]{HK07} with exponent
		one, and then the $L^1$ Poincar\'e inequality, gives
		\[
		\Phi_{w^\sharp}(R)\leq\|w^\sharp-m_R\|_{L^\infty(B(0,R))}
		\leq C\avint_{B(0,2R)}|w^\sharp-m_R|
		\leq CR\avint_{B(0,2R)}|\nabla w^\sharp|.
		\]
		Almost every $(x,t)$ is a Lebesgue point of $\nabla v$ in the interior of $\gamma_a(x)$. Thus
		$|\nabla v(x,t)|\leq N_a(\nabla v)(x)$. Reflection and H\"older's inequality give
		\[
		\int_{B(0,2R)}|\nabla v^\sharp|\leq 4R\int_{|x|<2R}N_a(\nabla v)(x)\,dx
		\leq C R^{n+1-n/p}\|N_a(\nabla v)\|_p.
		\]
		For $u^\sharp$, Cauchy--Schwarz gives the bound
		$\int_{B(0,2R)}|\nabla u^\sharp|\leq C R^{(n+1)/2}\|\nabla u\|_2$.
		Since $w^\sharp=v^\sharp-u^\sharp$, we conclude that
		\[
		\Phi_{w^\sharp}(R)\leq C\bigl(R^{1-n/p}\|N_a(\nabla v)\|_p
		+R^{(1-n)/2}\|\nabla u\|_2\bigr).
		\]
		For $R\geq1$ this is bounded by $CR^\beta$ with
		$\beta=\max\{0,1-n/p\}<1$, since $n\geq2$ and $p<\infty$.
		Lemma~\ref{lem:liouville-at-infinity} now gives that $w^\sharp$ is constant. Hence $\nabla v=\nabla u$. This contradicts \eqref{eq:main-failure}, since $N_a((A\nabla u)_{n+1})\leq\|A\|_\infty N_a(\nabla v)$ and $N_a(\nabla v)\in L^p$.
		
		For $p=\infty$ the same estimate gives only $\Phi_{w^\sharp}(R)\leq CR$ for $R\geq1$, whereas Lemma~\ref{lem:liouville-at-infinity} requires a bound by $CR^\beta$ with $\beta<1$.
	\end{proof}
	
	\begin{remark}
		The same argument applies to $\widetilde N_a(\nabla_xu)$. For every $i$ and $x\in S_i$, we have
		\[
		\lim_{t \to 0^+}
		\left(\avint_{B((x,t),t/2)}|\nabla_xu|^2\right)^{1/2}
		=|\nabla_x\Tr u(x)|.
		\]
		Indeed, continuity of the reflected tangential gradient gives the limit. Since $(x,t)\in\gamma_a(x)$ for each $t>0$, \eqref{eq:main-selected} shows that $\int_{\bigcup_i S_i}\widetilde N_a(\nabla_xu)\,dx=\infty$.
	\end{remark}
	
	\begin{remark}\label{rem:laplace-neumann}
		For $A=I_{n+1}$ one can read the Neumann estimate directly from the Fourier transform. Let $1<p<\infty$ and $g\in L^p(\R^n)\cap\dot H^{-1/2}(\R^n)$. The variational solution is given by
		\[
		\widehat u(\xi,t)=|\xi|^{-1}e^{-t|\xi|}\widehat g(\xi),\qquad
		\|\nabla u\|_2^2=\int_{\R^n}|\xi|^{-1}|\widehat g(\xi)|^2\,d\xi,
		\]
		and hence
		\[
		-u_t=P_tg,\qquad
		\partial_ku=P_tR_kg,\qquad
		1\leq k\leq n.
		\]
		If $\supp\widehat g\subset\{\xi:r\leq|\xi|\leq R\}$ for some $0<r<R<\infty$, these formulas follow
		from Plancherel's theorem and integration in $t$. To pass to general
		$g\in L^p\cap\dot H^{-1/2}$, take a real even cutoff
		$\phi\in C_c^\infty(\R^n)$ equal to one near zero and put
		$\widehat g_j(\xi)=\phi(\xi/j)(1-\phi(j\xi))\widehat g(\xi)$.
		Dominated convergence gives $g_j\to g$ in $\dot H^{-1/2}$. We also have
		$g_j\to g$ in $L^p$. Indeed, the multipliers $\phi(\cdot/j)$ form an
		approximation of the identity. The multipliers $\phi(j\,\cdot)$ are
		uniformly bounded on $L^p$ and tend strongly to zero there. Indeed, on $L^1\cap L^p$
		their convolution kernels have $L^p$ norm $C j^{-n(1-1/p)}$, and
		$L^1\cap L^p$ is dense in $L^p$. The variational bound gives convergence
		of the solution gradients in $L^2(\R^{n+1}_+)$. For each $t>0$, $P_tg_j\to P_tg$ and $P_tR_kg_j\to P_tR_kg$ in $L^p(\R^n)$, which gives the formulas for $\partial_tu$ and $\partial_ku$. The Poisson maximal inequality, followed by the $L^p$ estimates for $\M$ and $R_k$ \cite[Chapter~I, Section~1, and Chapter~III, Sections~1--2]{Stein}, gives
		\[
		N_a(\nabla u)\leq C_{n,a}\left(\M g+\sum_{k=1}^n\M(R_kg)\right),
		\qquad
		\|N_a(\nabla u)\|_p\leq C_{n,p,a}\|g\|_p.
		\]
		
		For arbitrary $g\in L^p(\R^n)$ we take $g_j\in C_c^\infty(\R^n)$ converging to $g$ in $L^p$ and let $u_j$ be the corresponding variational solutions for the identity matrix. Applying the preceding estimate to the differences, we have
		\[
		\|N(\nabla(u_j-u_k))\|_p\leq C_{n,p}\|g_j-g_k\|_p.
		\]
		If $|x-y|<t$, then $(y,t)\in\gamma_2(x)$. Thus
		$|\nabla(u_j-u_k)(y,t)|\leq |B(y,t)|^{-1/p}\|N(\nabla(u_j-u_k))\|_p$.
		Since $t$ is bounded below on each compact subset of the open half-space, the gradients are uniformly Cauchy there. After normalizing the solutions at one fixed interior point, we obtain a harmonic function $u$ with $\nabla u_j\to\nabla u$ locally uniformly. To pass to the boundary identity, we also need convergence of the gradients in $L^1$ on cylinders that meet the boundary. For every bounded cube $Q$ and $T>0$ we have
		\[
		\int_{Q\times(0,T)}|\nabla(u_j-u_k)|
		\leq T|Q|^{1-1/p}\|N(\nabla(u_j-u_k))\|_p.
		\]
		This gives local $L^1$ convergence of the gradients up to the boundary. Passing to the limit against each test function gives the weak Neumann identity with datum $g$, and Fatou's lemma yields \eqref{eq:Np}. This proves solvability in the sense of Definition~\ref{def:neumann-solvability}.
	\end{remark}

	\section{Single layer potentials}\label{sec:layers}
	We first show that for smooth densities the integral in \eqref{sl:kernel} gives the single layer
	defined in \eqref{sl:variational} (see also \cite[Section~6]{HMM15}). We will use both expressions
	when we estimate the truncations of the boundary kernel.
	
	\subsection{Kernel representation and slice estimates}\label{sl:section}
	Let $\cA$ be a real bounded uniformly elliptic matrix on $\R^{n+1}$. We use $\lambda$ for its ellipticity
	constant and $\Lambda$ for the bound on its operator norm. Unless otherwise indicated, constants in this
	subsection depend only on $n$, $\lambda$, and $\Lambda$. By Hofmann and Kim
	\cite[Theorem~3.1, \textup{(3.47)--(3.48)}]{HK07}, if $F\in C_c^\infty(\R^{n+1})$, the solution
	$v\in Y^{1,2}(\R^{n+1})$ of $L_{\cA}v=F$ is given by
	$v(z)=\int\Gamma_{\cA}(z,w)F(w)\,dw$.
	
	\begin{lemma}\label{sl:kernel-variational}
		For $f\in C_c^\infty(\R^n)$, the continuous representative of $\mS_{\cA}f$ in each open half-space is
		given by \eqref{sl:kernel}. If $\cA$ is smooth in one of the open half-spaces, the derivatives of
		$\mS_{\cA}f$ there are obtained by differentiating under the integral.
	\end{lemma}
	\begin{proof}
		We choose $0\leq\eta\in C_c^\infty((-1,1))$ with $\int\eta=1$ and set $\eta_\ve(t)=\ve^{-1}\eta(t/\ve)$ for $\ve>0$.
		Let $v_\ve\in Y^{1,2}(\R^{n+1})$ solve
		\[
		L_{\cA}v_\ve(x,t)=f(x)\eta_\ve(t).
		\]
		For $\varphi\in C_c^\infty(\R^{n+1})$, we integrate $\partial_t\varphi$ from $0$ to $s$ and apply
		Cauchy--Schwarz to get
		\[
		\|\varphi(\cdot,s)-\Tr\varphi\|_2
		\leq |s|^{1/2}\|\partial_t\varphi\|_{L^2(\R^{n+1})},
		\]
		and since $\int\eta_\ve=1$, it follows that
		\begin{align*}
			\left|\int_{\R^{n+1}}f(x)\eta_\ve(t)\varphi(x,t)\,dx\,dt
			-\langle f,\Tr\varphi\rangle\right|
			&\leq\|f\|_2\int\eta_\ve(s)|s|^{1/2}\,ds
			\|\nabla\varphi\|_2\\
			&\lesssim\ve^{1/2}\|f\|_2\|\nabla\varphi\|_2.
		\end{align*}
		By density we may take $\varphi=v_\ve-\mS_{\cA}f$ in this estimate. Subtracting the two
		equations and using ellipticity, we obtain
		\[
		\|\nabla(v_\ve-\mS_{\cA}f)\|_2
		\leq C_{\lambda,\eta}\ve^{1/2}\|f\|_2.
		\]
		By the Sobolev inequality $v_\ve\to\mS_{\cA}f$ in $L^{2(n+1)/(n-1)}(\R^{n+1})$. In particular,
		$v_{\ve_k}\to\mS_{\cA}f$ almost everywhere for some decreasing sequence $\ve_k \to 0^+$.
		
		On the other hand, fix $(x,t)$ with $t\ne0$ and take $\ve<|t|/2$. The function
		$\Gamma_{\cA}((x,t),\cdot)$ is uniformly continuous on $\supp f\times[-|t|/2,|t|/2]$.
		The representation of $v_\ve$ therefore gives
		\begin{align*}
			&\left|v_\ve(x,t)
			-\int_{\R^n}\Gamma_{\cA}((x,t),(y,0))f(y)\,dy\right|\\
			&\qquad\leq\|f\|_1
			\sup_{\substack{y\in\supp f\\|s|\leq\ve}}
			|\Gamma_{\cA}((x,t),(y,s))-\Gamma_{\cA}((x,t),(y,0))|
			\longrightarrow0.
		\end{align*}
		It follows that the integral in \eqref{sl:kernel} agrees with $\mS_{\cA}f$ almost everywhere. Both
		are continuous in each open half-space, the first by continuity of the kernel off the diagonal and
		compactness of $\supp f$, the second by interior H\"older regularity. This proves the identity everywhere
		off the plane. If $\cA$ is smooth in one half-space, interior estimates for
		$\Gamma_{\cA}(\cdot,(y,0))$ on compact subsets of that half-space are uniform in $y\in\supp f$.
		We may thus differentiate under the integral there.
	\end{proof}
	
	For general $f\in\dot H^{-1/2}(\R^n)$, the outward conormal derivatives from the two half-spaces satisfy
	\begin{equation}\label{sl:jump}
		\partial_{\nu_{\cA}}^{+}\mS_{\cA}f
		+\partial_{\nu_{\cA}}^{-}\mS_{\cA}f=f.
	\end{equation}
	Indeed, given $h\in\dot H^{1/2}(\R^n)$, the trace theorem gives a function in
	$Y^{1,2}(\R^{n+1}_+)$ with trace $h$. Its even extension belongs to $Y^{1,2}(\R^{n+1})$,
	as shown in the proof of Lemma~\ref{sl:reflection}.
	Testing \eqref{sl:variational} with this function and splitting the integral at $t=0$
	gives \eqref{sl:jump} by \eqref{sl:weak-conormals}. The outward normals are $-e_{n+1}$ above and
	$e_{n+1}$ below, which accounts for the plus sign.
	
	We next specify the representative used for the slices of $\nabla u$. Suppose $\cA$ is locally
	Lipschitz in both open half-spaces and $u\in Y^{1,2}_{\mathrm{loc}}(\R^{n+1})$ satisfies
	$L_{\cA}u=0$ in each of them. Then $\nabla^2 u$ is locally square integrable off the plane by the
	interior difference-quotient estimate \cite[Theorem~8.8]{GT}. In particular, for bounded $K\subset\R^n$
	and $0<s_0<s_1$, integration in $t$ gives an absolutely continuous representative of $\nabla u(\cdot,t)$
	in $L^2(K)$, with
	\[
	\|\nabla u(\cdot,t)-\nabla u(\cdot,s)\|_{L^2(K)}
	\leq |t-s|^{1/2}
	\|\partial_t\nabla u\|_{L^2(K\times(s_0,s_1))}
	\quad(s,t\in(s_0,s_1)).
	\]
	The representatives agree on overlapping cylinders by continuity. We use the same convention below
	the plane.
	
	\subsection{Reflection and the Neumann solution}\label{sl:reflected-section}
	
	Since $\mathcal J$ is orthogonal, conjugation by $\mathcal J$ in \eqref{sl:reflection-matrix} preserves ellipticity and the operator
	norm. The reflection $(x,t)\mapsto(x,-t)$ also preserves distances and sends Whitney balls to
	Whitney balls. Consequently the Carleson norms defined using $\nabla A$, $\partial_t A$ and
	the discrepancy from the identity are the same on the two sides.
	
	For the even reflection, the right-hand side of the weak identity in
	\cite[proof of Lemma~6.4, p.~463]{AM14} is $-2\langle\alpha,\Tr\varphi\rangle$,
	where $\alpha$ is the inward conormal derivative in \cite[(11), p.~451]{AM14}.
	In the following lemma the outward conormal derivative is $f/2$, so $\alpha=-f/2$. The right-hand side is $\langle f,\Tr\varphi\rangle$.
	
	\begin{lemma}\label{sl:reflection}
		Let $n\geq2$, let $A$ be a real bounded uniformly elliptic matrix in $\R^{n+1}_+$, and let $\mA$ be given
		by \eqref{sl:reflection-matrix}. Given $f\in\dot H^{-1/2}(\R^n)$, let $u\in Y^{1,2}(\R^{n+1}_+)$ solve
		\eqref{eq:variational} with datum $f/2$. Then
		\begin{equation}\label{sl:even-identity}
			\mS_{\mA}f(x,t)=u(x,|t|),\qquad
			\partial_{\nu_{\mA}}^{\pm}\mS_{\mA}f=f/2.
		\end{equation}
		The first identity holds for the continuous representatives off the plane, and the conormal identities
		hold in $\dot H^{-1/2}(\R^n)$.
	\end{lemma}
	\begin{proof}
		We set $u^{\sharp}(x,t)\coloneqq u(x,|t|)$. We first compute its weak derivatives on
		$\R^{n+1}$. Let $\varphi\in C_c^\infty(\R^{n+1})$ and choose a bounded cube $Q\subset\R^n$
		and $T>0$ so that $\supp\varphi\subset Q\times(-T,T)$. Since $u\in L^{2^*}$, $2^*>2$ and
		$\nabla u\in L^2$, we have $u\in W^{1,2}(Q\times(0,T))$. Thus, for almost every $x\in Q$, the
		function $t\mapsto u(x,t)$ belongs to $W^{1,2}(0,T)$ with weak derivative
		$\partial_tu(x,\cdot)$. Changing variables in the lower half-space and integrating
		by parts in $t$, we obtain
		\begin{align*}
			\int_{\R^{n+1}}u^{\sharp}\,\partial_t\varphi\,dx\,dt
			&=\int_{\R^{n+1}_+}u(x,t)\,\partial_t\bigl(\varphi(x,t)-\varphi(x,-t)\bigr)\,dx\,dt\\
			&=-\int_{\R^{n+1}_+}\partial_tu(x,t)\bigl(\varphi(x,t)-\varphi(x,-t)\bigr)\,dx\,dt.
		\end{align*}
		There is no boundary term at $t=0$, since the test function $\varphi(x,t)-\varphi(x,-t)$ vanishes there.
		Changing variables back in the term containing $\varphi(x,-t)$ gives
		$\partial_tu^{\sharp}(x,t)=\operatorname{sgn}(t)\,\partial_tu(x,|t|)$ as a weak derivative
		on $\R^{n+1}$. For $1\leq j\leq n$, we integrate by parts in $x_j$ against
		$\varphi(x,t)+\varphi(x,-t)$ for almost every $t>0$ and obtain $\partial_j u^{\sharp}(x,t)=\partial_j u(x,|t|)$.
		
		It follows by a change of variables that
		\[
		\|u^{\sharp}\|_{L^{2^*}(\R^{n+1})}^{2^*}
		=2\|u\|_{L^{2^*}(\R^{n+1}_+)}^{2^*}
		\qquad\text{and}\qquad
		\|\nabla u^{\sharp}\|_{L^2(\R^{n+1})}^2
		=2\|\nabla u\|_{L^2(\R^{n+1}_+)}^2.
		\]
		Both quantities are finite, so $u^{\sharp}\in Y^{1,2}(\R^{n+1})$. By the derivative
		formulas above and \eqref{sl:reflection-matrix}, for almost every $(x,t)\in\R^{n+1}_+$ we have
		\[
		\nabla u^{\sharp}(x,-t)=\mathcal{J}\nabla u(x,t),\qquad
		\mA(x,-t)\nabla u^{\sharp}(x,-t)
		=\mathcal{J} A(x,t)\nabla u(x,t).
		\]
		For the test function $\varphi$ above, we set $\psi\coloneqq\varphi+\varphi\circ\mathcal J$
		on $\R^{n+1}_+$. Since $\mathcal J$ is symmetric,
		$\nabla(\varphi\circ\mathcal J)(x,t)=\mathcal J\nabla\varphi(x,-t)$.
		Changing variables in the lower half-space and using \eqref{eq:variational}, we obtain
		\begin{align*}
			\int_{\R^{n+1}}\mA\nabla u^{\sharp}\cdot\nabla\varphi\,dx\,dt
			&=\int_{\R^{n+1}_+}\bigl(A\nabla u\cdot\nabla\varphi
			+\mathcal J A\nabla u\cdot\nabla\varphi(x,-t)\bigr)\,dx\,dt\\
			&=\int_{\R^{n+1}_+}A\nabla u\cdot\nabla\psi\,dx\,dt
			=\tfrac12\langle f,\Tr\psi\rangle
			=\langle f,\Tr\varphi\rangle,
		\end{align*}
		where the last equality follows from $\psi(x,0)=2\varphi(x,0)$. By density, this gives
		\eqref{sl:variational}. Uniqueness and interior continuity give
		$\mS_{\mA}f=u^{\sharp}$ at every point off the plane. In the upper half-space its outward conormal derivative is $f/2$ by
		the choice of $u$. By \eqref{sl:jump}, its outward conormal derivative in the lower half-space is also $f/2$.
	\end{proof}
	
	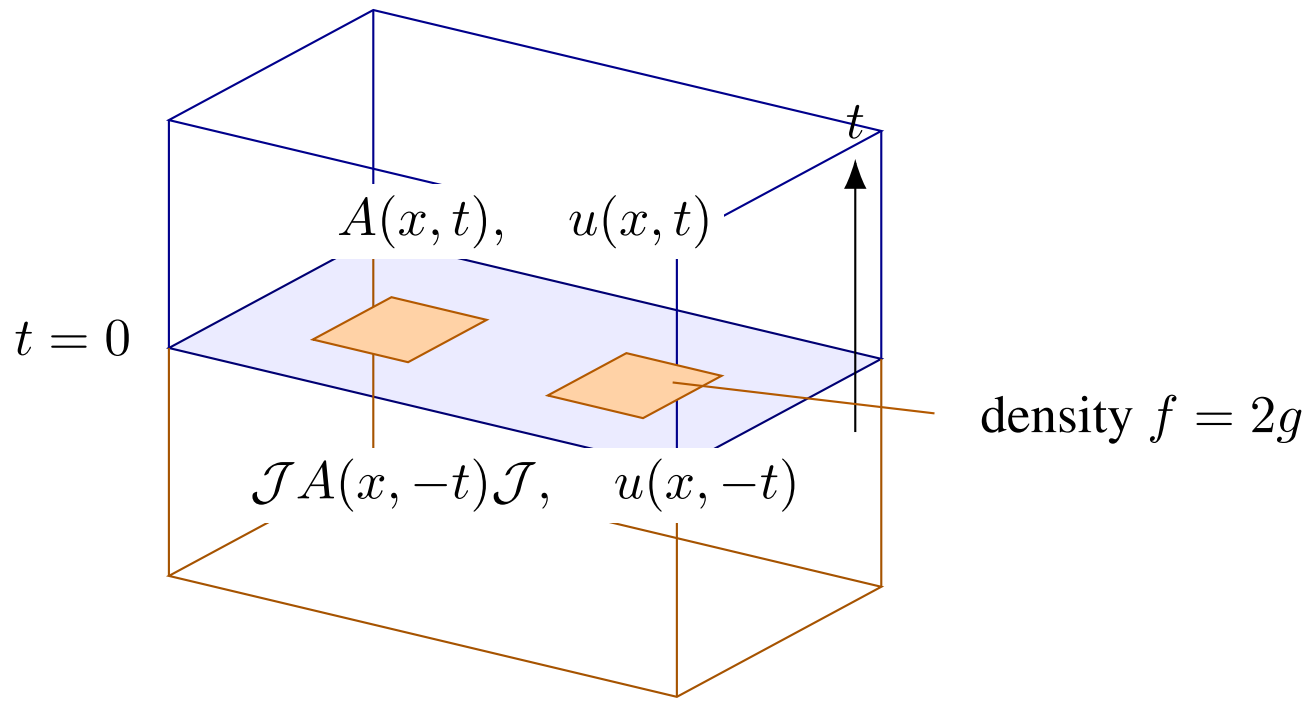
\begin{figure}[htbp]
		\centering
		{ \begin{tikzpicture}[font=\fontsize{10}{12}\selectfont,
				x={(.84cm,-.20cm)},y={(.52cm,.28cm)},z={(0cm,.67cm)}]
				\path[fill=blue!8,draw=blue!45!black]
				(-2,-1.3,0)--(2,-1.3,0)--(2,1.3,0)--(-2,1.3,0)--cycle;
				\draw[blue!55!black] (-2,-1.3,2.25)--(2,-1.3,2.25)--(2,1.3,2.25)--(-2,1.3,2.25)--cycle;
				\draw[blue!55!black] (-2,-1.3,0)--(-2,-1.3,2.25)
				(2,-1.3,0)--(2,-1.3,2.25)
				(2,1.3,0)--(2,1.3,2.25)
				(-2,1.3,0)--(-2,1.3,2.25);
				\draw[orange!65!black] (-2,-1.3,-2.25)--(2,-1.3,-2.25)--(2,1.3,-2.25)--(-2,1.3,-2.25)--cycle;
				\draw[orange!65!black] (-2,-1.3,0)--(-2,-1.3,-2.25)
				(2,-1.3,0)--(2,-1.3,-2.25)
				(2,1.3,0)--(2,1.3,-2.25)
				(-2,1.3,0)--(-2,1.3,-2.25);
				\path[draw=orange!70!black,fill=orange!35]
				(-1.3,-.6,0)--(-.55,-.6,0)--(-.55,.4,0)--(-1.3,.4,0)--cycle;
				\path[draw=orange!70!black,fill=orange!35]
				(.55,-.6,0)--(1.3,-.6,0)--(1.3,.4,0)--(.55,.4,0)--cycle;
				\node[fill=white,inner sep=2pt] at (0,0,1.3) {$A(x,t),\quad u(x,t)$};
				\node[fill=white,inner sep=2pt] at (0,0,-1.3)
				{$\mathcal{J} A(x,-t)\mathcal{J} ,\quad u(x,-t)$};
				\draw[-{Latex[length=2mm]}] (2.6,0,0)--(2.6,0,2.7)
				node[above] {$t$};
				\node[anchor=west] at (3.3,.2,.25) {density $f=2g$};
				\draw[orange!70!black] (3.1,.2,.25)--(1.1,.1,0);
				\node[anchor=east] at (-2.2,-1.2,0) {$t=0$};
		\end{tikzpicture}}
		\caption{The matrix $\mA$ and the even extension of the Neumann solution from the upper half-space in dimension three.
				The density on the plane is $2g$, and the outward conormal derivative is $g$ on each side.
				The two shaded regions indicate the supports of the two terms in \eqref{eq:fixed-datum}.}
		\label{sl:reflection-figure}
	\end{figure}
	
	\subsection{Proof of Theorem~\ref{sl:negative}}\label{lp:failure-section}
	
	Let $A=A^{(i)}$, $i\in\{1,2\}$, and $f=2g$ be as in Theorem~\ref{sl:negative}. In addition to
	\eqref{sl:negative-maximal}, we have, for every $a>1$ and $1\leq p\leq\infty$,
	\begin{equation}\label{sl:negative-normal-maximal}
		\widetilde N_a^\pm(\nabla_x\mS_{\mA}f),\quad
		\widetilde N_a^\pm(\partial_t\mS_{\mA}f),\quad
		\widetilde N_a^\pm\bigl((\mA\nabla\mS_{\mA}f)_{n+1}\bigr),\quad
		\widetilde N_a^\pm(\nabla\mS_{\mA}f)\notin L^p(\R^n).
	\end{equation}
	The distribution $\nabla_x\Tr\mS_{\mA}f$ has no representative in $L^p(\R^n)$, and
	\begin{equation}\label{sl:negative-slices}
		\liminf_{t \to 0^+}\|\nabla_x\mS_{\mA}f(\cdot,\pm t)\|_{L^p(\R^n)}
		=\liminf_{t \to 0^+}\|\nabla\mS_{\mA}f(\cdot,\pm t)\|_{L^p(\R^n)}=\infty,
		\qquad 1\leq p\leq\infty.
	\end{equation}
	The ordinary normal derivative also satisfies
	\begin{equation}\label{sl:negative-normal-slices}
		\sup_{t\ne0}\|\partial_t\mS_{\mA}f(\cdot,t)\|_2=\infty.
	\end{equation}
	
	\begin{proof}
		Fix $i\in\{1,2\}$ and let $u_g$ be the solution of \eqref{eq:variational} for $A=A^{(i)}$ and the
		datum $g$ in Proposition~\ref{prop:construction}. By Lemma~\ref{sl:reflection},
		$\mS_{\mA}(2g)=u_g^{\sharp}$. Since $\mathcal Je_{n+1}=-e_{n+1}$, for $t>0$ we have
		\[
		\bigl(\mA\nabla\mS_{\mA}(2g)\bigr)_{n+1}(x,\pm t)=\pm (A\nabla u_g)_{n+1}(x,t).
		\]
		The last equality in \eqref{sl:negative-maximal} now follows from Theorem~\ref{thm:main}. For $t>0$ we also have
		\[
		\partial_t\mS_{\mA}(2g)(x,-t)=-\partial_tu_g(x,t).
		\]
		Lemma~\ref{lem:normal-components} gives infinite $L^p$ norm for the ordinary and averaged maximal
		functions of $\partial_tu_g$ and for the averaged maximal function of $(A\nabla u_g)_{n+1}$.
		It also gives \eqref{sl:negative-normal-slices}. Reflection preserves the absolute values
		of these functions and the averages of their squares over corresponding balls.
		For the tangential estimates, we use the pairwise disjoint compact sets $S_j$ given by
		Proposition~\ref{prop:construction}. In both examples we have
		\[
		E\coloneqq \bigcup_jS_j,\qquad |E|<\infty,\qquad
		\int_{S_j}|\nabla_x\Tr u_g|\,dx\geq1.
		\]
		Moreover, $E\subset B(0,1)$ if $i=2$. In a neighborhood of $S_j\times\{0\}$ we have $A=I_{n+1}$
		and $g=0$. Hence $u_g^{\sharp}$ is harmonic there, and the distributional tangential derivative of
		its trace agrees with the classical derivative on $S_j$. It follows that
		\[
		\int_E|\nabla_x\Tr u_g|\,dx
		=\sum_j\int_{S_j}|\nabla_x\Tr u_g|\,dx=\infty.
		\]
		An $L^p(\R^n)$ representative of $\nabla_x\Tr\mS_{\mA}(2g)$ would therefore coincide with
		$\nabla_x\Tr u_g$ almost everywhere on $E$. This is impossible for $1\leq p\leq\infty$, since
		$L^p(E)\subset L^1(E)$.
		
		For $x\in S_j$, smoothness of $u_g^{\sharp}$ also gives
		\[
		\lim_{t \to 0^+}|\nabla_xu_g(x,t)|
		=\lim_{t \to 0^+}
		\left(\avint_{B((x,t),t/2)}|\nabla_xu_g|^2\right)^{1/2}
		=|\nabla_x\Tr u_g(x)|.
		\]
		For all sufficiently small $t>0$, $(x,t)$ is an interior point of $\gamma_a(x)$ at which the
		gradient is continuous. Taking essential suprema and letting $t \to 0^+$, we obtain on $E$
		\[
		N_a(\nabla_xu_g)\geq|\nabla_x\Tr u_g|,
		\qquad
		\widetilde N_a(\nabla_xu_g)\geq|\nabla_x\Tr u_g|.
		\]
		The right-hand side has infinite integral on $E$, so neither maximal function is in $L^p$.
		The same holds for the full gradient and, by reflection, for the corresponding maximal functions
		in the lower half-space.
		
		For $1\leq p<\infty$, we obtain by Fatou's lemma
		\[
		\liminf_{t \to 0^+}\int_{\R^n}|\nabla_xu_g(x,t)|^p\,dx
		\geq\int_E|\nabla_x\Tr u_g|^p\,dx=\infty.
		\]
		For $p=\infty$, let $K>0$. Since $\nabla_x\Tr u_g\notin L^\infty(E)$, we may choose $x_K\in S_j$ for some
		$j$ with $|\nabla_x\Tr u_g(x_K)|>2K$. By continuity of the gradient near $(x_K,0)$, there exist
		$r_K,t_K>0$ such that
		\[
		|\nabla_xu_g(x,t)|>K
		\qquad\bigl(x\in B(x_K,r_K),\ 0<t<t_K\bigr).
		\]
		Hence $\|\nabla_xu_g(\cdot,t)\|_\infty\geq K$ for every $0<t<t_K$. Since $|\nabla u_g|\geq|\nabla_xu_g|$,
		this proves \eqref{sl:negative-slices} in the upper half-space. The lower slice norms equal the upper ones by
		reflection. The Carleson bounds on the reflected coefficient follow from
		Subsection~\ref{sl:reflected-section}.
	\end{proof}
	
	For $i=2$ and $x\ne0$, $A$ is the identity in a neighborhood of $(x,0)$. There $u_g$ is harmonic
	with smooth Neumann datum $g$. Since the upper outward normal is $-e_{n+1}$, local boundary
	regularity and the identity $\mS_{\mA}(2g)=u_g^{\sharp}$ give
	\[
	\lim_{t \to 0^+}\partial_t\mS_{\mA}(2g)(x,t)=-g(x),\qquad
	\lim_{t \to 0^-}\partial_t\mS_{\mA}(2g)(x,t)=g(x)
	\quad\text{for a.e. }x\in\R^n.
	\]
	Thus the one-sided limits of $\partial_t\mS_{\mA}(2g)$ agree almost everywhere with the smooth
	functions $-g$ and $g$, respectively, although its non-tangential maximal function has infinite
	$L^p$ norm for every $1\leq p\leq\infty$ and its $L^2$ slice norms are unbounded.
	
	On $E$, the truncated integral below equals $\nabla_x\Tr\mS_{\mA}(2g)$ whenever $0<\ve<1/2$.
	
	\begin{corollary}\label{sl:boundary-truncations}
		Let $n\geq2$, $i\in\{1,2\}$, and let $A=A^{(i)}$, $g$, and $S_j=S_j^{(i)}$ be as in
		Proposition~\ref{prop:construction}. We set $f\coloneqq 2g$ and $E\coloneqq\bigcup_j S_j$, and take $\mA$
		as in \eqref{sl:reflection-matrix}. For $x\in E$ and $\ve>0$, let
		\[
		T_\ve f(x)
		\coloneqq \int_{|x-y|>\ve}
		\nabla_x\Gamma_{\mA}((x,0),(y,0))f(y)\,dy.
		\]
		Then
		\[
		\|T_\ve f\|_{L^2(E)}=\infty\qquad\text{for every }0<\ve<1/2.
		\]
		In particular, the truncations on smooth compactly supported densities admit no uniform $L^2(\R^n)\to L^2(E)$ bound.
	\end{corollary}
	\begin{proof}
		By \eqref{eq:fixed-datum}, $\supp f\subset B(3e_1,1/8)\cup B(-3e_1,1/8)$. For $i=1$,
		\eqref{eq:initial-balls} gives $S_j\subset B(4je_1,2^{-j-10})$ for every $j\geq1$, and for $i=2$ we have $E\subset B(0,1)$. Hence,
		\[
		\dist(E,\supp f)\geq
		\begin{cases}
		1-2^{-11}-1/8,&i=1,\\
		3-1-1/8,&i=2,
		\end{cases}
		\qquad\text{and in both cases }\dist(E,\supp f)>1/2.
		\]
		Fix $j$ and choose a bounded neighborhood of $S_j\times\{0\}$ where $\mA=I$ and whose closure
		is disjoint from $\supp f\times\{0\}$. For $y\in\supp f$, $\Gamma_{\mA}(\cdot,(y,0))$ is harmonic
		there and uniformly bounded in $y$ by the fundamental solution estimate of Hofmann and Kim
		\cite[Theorem~3.1, \textup{(3.55)}]{HK07}. The interior gradient estimate gives a uniform bound
		on its gradient on a smaller neighborhood of $S_j\times\{0\}$. Thus we may differentiate
		\eqref{sl:kernel} in $x$ and let $t\to0^+$. For every $x\in S_j$ and $0<\ve<1/2$, we have
		$|x-y|>\ve$ for every $y\in\supp f$. It follows that
		\[
		T_\ve f(x)=\nabla_x\Tr\mS_{\mA}f(x).
		\]
		This holds on all of $E$, and $\nabla_x\Tr\mS_{\mA}f\notin L^2(E)$.
	\end{proof}
	
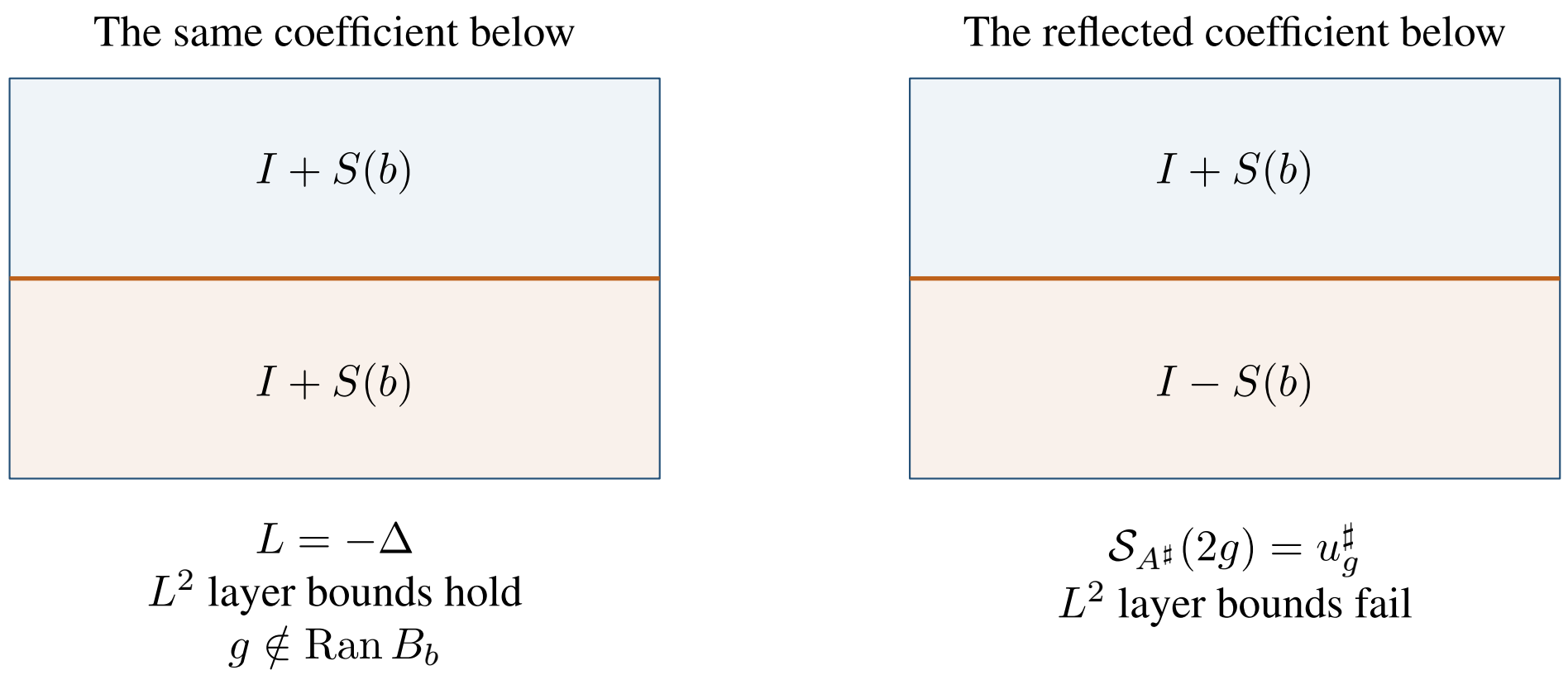
\begin{figure}[htbp]
		\centering
		\begin{tikzpicture}[>=Latex,font=\small,line join=round]
			\begin{scope}[shift={(-3.6,0)}]
				\fill[lightFill] (-2.6,0) rectangle (2.6,1.6);
				\fill[mathOrange!9] (-2.6,-1.6) rectangle (2.6,0);
				\draw[mathBlue] (-2.6,-1.6) rectangle (2.6,1.6);
				\draw[mathOrange,line width=1pt] (-2.6,0)--(2.6,0);
				\node at (0,.85) {$I+S(b)$};
				\node at (0,-.85) {$I+S(b)$};
				\node[above,align=center] at (0,1.72) {The same coefficient below};
				\node[below,align=center] at (0,-1.82)
				{$L=-\Delta$\\$L^2$ layer bounds hold\\$g\notin\operatorname{Ran}B_b$};
			\end{scope}
			\begin{scope}[shift={(3.6,0)}]
				\fill[lightFill] (-2.6,0) rectangle (2.6,1.6);
				\fill[mathOrange!9] (-2.6,-1.6) rectangle (2.6,0);
				\draw[mathBlue] (-2.6,-1.6) rectangle (2.6,1.6);
				\draw[mathOrange,line width=1pt] (-2.6,0)--(2.6,0);
				\node at (0,.85) {$I+S(b)$};
				\node at (0,-.85) {$I-S(b)$};
				\node[above,align=center] at (0,1.72) {The reflected coefficient below};
				\node[below,align=center] at (0,-1.82)
				{$\mS_{\mA}(2g)=u_g^{\sharp}$\\$L^2$ layer bounds fail};
			\end{scope}
		\end{tikzpicture}
		\caption{Two extensions of $A^{(1)}$ to the lower half-space. On the left, $g$ is outside the range of
			$B_b$. On the right, the single layer with density $2g$ is the even extension of the Neumann solution with datum $g$.}
		\label{lp:two-extensions}
	\end{figure}

	\begin{remark}[The choice of coefficients in the lower half-space]\label{rem:two-extensions}
		Let $A^{(1)}=I+S(b)$ and $g$ be as in Theorem~\ref{thm:main}. Both $A_b$ in \eqref{eq:unchanged_extension}
		and $\mA$ in \eqref{sl:reflection-matrix} agree with $A^{(1)}$ in the upper half-space.
		
		Since $b$ is independent of $t$ and $\div_xb=0$, the contribution of $S(b)$ to the bilinear form
		on $\R^{n+1}$ is zero. Hence, $L_{A_b}=-\Delta$ in $\R^{n+1}$ and
		the Fourier representation in Remark~\ref{rem:laplace-neumann} gives
		\begin{equation}\label{sl:unchanged-formula}
			\mS_{A_b}f=\tfrac12P_{|t|}(-\Delta_x)^{-1/2}f,
			\qquad f\in L^2(\R^n)\cap\dot H^{-1/2}(\R^n).
		\end{equation}
		For $t>0$, differentiation gives
		\[
		\partial_t\mS_{A_b}f(\cdot,t)=-\tfrac12P_tf,\qquad
		\partial_j\mS_{A_b}f(\cdot,t)=\tfrac12P_tR_jf,\quad 1\leq j\leq n,
		\]
		where $R_j=\partial_j(-\Delta_x)^{-1/2}$. The Poisson maximal estimate and the $L^2$ boundedness of
		the Riesz transforms give the maximal, slice and tangential trace bounds, first on
		$L^2\cap\dot H^{-1/2}$ and then by density on $L^2$.
		To compute the upper outward conormal derivative, put $v=\mS_{A_b}f$ with
		$f\in L^2\cap\dot H^{-1/2}$. For $\varphi\in C_c^\infty(\R^{n+1})$, integrate by parts in
		$\{t>\varepsilon\}$. Since $L_{A_b}v=0$ there, we obtain
		\[
		\int_{t>\varepsilon} A_b\nabla v\cdot\nabla\varphi
		=\int_{\R^n}\bigl(-\partial_t v+b\cdot\nabla_xv\bigr)(x,\varepsilon)
		\varphi(x,\varepsilon)\,dx.
		\]
		The derivatives in the preceding formulas converge in $L^2$ as $t \to 0^+$, and $b$ is bounded.
		Letting $\varepsilon \to 0^+$ gives the upper conormal derivative
		\[
		B_bf=\tfrac12\Bigl(f+\sum_{j=1}^n b_jR_jf\Bigr).
		\]
		This formula defines a bounded operator on $L^2$. We claim that $g$ is outside its range.
		Indeed, if $B_bf=g$ for some $f\in L^2$, then
		$\sum_{j=1}^n b_jR_jf$ is in $L^2$ and vanishes off a set of finite measure. By H\"older's
		inequality and \eqref{eq:negative},
		\[
		\sum_{j=1}^n b_jR_jf
		\in L^{2n/(n+1)}(\R^n)\subset\dot H^{-1/2}(\R^n).
		\]
		The identity $f=2g-\sum_{j=1}^n b_jR_jf$ now gives $f\in\dot H^{-1/2}$. Thus, the restriction of
		\eqref{sl:unchanged-formula} to the upper half-space is in $Y^{1,2}$ and has Neumann datum $g$ for
		$A^{(1)}$. By uniqueness in \eqref{eq:variational} it is the solution of Theorem~\ref{thm:main}. Its
		gradient satisfies the $L^2$ maximal bound by the Poisson estimate, contradicting that theorem.
		
		As observed in \cite[p.~473]{AM14}, the reflected matrix can depend on $t$
		even when the coefficients are independent of $t$ in the upper half-space. Here
		$\mA=I+\operatorname{sgn}(t)S(b)$, so the mixed entries change sign across the plane. Now Lemma~\ref{sl:reflection} gives $\mS_{\mA}(2g)=u_g^{\sharp}$,
		and the layer estimates fail for $2g$ by Theorem~\ref{sl:negative}
		(see Figure~\ref{lp:two-extensions}).
	\end{remark}
	
	\FloatBarrier
	
	\subsection{A small DKP norm}\label{sl:small-section}
	
	Let $n\geq2$, and let $A$ be a real, locally Lipschitz matrix in $\R^{n+1}_+$. We write $\lambda>0$ and
	$\Lambda<\infty$ for the constants in
	\begin{equation}\label{sl:ellipticity}
		A(x,t)\xi\cdot\xi\geq\lambda|\xi|^2,\qquad
		|A(x,t)\xi|\leq \Lambda|\xi|
		\quad\text{for a.e. }(x,t)\in\R^{n+1}_+\text{ and every }\xi\in\R^{n+1}.
	\end{equation}
	We denote by $u_g$ the solution of \eqref{eq:variational} for $g\in\dot H^{-1/2}(\R^n)$,
	and use the measure $\mu_A$ defined in \eqref{eq:dkp-definition}. By the trace estimate and the Lax--Milgram theorem,
	\[
	\|\nabla u_g\|_2\leq C_{n,\lambda}\|g\|_{\dot H^{-1/2}}.
	\]
	Unless indicated otherwise, constants depend only on $n$, $\lambda$, and $\Lambda$, and in the maximal
	estimates also on the aperture.
	
	For the Neumann estimate we use \cite[Theorem~7.1]{DPR17} and
	\cite[Theorem~2.2]{DP24}.
	
	\begin{theorem}[Neumann estimate for small DKP norm]\label{sl:small-neumann-input}
		Let $n\geq2$, $\lambda>0$, and $\Lambda<\infty$. There exists $\ve_0=\ve_0(n,\lambda,\Lambda)>0$ with
		the following property. Assume that $A$ is real and locally Lipschitz in $\R^{n+1}_+$, that
		\eqref{sl:ellipticity} holds, and that $\|\mu_A\|_{\mC}<\ve_0$. For
		every $g\in L^2(\R^n)\cap\dot H^{-1/2}(\R^n)$, the weak Neumann solution $u_g\in Y^{1,2}(\R^{n+1}_+)$ then
		satisfies
		\begin{equation}\label{sl:small-neumann}
			\|N_a(\nabla u_g)\|_2\leq C_{n,\lambda,\Lambda,a}\|g\|_2
			\qquad(a>1).
		\end{equation}
		The threshold $\ve_0$ is independent of the aperture $a$.
	\end{theorem}
	\begin{proof}
		We take $p=2$ and $r_0=\infty$ in \cite[Theorem~2.2]{DP24}. The boundary Lipschitz constant is
		zero. The datum in \cite[Definition~1.7]{DP24} belongs to $L^2\cap\dot H^{-1/2}$, and the weak
		identity there is \eqref{eq:variational} for the variational solution $u_g$.
		For the oscillation of the coefficients we have
		\[
		\operatorname{osc}_{B((x,t),t/2)}A
		\leq t\esssup_{B((x,t),t/2)}|\nabla A|,
		\]
		so the oscillation measure in that theorem has Carleson norm at most $C_n\|\mu_A\|_{\mC}$.
		For $\ve_0$ sufficiently small, the cited theorem gives
		$\|\widetilde N_{a_0}(\nabla u_g)\|_2\lesssim\|g\|_2$ at a fixed aperture $a_0>1$.
		By change of aperture and \eqref{eq:gradient-maximal-comparison}, we obtain
		\eqref{sl:small-neumann} for every $a>1$ with the same choice of $\ve_0$.
	\end{proof}
	
	\begin{proof}[Proof of Theorem~\ref{sl:small-layer}]
		Let $f\in L^2\cap\dot H^{-1/2}$ and $u=\mS_{\mA}f$. By Lemma~\ref{sl:reflection}, $u$ is the
		even extension of the Neumann solution with datum $f/2$. Hence \eqref{sl:small-neumann} gives
		\eqref{sl:small-maximal} on both sides of the plane. The estimate for $\widetilde N_a^\pm$
		follows from \eqref{eq:gradient-maximal-comparison}.
		Fix $a=2$. For almost every $(x,t)\in\R^{n+1}_+$, we have
		$|\nabla u(x,t)|\leq N_a(\nabla u)(x)$, since $(x,t)$ is a Lebesgue point of $\nabla u$ in the
		interior of $\gamma_a(x)$. Thus Fubini's theorem and \eqref{sl:small-maximal} give
		\[
		\|\nabla u(\cdot,t)\|_2\lesssim\|f\|_2
		\qquad\text{for almost every }t>0.
		\]
		Given $t_0>0$, we choose a sequence of such heights $t_j\to t_0$. For each bounded horizontal set $K$, the
		continuity of the slices in $L^2(K)$ from Subsection~\ref{sl:section} gives
		\[
		\|\nabla u(\cdot,t_0)\|_{L^2(K)}
		=\lim_j\|\nabla u(\cdot,t_j)\|_{L^2(K)}
		\lesssim\|f\|_2.
		\]
		Letting $K$ increase to $\R^n$, we obtain the bound at $t_0$. Reflection gives it for every negative
		height as well.
		
		By integrating $\partial_tu$ on the vertical lines above $K$ and applying
		Cauchy--Schwarz, we have
		\[
		\|u(\cdot,t)-\Tr u\|_{L^2(K)}^2
		\leq t\int_{K\times(0,t)}|\partial_tu|^2\longrightarrow0.
		\]
		For $\varphi\in C_c^\infty(\R^n;\R^n)$, we take $K$ containing $\supp\varphi$. Integrating by parts in $x$ and
		using the slice bound gives
		\begin{align*}
			|\langle\nabla_x\Tr u,\varphi\rangle|
			&=\lim_{t \to 0^+}
			\left|\int_{\R^n}\nabla_xu(x,t)\cdot\varphi(x)\,dx\right|\\
			&\lesssim\|f\|_2\|\varphi\|_2.
		\end{align*}
		Hence, $\nabla_x\Tr u\in L^2(\R^n)$ and $\|\nabla_x\Tr u\|_2\leq C\|f\|_2$, as asserted in
		\eqref{sl:small-slice-trace}. Moreover, Lemma~\ref{sl:reflection} identifies both outward conormal
		derivatives with $f/2$, so each has $L^2$ norm $\|f\|_2/2$.
		
		Given $f\in L^2(\R^n)$, choose $f_j\in L^2(\R^n)\cap\dot H^{-1/2}(\R^n)$ with $f_j\to f$ in $L^2(\R^n)$.
		This is possible since, by \eqref{eq:negative}, the intersection contains $C_c^\infty(\R^n)$, which is dense in $L^2(\R^n)$.
		Set $u_j=\mS_{\mA}f_j-\avint_{B(0,1)}\mS_{\mA}f_j$, where $B(0,1)\subset\R^{n+1}$.
		Integrating the slice bound for $u_j-u_k$ over $-R<t<R$, we get
		\begin{equation}\label{sl:completion-slab}
			\int_{\R^n\times(-R,R)}|\nabla(u_j-u_k)|^2\,dx\,dt
			\leq CR\|f_j-f_k\|_2^2
			\qquad(R>0).
		\end{equation}
		For $R>1$, put $v=u_j-u_k$ and $m_R=\avint_{B(0,R)}v$. Since $v$ has mean zero on
		$B(0,1)$, we have $|m_R|\leq\avint_{B(0,1)}|v-m_R|$. Thus
		$\|v\|_{L^{2^*}(B(0,R))}\leq C_R\|v-m_R\|_{L^{2^*}(B(0,R))}$.
		The Sobolev--Poincar\'e inequality and \eqref{sl:completion-slab} now give
		\[
		\|u_j-u_k\|_{Y^{1,2}(B(0,R))}
		\leq C_R\|f_j-f_k\|_2\longrightarrow0
		\qquad(R>1).
		\]
		Thus the normalized functions $u_j$ converge in $Y^{1,2}$ on each $B(0,R)$ to a function $u$.
		Applying the same estimate to two choices of $f_j$ gives the same normalized limit. Changing
		the normalization adds a constant.
		For $\varphi\in C_c^\infty(\overline{\R^{n+1}_\pm})$, local $L^2$ convergence of the gradients and
		convergence of the densities in $L^2(\R^n)$ give
		\begin{align*}
			\int_{\R^{n+1}_\pm}\mA\nabla u\cdot\nabla\varphi\,dx\,dt
			&=\lim_j\int_{\R^{n+1}_\pm}\mA\nabla u_j\cdot\nabla\varphi\,dx\,dt\\
			&=\frac12\lim_j\int_{\R^n}f_j(x)\varphi(x,0)\,dx
			=\frac12\int_{\R^n}f(x)\varphi(x,0)\,dx.
		\end{align*}
		This gives the outward conormal datum $f/2$ on each side. For $\varphi\in C_c^\infty(\R^{n+1})$, adding the two identities gives
		\eqref{sl:variational} with density $f$ for this test function.
		
		Along a subsequence, $\nabla u_j\to\nabla u$ almost everywhere. Taking essential suprema in each cone and
		then applying Fatou's lemma on the boundary gives
		\[
		N_a^\pm(\nabla u)\leq\liminf_jN_a^\pm(\nabla u_j),
		\qquad
		\|N_a^\pm(\nabla u)\|_2
		\leq\liminf_j\|N_a^\pm(\nabla u_j)\|_2
		\lesssim\|f\|_2.
		\]
		For fixed $j$, applying the same argument to the sequence $u_k-u_j$ as $k\to\infty$ gives
		\[
		\|N_a^\pm(\nabla(u-u_j))\|_2
		\lesssim\|f-f_j\|_2\longrightarrow0.
		\]
		The averaged bound follows again from \eqref{eq:gradient-maximal-comparison}. For the slices, the maximal estimate gives
		$\|\nabla u(\cdot,t)\|_2\lesssim\|f\|_2$ for almost every $t\ne0$ by Fubini's theorem.
		For any $t_0\ne0$, take such heights $t_k\to t_0$ on the same side of the plane.
		By Subsection~\ref{sl:section}, $\nabla u(\cdot,t_k)\to\nabla u(\cdot,t_0)$ in $L^2(B(0,R))$
		for every $R>0$. Passing first to the limit in $k$ and then letting $R\to\infty$ gives the
		bound at $t_0$. By \eqref{sl:small-slice-trace} applied to $u_j-u_k$,
		the sequence $\nabla_x\Tr u_j$ is Cauchy in $L^2(\R^n;\R^n)$.
		Local convergence in $Y^{1,2}$ gives $\Tr u_j\to\Tr u$ in $L^2_{\mathrm{loc}}(\R^n)$.
		Thus, for $\varphi\in C_c^\infty(\R^n;\R^n)$, integration by parts gives
		\[
		\lim_j\int_{\R^n}\nabla_x\Tr u_j\cdot\varphi
		=-\lim_j\int_{\R^n}\Tr u_j\,\div\varphi
		=-\int_{\R^n}\Tr u\,\div\varphi.
		\]
		It follows that $\nabla_x\Tr u_j\to\nabla_x\Tr u$ in $L^2(\R^n;\R^n)$ and
		\[
		\|\nabla_x\Tr u\|_2=\lim_j\|\nabla_x\Tr u_j\|_2
		\leq C\lim_j\|f_j\|_2=C\|f\|_2.
		\]
		
		If $f\in L^2(\R^n)\cap\dot H^{-1/2}(\R^n)$, we may take $f_j=f$ for every $j$. Since the limit is independent of the
		approximating sequence, $\nabla u=\nabla\mS_{\mA}f$.
	\end{proof}
	
	\begin{remark}
		In \eqref{sl:small-slice-trace}, $\nabla_x\Tr\mS_{\mA}f$ is the distributional tangential
		derivative of the trace. The outward conormal derivatives on the two sides are given by the
		weak identities above, which only involve the coefficients in the open half-spaces.
	\end{remark}

\end{document}